\documentclass[10pt]{amsart}
\usepackage{amssymb,latexsym,amsmath}
\usepackage[english]{babel}
\usepackage{latexsym}
\usepackage[all]{xy}
\date{}
\font\mm=msbm10
\def\bbb#1{{\hbox{{\mm #1}}}}
\def\h{\hfill}
\DeclareSymbolFont{largesymbolsA}{U}{pxexa}{m}{n}
\DeclareMathSymbol{\bigtimes}{\mathop}{largesymbolsA}{16}
\begin{document}
\title[Symmetric polynomials and exterior power]
{Symmetric polynomials and exterior power\\
of a polynomial ring in one variable}
\author{T. R. Seifullin}
\email{timur\_sf@mail.ru}
\address{V. M. Glushkov Institute of Cybernetics\\
National Academy of Sciences of Ukraine} 
\maketitle

\setcounter{MaxMatrixCols}{20}

\begin{abstract}
\noindent  
We consider exterior power of a polynomial ring in one variable 
as a module over a ring of symmetric polynomials. It is obtained 
expliáit expressions of symmetric polynomials via elementary symmetric 
polynomials.
\end{abstract}

\section*{Introduction}

\noindent
In this article we consider the $r$-th exterior power and 
the $r$-th symmetric tensors of the polynomal ring in one variable. We consider   
the second as an algebra,  and the first as a module over second. 
The $r$-th symmetric tensors of the polynomal ring in one variable as an algebra is 
isomorphic to the algebra of symmetric polynomials in $r$ vatiables. 
By this consideration it was obtained explicit expression of any symmetric 
polynomial via elementary symmetric polynomials in 1 of theorem 2.2. 
A similar consideration is available in the works of authors 
D. Laksov and A. Thorup. 
In the work [6] of these authors by using similar consideration it was 
obtained formula in (2) of theorem 0.1 from which can be obtained another 
formula that explicitly expresses any symmetric polynomial via elementary 
symmetric polynomials. 
In writing this article I relied on works of N. Bourbaki [1-5].  

\section{Preliminaries}
%

 \noindent 
{\bf Notation 1.1.} Denote by \medskip

 $(a_i|_{i{=}p,p{-}1} ) =  ()$, $(a_i|_{i{=}p,p} ) =  (a_p)$, 
\smallskip

 $(a_i|_{i{=}p,r} ) =  (a_i|_{i{=}p,q} ,a_i|_{i{=}q{+}1,r} 
)$, if $p{\leq}q$, $q{+}1{\leq}r$. \medskip\\
Denote by \medskip

 $(a_{{\leq}r} ) =  (a_i|_{i{=}1,r} )$. \medskip
\\
{\bf Notation 1.2.} Denote by \medskip

 $\| a_i|_{i{=}p,p{-}1} \|  =  \| \| $, $\| a_i|_{i{=}p,p} 
\|  =  \| a_p\| $, \smallskip

 $\| a_i|_{i{=}p,q} \|  =  \| \begin{matrix}a_i|_{i{=}p,r} 
 & a_i|_{i{=}r{+}1,q} \end{matrix} \| $, if $p{\leq}r$, $r{+}1{\leq}q$. 
\medskip
\\
We will   say that $\| a_i|_{i{=}p,q} \| $ is a vector 
of the size $|_{q{-}p{+}1} $. \medskip\\
Denote by \medskip

 $\| a_{{\leq}r} \|  =  \| a_i|_{i{=}1,r} \| $. \medskip
\\
{\bf Notation 1.3.} Denote by \medskip

 $\| a^j|^{j{=}k,k{-}1} \|  =  \| \| $, $\| a^j|^{j{=}k,k} 
\|  =  \| a^k\| $, \medskip

 $\| a^j|^{j{=}k,l} \|  =  \left\| \begin{matrix}a_i|^{j{=}k,m} 
  \h \cr
a_i|^{j{=}m{+}1,l} \h\end{matrix} \right\| $, if $k{\leq}m$, 
$m{+}1{\leq}l$. \medskip\\
We will   say that $\| a^j|^{j{=}k,l} \| $ is a covector 
of the size $|^{l{-}k{+}1} $. \medskip\\
Denote by \medskip

 $\| a^{{\leq}m} \|  =  \| a^i|^{i{=}1,m} \| $. \medskip\\
{\it There holds \medskip

 $\| a^j_i|_{i{=}p,q} |^{j{=}k,l} \|  =  \| a^j_i|^{j{=}k,l} 
|_{i{=}p,q} \| $, if $p{\leq}q{+}1$, $k{\leq}l{+}1$. } \medskip
\\
{\bf Notation 1.4.} Denote by \medskip

 $\| a^j_i|^{j{=}k,l} _{i{=}p,q} \|  =  \| a^j_i|_{i{=}p,q} 
|^{j{=}k,l} \|  =  \| a^j_i|^{j{=}k,l} |_{i{=}p,q} \| $, if 
$p{\leq}q{+}1$, $k{\leq}l{+}1$. \medskip\\
We will   say that $\| a^j_i|^{j{=}k,l} _{i{=}p,q} \| $ 
is a matrix of the size $|^{l{-}k{+}1} _{q{-}p{+}1} $. \bigskip
\hfil\vfill\eject

\noindent 
{\bf Notation 1.5.} Denote by \medskip

 $(a|_{{\times}0} ) =  ()$, $(a|_{{\times}1} ) =  (a)$, \smallskip

 $(a|_{{\times}(p{+}q)} ) =  (a|_{{\times}p} ,a|_{{\times}q} 
)$, if $p{\geq}1$, $q{\geq}1$. \medskip
\\
{\bf Notation 1.6.} Denote by \medskip

 $\| a|_{{\times}0} \|  =  \| \| $, $\| a|_{{\times}1} \|  
=  \| a\| $, \smallskip

 $\| a|_{{\times}(p{+}q)} \|  =  \| \begin{matrix}a|_{{\times}p} 
 & a|_{{\times}q} \end{matrix} \| $, if $p{\geq}1$, $q{\geq}1$. 
\medskip
\\
{\bf Notation 1.7.} Denote by \medskip

 $\| a|^{{\times}0} \|  =  \| \| $, $\| a|^{{\times}1} \|  
=  \| a\| $, \medskip

 $\| a|^{{\times}(k{+}l)} \|  =  \left\| \begin{matrix}a|^{{\times}k} 
\h \cr
a|^{{\times}l} \h\end{matrix} \right\| $, if $k{\geq}1$, 
$l{\geq}1$. \medskip\\
{\it There holds \medskip

 $\| a|_{{\times}k} |^{{\times}p} \|  =  \| a|^{{\times}p} 
|_{{\times}k} \| $, if $k{\geq}0$, $p{\geq}0$.} \bigskip
\\
{\bf Notation 1.8.} Denote by \medskip

 $\| a|^{{\times}p} _{{\times}k} \|  =  \| a|_{{\times}k} |^{{\times}p} 
\|  =  \| a|^{{\times}p} |_{{\times}k} \| $, if $k{\geq}0$, $p{\geq}0$.\bigskip 
\\ 
{\bf Notation 1.9.} Let $\top $  be an associative 
operation, denote by \medskip

 $\top (a_i|_{i{=}p,p} ) =  \top (a_p) =  a_p$, \smallskip

 $\top (a_i|_{i{=}p,r} ) =  \top (\top (a_i|_{i{=}p,q} ),\top 
(a_i|_{i{=}q{+}1,r} ))$, if $p{\leq}q$, $q{+}1{\leq}r$. \medskip
\\
Let $\top $  be an associative operation with an unity ${\bf 1}$, 
denote by \medskip

 $\top (a_i|_{i{=}p,p{-}1} ) =  \top () =  {\bf 1}$. \medskip
\\
{\bf Notation 1.10.} If $\top $   be an associative  and 
 commutative   operation,  let $(a_i|_{i{\in}I} )$  be a family of
elements, where $I$ be a finite set. \medskip
\\
Denote by \medskip

 $\top (a_i|_{i{\in}I} ) =  \top (a_p) =  a_p$, if $I{=}\{p\}$, 
\smallskip

 $\top (a_i|_{i{\in}I'{\cup}I''} ) =  \top (\top (a_i|_{i{\in}I'} 
),\top (a_i|_{i{\in}I''} ))$, if $I'{\cap}I''{=}\emptyset $. 
\medskip
\\
If $\top $  be an associative and commutative   operation with 
an unity ${\bf 1}$, denote by \medskip

 $\top (a_i|_{i{\in}I} ) =  \top () =  {\bf 1}$, if $I{=}\emptyset 
$. \medskip
\\
{\bf Notation 1.11.} Let $\top $  be an binary operation. 
Let  $l{=}\| l^j|^{j{=}k,l} \| $.  Denote by \medskip

 $\| \begin{matrix}l & \top | & a_i|_{i{=}p,q} \end{matrix} 
\|  =  \| l^j\top a_i|^{j{=}k,l} _{i{=}p,q} \| $, if $p{\leq}q{+}1$, 
$k{\leq}l{+}1$. \medskip
\\
{\bf Notation 1.12.}  Let  $\top $   be an binary  operation. 
 Let  $a{=}\| a_i|_{i{=}p,q} \| $,  where $p{+}1{\leq}q$. Denote by 
\medskip

 $\| a\| \top b =  \| a_i\top b|_{i{=}p,q} \| $.\bigskip
\\
{\bf Notation  1.13.}  Let  $p,q{\in}{\bbb 
Z}$,  and  $q{\leq}p{-}1$.  Denote by   $[q,p]{=}\{l|l{\in}{\bbb 
Z},q{\leq}l{\leq}p\}$. \bigskip 
\\
{\bf Notation 1.14.} Let $P$  be a statement. We shall write 
$?P{=}1$,  if  $P$  is true; $?P{=}0$, if $P$ is 
false.\bigskip 
\\
{\bf Notation 1.15.} Let $A$  be a set. 
Denote by ${\bf 1}_A$ the identity map $A \rightarrow  A$, $a \mapsto  a$. \medskip
\\
{\bf Notation 1.16.} Let $A$, $B$  be sets,  $B{\subseteq}A$. 
Denote by  ${\bf 1}_{B{\subseteq}A} $ the  embedding map 
$B \rightarrow  A$, $b \mapsto  b$. 
\hfil\vfill\eject

\noindent
{\bf Notation 1.17.} Let ${\bf A}$  be an associative ring 
with an unity. Denote by $0_{\bf A}$ a zero element of the ring ${\bf A}$, 
denote by $1_{\bf A}$ an unity element of the ring ${\bf A}$. \medskip
\\
{\bf Notation  1.18.}  Let  ${\bf A}$   be an associative 
 ring,  ${\bf M}$   be a module  over  ${\bf A}$, denote by $0_{\bf M}$ 
a zero element of the module ${\bf M}$. \medskip
\\
Let $M$  be  a subset of the  set ${\bf M}$,  
denote by $\mathop{\rm Span}(M)$  the set of  
all linear combinations of elements of the set $M$,  
i.  e.  the set of all  elements 
of the form $\sum ({\bf a}_q{\cdot}{\bf m}_q|q{\in}Q)$, 
where $Q$  be a finite set, $({\bf a}_q|q{\in}Q)$  
be elements of the ring ${\bf A}$, $({\bf m}_q|q{\in}Q)$  
be elements of the set $M$, and call it a linear closure 
of the set $M$ over the ring ${\bf A}$.\bigskip 
\\
{\bf Notation 1.19.} Denote by $!r$ the set of all invertible 
maps  $[1,r]  \rightarrow $ $[1,r]$. \medskip
\\
Denote by ${\varepsilon} _r$ the identity map $[1,r] 
\rightarrow  [1,r]$, $i \mapsto  i$. \medskip
\\
A transposition we call element ${\tau} {\in}!r$ such  that for 
some $i',i''{\in}[1,r]$ such that $i'\not=i"$ there holds ${\tau} (i'){=}i''$, 
${\tau} (i''){=}i'$, for any $i{\in}[1,r]$  such  that  
$i{\not=}i'$,  $i{\not=}i''$  there holds ${\tau} (i){=}i$. \medskip
\\
{\it  $!r$  is  a group  under  operations  compositions 
  maps, an unit element of this group is ${\varepsilon} _r$. 
The group $!r$ generated by transpositions.} \medskip
\\
{\bf Definition 1.1.} {\it There exists an unique homomorphism  of groups 
group $!r$  into  subgroup  $\{{-}1,1\}$ of multiplicative group  
integers ${\bbb Z}$ such  that ${\tau} \mapsto {-}1$
for any transposition ${\tau} {\in}!r$. } \medskip 
\\ 
Denote this homomorphism by $\mathop{\rm sgn}$. \bigskip
\\
Let ${\bf R}$  be a commutative ring with a unity. 
\medskip
\\
{\bf Convention 1.1.} Throughout this article, unless otherwise mentioned, 
by module  we  shall mean module over ${\bf R}$, by tensor 
 product  we   shall  mean tensor product over ${\bf R}$, 
by linear map we shall mean linear  map  over  ${\bf R}$,  
by  bilinear  map  we  shall  mean bilinear map over ${\bf R}$, 
by algebra we shall mean algebra over ${\bf R}$,  
by linear closure we  shall mean linear closure over ${\bf R}$, 
by homomorphism of algebras shall mean homomorphism of algebras over ${\bf R}$.
\medskip
\\
Let ${\bf M}$ be a module. \medskip\\
{\bf Notation 1.20.} Denote by $\mathop{\sf T}^r({\bf M}) 
=  \bigotimes ({\bf M}|_{i{=}1,r} )$. \medskip\\
Denote by the map \medskip

 $\mathop{\sf T}^r_{\bf M}{:}\ \bigtimes ({\bf M}|_{i{=}1,r} 
) \rightarrow  \bigotimes ({\bf M}|_{i{=}1,r} )$, $\bigtimes 
({\bf m}_i|_{i{=}1,r} ) \mapsto  \bigotimes ({\bf m}_i|_{i{=}1,r} )$. \medskip
\\
{\bf Definition 1.2.} Let ${\tau} {\in}!r$. Define the action of 
${\tau} $ on  $\mathop{\sf T}^r({\bf M})$  by \medskip

 ${\tau} _{\bf M}{:}\ \mathop{\sf T}^r({\bf M}) \rightarrow 
 \mathop{\sf T}^r({\bf M})$, $\bigotimes (\mathfrak{ m} _i|_{i{=}1,r} 
) \mapsto  \bigotimes (\mathfrak{ m} _{{\tau} ^{{-}1} (i)} |_{i{=}1,r} 
)$. \medskip\\
{\it 1) For any ${\tau} ',{\tau} ''{\in}!r$, $\mathfrak{ 
m} {\in}\mathop{\sf T}^r({\bf M})$ there holds \medskip

 ${\tau} '\, ({\tau} ''\, \mathfrak{ m} ) =  ({\tau} '\, {\tau} 
'')\, \mathfrak{ m} $, \medskip\\
2) For  $\mathfrak{ m} {\in}\mathop{\sf T}^r({\bf M})$ 
there holds \medskip

 ${\varepsilon} _r\, \mathfrak{ m}  =  \mathfrak{ m} $. }
\hfil\vfill\eject
 
 \noindent {\bf Notation 1.21.} Denote by $\mathop{\sf TS}^r({\bf M})$  
the set of all elements $\mathfrak{ m} {\in}\mathop{\sf 
T}^r({\bf M})$ such  that for any transposition ${\tau} {\in}!r$ 
there holds \medskip

{\it ${\tau} \, \mathfrak{ m}  =  \mathfrak{ m} $. \medskip\\
$\mathop{\sf TS}^r({\bf M})$ is a submodule of 
the module $\mathop{\sf T}^r({\bf M})$.} Denote by the map \medskip

$\mathop{\sf TS}^r_{\bf M}{:}\ \mathop{\sf TS}^r({\bf M}) 
\rightarrow  \mathop{\sf T}^r({\bf M})$, $\mathfrak{ m}  \mapsto 
 \mathfrak{ m} $. \medskip
\\
{\it Note that in cases $r{=}0$ and $r{=}1$ there holds $\mathop{\sf 
TS}^r({\bf M}) =  \mathop{\sf T}^r({\bf M})$. } \bigskip
\\
{\bf Notation 1.22.} Denote by $\mathop{\sf TA}^r({\bf M})$ 
the set of all elements $\mathfrak{ m} {\in}\mathop{\sf T}^r({\bf M})$ 
such  that for any transposition ${\tau} {\in}!r$ 
there holds \medskip

 ${\tau} \, \mathfrak{ m}  =  {-}\mathfrak{ m} $. \medskip
\\
{\it $\mathop{\sf TA}^r({\bf M})$ is a submodule of the module 
$\mathop{\sf T}^r({\bf M})$.} Denote by the map \medskip

 $\mathop{\sf TA}^r_{\bf M}{:}\ \mathop{\sf TA}^r({\bf M}) 
\rightarrow  \mathop{\sf T}^r({\bf M})$, $\mathfrak{ m}  \mapsto 
 \mathfrak{ m} $. \medskip\\
{\it Note that in cases $r{=}0$ and $r{=}1$ there holds $\mathop{\sf 
TA}^r({\bf M}) =  \mathop{\sf T}^r({\bf M})$. } \bigskip
\\
{\bf Notation 1.23.}  Denote by  $\mathop{\sf AT}^r({\bf M})$  
the submodule of the module $\mathop{\sf T}^r({\bf M})$  generated by 
elements $\mathfrak{ m} {\in}\mathop{\sf T}^r({\bf M})$ such 
 that for some transposition  ${\tau} {\in}!r$  there holds 
\medskip

 ${\tau} \, \mathfrak{ m}  =  {-}\mathfrak{ m} $, \medskip
\\
{\it Note that in cases $r{=}0$ and $r{=}1$ there holds $\mathop{\sf 
AT}^r({\bf M}) =  \{0_{\mathop{\sf T}^r({\bf M})} \}$. } \bigskip
\\
{\bf Notation 1.24.} Denote by $\mathop{\sf S}^r({\bf M}) 
 =   \mathop{\sf T}^r({\bf M})/\mathop{\sf AT}^r({\bf M})$.  
Denote by the map  \medskip

 $\mathop{\sf S}^r_{\bf M}{:}\ \mathop{\sf T}^r({\bf M}) \rightarrow 
 \mathop{\sf S}^r({\bf M})$, $\mathfrak{ m}  \mapsto  \mathfrak{ 
m} {+}\mathop{\sf AT}^r({\bf M})$. \medskip\\
{\it Note that in cases $r{=}0$ and $r{=}1$ there holds $\mathop{\sf 
S}^r({\bf M}) =  \mathop{\sf T}^r({\bf M})$.} \bigskip
\\
{\bf Notation  1.24.}  Denote by  $\mathop{\sf ST}^r({\bf M})$  
the submodule  of the module  $\mathop{\sf T}^r({\bf M})$  generated by 
elements $\mathfrak{ m} {\in}\mathop{\sf T}^r({\bf M})$ such 
 that for some transposition  ${\tau} {\in}!r$  there holds 
\medskip

 ${\tau} \, \mathfrak{ m}  =  \mathfrak{ m} $. \medskip
\\
{\it Note that in cases $r{=}0$ and $r{=}1$  there holds 
 $\mathop{\sf ST}^r({\bf M})  =   \{0_{\mathop{\sf T}^r({\bf M})} \}$.} \bigskip
\\
{\bf Notation 1.26.} Denote by $\mathop{\sf A}^r({\bf M}) 
=  \mathop{\sf T}^r({\bf M})/\mathop{\sf ST}^r({\bf M})$,  or 
 $\bigwedge  ^r({\bf M})$.  Denote by the map \medskip

 $\mathop{\sf A}^r_{\bf M}{:}\ \mathop{\sf T}^r({\bf M}) \rightarrow 
 \mathop{\sf A}^r({\bf M})$, $\mathfrak{ m}  \mapsto  
\mathfrak{ m} {+}\mathop{\sf ST}^r({\bf M})$. \medskip
\\
{\it Note that in cases $r{=}0$  and  $r{=}1$  there holds 
 $\mathop{\sf A}^r({\bf M})  =   \mathop{\sf T}^r({\bf M})$.} \medskip
%

 \noindent {\bf Notation 1.27.} Let ${\bf m}_i{\in}{\bf 
M}$ for all $i{=}1,r$, define \medskip

 $\bigwedge  ({\bf m}_i|_{i{=}1,r} ) =  \mathop{\sf A}^r_{\bf 
M}(\bigotimes ({\bf m}_i|_{i{=}1,r} )) =  \mathop{\sf A}^r_{\bf 
M}(\mathop{\sf T}^r_{\bf M}({\bf m}_i|_{i{=}1,r} ))$, \medskip

 $\prod ({\bf m}_i|_{i{=}1,r} ) =  \mathop{\sf S}^r_{\bf M}(\bigotimes 
({\bf m}_i|_{i{=}1,r} )) =  \mathop{\sf S}^r_{\bf M}(\mathop{\sf 
T}^r_{\bf M}({\bf m}_i|_{i{=}1,r} ))$.

\hfil\vfill\eject
 
 \noindent {\bf Notation 1.28.} Let $\mathfrak{ m} {\in}{\bf 
M}$. Denote by $\mathop{\sf T}^r({\bf m}) =  {\bf m}^{{\otimes}r} 
$. \medskip
\\
{\bf Proposition 1.1.} {\it Let  $\mathfrak{ m} {\in}{\bf M}$.  
For  any  ${\bf m}{\in}{\bf M}$,  any  ${\tau} {\in}!p$ 
 there holds \medskip

 ${\tau} \, \mathop{\sf T}^r({\bf m}) =  \mathop{\sf T}^r({\bf 
m})$, \medskip\\
i. e. $\mathop{\sf T}^r({\bf m}){\in}\mathop{\sf TS}^r({\bf 
M})$. } \medskip
\\
{\bf Proposition 1.2.} {\it For any ${\bf m}{\in}{\bf M}$, 
$\mathfrak{ m} '{\in}\mathop{\sf T}^{r'} ({\bf M})$,  $\mathfrak{ 
m} ''{\in}\mathop{\sf T}^{r''} ({\bf M})$,  where  $r'{+}r''{+}2{=}r$, 
there holds \medskip

 $\mathfrak{ m} '{\otimes}{\bf m}{\otimes}{\bf m}{\otimes}\mathfrak{ m} '' 
\in  \mathop{\sf ST}^r({\bf M})$. \medskip
\\
If ${\bf M}$ is a free module, then $\mathop{\sf ST}^r({\bf M})$ 
lineary  generated by  elements of the form \medskip

$\mathfrak{ m} '{\otimes}{\bf m}{\otimes}{\bf m}{\otimes}\mathfrak{ m}$, \medskip
\\
where ${\bf m}{\in}{\bf M}$, $\mathfrak{ m} '{\in}\mathop{\sf T}^{r'} 
({\bf M})$, $\mathfrak{ m} ''{\in}\mathop{\sf T}^{r''} ({\bf 
M})$, $r'{+}r''{+}2{=}r$. } \medskip
\\
{\bf Proposition 1.3.} {\it For any ${\bf m}',{\bf m}''{\in}{\bf 
M}$, $\mathfrak{ m} '{\in}\mathop{\sf T}^{r'} ({\bf M})$,  $\mathfrak{ 
m} ''{\in}\mathop{\sf T}^{r''} ({\bf M})$,  where  $r'{+}r''{+}2{=}r$, 
there holds\medskip

 $\mathfrak{ m} '{\otimes}({\bf m}'{\otimes}{\bf m}''{-}{\bf 
m}''{\otimes}{\bf m}'){\otimes}\mathfrak{ m} '' \in  \mathop{\sf 
AT}^r({\bf M})$. \medskip\\
If ${\bf M}$ is a free module, then $\mathop{\sf AT}^r({\bf M})$ 
is lineary generated by  elements of the form \medskip

 $\mathfrak{ m} '{\otimes}({\bf m}'{\otimes}{\bf m}''{-}{\bf 
m}''{\otimes}{\bf m}'){\otimes}\mathfrak{ m} ''$, \medskip\\
where  ${\bf m}{\in}{\bf M}$,  $\mathfrak{ m} '{\in}\mathop{\sf 
T}^{r'} ({\bf M})$,  $\mathfrak{ m} ''{\in}\mathop{\sf T}^{r''} 
({\bf M})$,  $r'{+}r''{+}2{=}r$.  } \medskip
\\
Let ${\bf M}$  be a module. \medskip
\\
{\bf Definition 1.3.} Denote by 
the map \medskip

 $\mathop{\sf Atr}^r_{\bf M}{:}\ \mathop{\sf T}^r({\bf M}) 
\rightarrow  \mathop{\sf T}^r({\bf M})$, $\mathfrak{ m}  \mapsto 
 \sum (\mathop{\rm sgn}({\tau} ){\cdot}({\tau} \, \mathfrak{ 
m} )|{\tau} {\in}!p)$. \bigskip\\
{\bf Proposition-notation 1.4.} {\it The map $\mathop{\sf 
Atr}^r_{\bf M}$ is a linear. \medskip\\
There holds \medskip

 $\mathop{\sf Atr}^r_{\bf M}\, \mathop{\sf T}^r({\bf M}) \subseteq 
 \mathop{\sf TA}^r({\bf M})$, \medskip

 $\mathop{\sf Atr}^r_{\bf M}\, \mathop{\sf ST}^r({\bf M}) \subseteq 
 \{0_{\mathop{\sf TA}^r({\bf M})} \}$. } \medskip\\
Hence, {\it the linear map $\mathop{\sf Atr}^r_{\bf M}$ induces a linear map \medskip

 $\mathop{\sf A}^r({\bf M}) \rightarrow  \mathop{\sf TA}^r({\bf M})$, } \medskip
\\
which denote by $\mathop{\sf atr}^r_{\bf M}$. \medskip
\\
{\bf Proposition 1.5.} {\it If a module ${\bf M}$ is 
free, then the map  $\mathop{\sf atr}^r_{\bf M}$ is invertible.}
\hfil\vfill\eject
 
 \noindent Let ${\bf N}$, ${\bf L}$  be modules over ${\bf 
R}$ and ${\phi} {:}\ {\bf N} \rightarrow  {\bf L}$  be a linear map. \medskip
\\
{\bf Notation 1.29.} Denote by $\mathop{\sf T}^r({\phi} )$ the tensor 
power maps ${\phi} {:}\ {\bf N} \rightarrow 
 {\bf L}$, i. e. the linear map 
$\mathop{\sf T}^r({\bf N}) \rightarrow  \mathop{\sf T}^r({\bf L})$ 
such that \medskip

 $\mathop{\sf T}^r({\phi} )\, (\bigotimes ({\bf n}_i|_{i{=}1,r} )) =  
\bigotimes ({\phi} ({\bf n}_i)|_{i{=}1,r} )$, \medskip
\\
where ${\bf 
n}_i{\in}{\bf N}$ for $i{=}1,r$. \medskip
\\
{\bf Proposition 1.6.} {\it Let ${\bf N}$, ${\bf L}$  be modules over  
${\bf R}$  and  ${\phi} {:}\  {\bf N}  \rightarrow   {\bf L}$  be a linear map. 
Let ${\tau} {\in}!r$, $\mathfrak{ n} {\in}\mathop{\sf T}^r({\bf N})$, then 
\medskip

 ${\tau} \, \mathop{\sf T}^r({\phi} )\, (\mathfrak{ n} ) = 
 \mathop{\sf T}^r({\phi} )\, ({\tau} \, \mathfrak{ n} )$. } \medskip
\\
{\bf Proposition 1.7.} {\it Let ${\bf N}$, ${\bf L}$  be  
modules, ${\phi} {:}\ {\bf N} \rightarrow  {\bf L}$   be a linear map. Tnen \medskip

 1. $\mathop{\sf T}^r({\phi} )(\mathop{\sf TS}^r({\bf N})) 
\subseteq  \mathop{\sf TS}^r({\bf L}^r)$. \medskip

 2. $\mathop{\sf T}^r({\phi} )(\mathop{\sf TA}^r({\bf N})) 
\subseteq  \mathop{\sf TA}^r({\bf L}^r)$. \medskip

 3. $\mathop{\sf T}^r({\phi} )(\mathop{\sf ST}^r({\bf N})) 
\subseteq  \mathop{\sf ST}^r({\bf L}^r)$. }\\
\\
{\bf Proposition-notation 1.8.} {\it Let ${\bf N}$, ${\bf L}$ 
 be modules, ${\phi} {:}\ {\bf N} \rightarrow   {\bf L}$ be a linear map.\medskip
\\
Then the linear map $\mathop{\sf T}^r({\phi} ){:}\ \mathop{\sf T}^r({\bf N}) 
\rightarrow  \mathop{\sf T}^r({\bf L})$ induces linear  maps: \medskip

 1) $\mathop{\sf TS}^r({\bf N}) \rightarrow  \mathop{\sf TS}^r({\bf 
L})$, } which denote by $\mathop{\sf TS}^r({\phi} )$, {\it 
\medskip

 2) $\mathop{\sf TA}^r({\bf N}) \rightarrow  \mathop{\sf TA}^r({\bf 
L})$, } which denote by $\mathop{\sf TA}^r({\phi} )$, {\it 
\medskip

 3) $\mathop{\sf A}^r({\bf N}) \rightarrow  \mathop{\sf A}^r({\bf 
L})$. } which denote by $\mathop{\sf A}^r({\phi} )$.\bigskip 
\\
Let ${\bf A}$ an associative algebra over ${\bf R}$ with 
a product \medskip

 ${\pi} _{\bf A}{:}\ {\bf A}{\times}{\bf A} \rightarrow   {\bf A}$, \medskip
\\
and with a unity $1_{\bf A}$. \medskip
\\
{\bf Proposition 1.9.} {\it $\mathop{\sf T}^r({\bf A}) =  \bigotimes 
({\bf A}|_{i{=}1,r} )$ is an  associative  algebra, as 
tensor product of associative algebras with units.} \medskip
\\
{\bf Notation 1.30.} Denote by $\mathop{\sf T}^r({\pi} )$ 
the product in $\mathop{\sf T}^r({\bf A})$. \medskip\\
{\it Let ${\bf a}'_i,{\bf a}''_i{\in}{\bf A}$ for $i{=}1,r$. 
There holds \medskip

 $(\bigotimes ({\bf a}'_i|_{i{=}1,r} )){\cdot}(\bigotimes ({\bf a}''_i|_{i{=}1,r} )) =  \bigotimes ({\bf a}'_i{\cdot}{\bf a}''_i|_{i{=}1,r} )$.}\medskip 
\\ 
Denote by $1_{\mathop{\sf T}^r({\bf A})}$ the unity of 
$\mathop{\sf T}^r({\bf A})$, {\it there holds \medskip 

$1_{\mathop{\sf T}^r({\bf A})} = \mathop{\sf T}^r(1_{\bf A})$.} \medskip
\\
{\bf Proposition 1.10.} {\it If  ${\bf A}$  is a commutative  
 algebra,  then  $\mathop{\sf T}^r({\bf A})$ is a commutative  
algebra.}\medskip
\\
{\it For any ${\tau} {\in}!r$. $\mathfrak{ a} ',\mathfrak{ 
a} ''{\in}\mathop{\sf T}^r({\bf A})$ there holds \medskip

 ${\tau} \, (\mathfrak{ a} '{\cdot}\mathfrak{ a} '') =  ({\tau} 
\, \mathfrak{ a} '){\cdot}({\tau} \, \mathfrak{ a} '')$. } \medskip
\\
{\it Let ${\tau} {\in}!r$. Then \medskip

 ${\tau} \, 1_{\mathop{\sf T}^r({\bf A})}  =  1_{\mathop{\sf 
T}^r({\bf A})} $.}   \medskip\\
Hence, {\it \medskip

 $1_{\mathop{\sf T}^r({\bf A})}  \in  \mathop{\sf TS}^r({\bf 
A})$. } \medskip\\
\hfil\vfill\eject
 
 \noindent {\bf Proposition 1.11.} {\it There holds  \medskip

 1) $\mathop{\sf TS}^r({\bf A}){\cdot}\mathop{\sf TS}^r({\bf 
A}) \subseteq  \mathop{\sf TS}^r({\bf A})$; \medskip

 2) $\mathop{\sf TS}^r({\bf A}){\cdot}\mathop{\sf TA}^r({\bf 
A}) \subseteq  \mathop{\sf TA}^r({\bf A})$; \medskip

 3) $\mathop{\sf TS}^r({\bf A}){\cdot}\mathop{\sf ST}^r({\bf 
A}) \subseteq  \mathop{\sf ST}^r({\bf A})$. }\\
\\
{\bf Proposition 1.12.} {\it Let ${\bf A}$  be an associative 
algebra with an unity $1_{\bf A}$. \medskip\\
1. $\mathop{\sf TS}^r({\bf A})$ with the product  induces by 
the product  in  $\mathop{\sf T}^r({\bf A})$  is an associative 
algebra with the unity \medskip

 $1_{\mathop{\sf TS}^r({\bf A})}  =  1_{\mathop{\sf T}^r({\bf 
A})} $, \medskip\\
and is a subalgebra of the algebra $\mathop{\sf T}^r({\bf A})$. 
\medskip\\
If ${\bf A}$ is a commutative   algebra,  then  
$\mathop{\sf TS}^r({\bf A})$  is a commutative  algebra. \medskip\\
2. $\mathop{\sf TA}^r({\bf A})$ with the product  induces by 
the product  in  $\mathop{\sf T}^r({\bf A})$ is an unitary 
module over algebra $\mathop{\sf TS}^r({\bf A})$. \medskip\\
3. $\mathop{\sf A}^r({\bf A})$ with the product  induces by 
the product  in  $\mathop{\sf T}^r({\bf A})$  is an unitary 
module over algebra $\mathop{\sf TS}^r({\bf A})$.}\medskip 
\\
{\bf Proposition 1.13.} {\it The map $\mathop{\sf atr}^r_{\bf 
A}{:}\ \mathop{\sf A}^r({\bf A}) \rightarrow  \mathop{\sf TA}^r({\bf 
A})$ is a homomorphism of modules over $\mathop{\sf TS}^r({\bf 
A})$, i. e. for any $\mathfrak{ a} {\in}\mathop{\sf TS}^r({\bf 
A})$, $\mathfrak{ m} {\in}\mathop{\sf A}^r({\bf A})$ there holds 
\medskip

 $\mathop{\sf atr}^r_{\bf A}\, (\mathfrak{ a} {\cdot}\mathfrak{ 
m} ) =  \mathfrak{ a} {\cdot}(\mathop{\sf atr}^r_{\bf A}\, (\mathfrak{ 
m} ))$. }\medskip\\
{\bf Proof.} Let $\mathfrak{ a} {\in}\mathop{\sf 
TS}^r({\bf A})$, $\mathfrak{ m} {\in}\mathop{\sf A}^r({\bf A})$. 
Since  $\mathfrak{ a} {\in}\mathop{\sf TS}^r({\bf A})$,  then 
 for  any transposition ${\tau} {\in}!r$ there holds 
${\tau} \, \mathfrak{ a}  =  \mathfrak{ a} $. We have \medskip

 $\mathop{\sf atr}^r_{\bf A}\, (\mathfrak{ a} {\cdot}\mathfrak{ 
m} ) =  \sum (\mathop{\rm sgn}({\tau} ){\cdot}({\tau} \, (\mathfrak{ 
a} {\cdot}\mathfrak{ m} ))|{\tau} {\in}!r) =  \sum (\mathop{\rm 
sgn}({\tau} ){\cdot}(({\tau} \, \mathfrak{ a} ){\cdot}({\tau} 
\, \mathfrak{ m} ))|{\tau} {\in}!r)$  \medskip

 \qquad $=  \sum (\mathop{\rm sgn}({\tau} ){\cdot}(\mathfrak{ 
a} {\cdot}({\tau} \, \mathfrak{ m} ))|{\tau} {\in}!r) =  \mathfrak{ 
a} {\cdot}(\sum (\mathop{\rm sgn}({\tau} ){\cdot}({\tau} \, \mathfrak{ 
m} )|{\tau} {\in}!r)) =  \mathfrak{ a} {\cdot}(\mathop{\sf atr}^r_{\bf 
M}\, (\mathfrak{ m} ))$. \medskip\\
Hence, \medskip

 $\mathop{\sf atr}^r_{\bf A}\, (\mathfrak{ a} {\cdot}\mathfrak{ 
m} ) =  \mathfrak{ a} {\cdot}(\mathop{\sf atr}^r_{\bf A}\, (\mathfrak{ 
m} ))$. \medskip\\
Since $\mathop{\sf atr}^r_{\bf A}$ is a linear map, 
then $\mathop{\sf atr}^r_{\bf A}$  is a homomorphism of modules 
over $\mathop{\sf TS}^r({\bf A})$.\bigskip
\\ 
 \noindent Let ${\bf M}$  be a semigraded module.\medskip
\\
{\bf Notation 1.31.} Denote by $\mathop{\sf T}^r({\bf M})^{{\leq}d} 
 =  \mathop{\rm Span}(\mathop{\sf T}^r_{\bf M}\, ({\bf M}^{{\leq}d} 
|{\times}r))$. \medskip
\\
{\it $\mathop{\sf T}^r({\bf M})$  with a semigraded  
$(\mathop{\sf T}^r({\bf M})^{{\leq}d} |d{\in}{\bbb Z})$  
is a semigraded module. }\medskip 
\\
{\bf Notation 1.32.} Denote by $\mathop{\sf TS}^r({\bf M})^{{\leq}d} 
 =  \mathop{\sf T}^r({\bf M})^{{\leq}d} {\cap}\mathop{\sf TS}^r({\bf M})$. \medskip
\\
{\it $\mathop{\sf TS}^r({\bf M})$ with a semigraded  
$(\mathop{\sf TS}^r({\bf M})^{{\leq}d} |d{\in}{\bbb Z})$  
is a semigraded module. }\medskip 
\\
{\bf Notation 1.33.} Denote by $\mathop{\sf TA}^r({\bf M})^{{\leq}d} 
 =  \mathop{\sf T}^r({\bf M})^{{\leq}d} {\cap}\mathop{\sf TA}^r({\bf M})$. \medskip
\\
{\it $\mathop{\sf TA}^r({\bf M})$ with a semigraded  
$(\mathop{\sf TA}^r({\bf M})^{{\leq}d} |d{\in}{\bbb Z})$  
is a semigraded module. }\medskip 
\\
{\bf Notation 1.34.} Denote by $\mathop{\sf S}^r({\bf M})^{{\leq}d} 
 =  (\mathop{\sf T}^r({\bf M})^{{\leq}d} {+}\mathop{\sf AT}^r({\bf M}))/\mathop{\sf AT}^r({\bf M})$. \medskip
\\
{\it $\mathop{\sf S}^r({\bf M})$  with a semigraded  
$(\mathop{\sf S}^r({\bf M})^{{\leq}d} |d{\in}{\bbb Z})$  
is a semigraded module. }\medskip
\\
{\bf Notation 1.35.} Denote by $\mathop{\sf A}^r({\bf M})^{{\leq}d} 
 =  (\mathop{\sf T}^r({\bf M})^{{\leq}d} {+}\mathop{\sf ST}^r({\bf M}))
/\mathop{\sf ST}^r({\bf M})$. \medskip
\\
{\it $\mathop{\sf A}^r({\bf M})$  with a semigraded  
$(\mathop{\sf A}^r({\bf M})^{{\leq}d} |d{\in}{\bbb Z})$ 
is a semigraded module. }\medskip
\\
{\bf Proposition 1.14.} {\it There holds \medskip

 $\mathop{\sf S}^r({\bf M})^{{\leq}d}  =  \{\mathfrak{ m} '|\mathfrak{ 
m} '{=}\mathop{\sf S}^r_{\bf M}(\mathfrak{ m} ), \mathfrak{ m} 
{\in}\mathop{\sf T}^r({\bf M})^{{\leq}d} \}$, \medskip

 $\mathop{\sf A}^r({\bf M})^{{\leq}d}  =  \{\mathfrak{ m} '|\mathfrak{ 
m} '{=}\mathop{\sf A}^r_{\bf M}(\mathfrak{ m} ), \mathfrak{ m} 
{\in}\mathop{\sf T}^r({\bf M})^{{\leq}d} \}$, \medskip

 $\mathop{\sf TS}^r({\bf M})^{{\leq}d}  =  \{\mathfrak{ m} 
'|\mathop{\sf TS}^r_{\bf M}(\mathfrak{ m} '){=}\mathfrak{ m} 
, \mathfrak{ m} {\in}\mathop{\sf T}^r({\bf M})^{{\leq}d} \}$, 
\medskip 

 $\mathop{\sf TA}^r({\bf M})^{{\leq}d}  =  \{\mathfrak{ m} 
'|\mathop{\sf TA}^r_{\bf M}(\mathfrak{ m} '){=}\mathfrak{ m} 
, \mathfrak{ m} {\in}\mathop{\sf T}^r({\bf M})^{{\leq}d} \}$. 
}\medskip 
\\
{\bf Proposition  1.15.}  {\it  The map  
$\mathop{\sf atr}^r_{\bf M}{:}\  
\mathop{\sf A}^r({\bf M})  \rightarrow    \mathop{\sf AT}^r({\bf M})$
is a semihomogeneouse linear map of modules of the degree ${\leq}0$, i. e.  
there holds \medskip

 $\mathop{\sf atr}^r_{\bf M}(\mathop{\sf A}^r({\bf M})^{{\leq}d} 
) \subseteq  \mathop{\sf AT}^r({\bf M})^{{\leq}d} $. }\\
\hfil\vfill\eject
 
 \noindent {\bf Proposition 1.16.} {\it Let ${\bf m}_i{\in}{\bf 
M}^{{\leq}d} $ for all $i{=}1,r$. Then \medskip

 $\bigwedge  ({\bf m}_i|_{i{=}1,r} ) =  \mathop{\sf A}^r_{\bf 
M}(\bigotimes ({\bf m}_i|_{i{=}1,r} )) \in  \mathop{\sf A}^r({\bf 
M})^{{\leq}d} $, \medskip

 $\prod ({\bf m}_i|_{i{=}1,r} ) =  \mathop{\sf S}^r_{\bf M}(\bigotimes 
({\bf m}_i|_{i{=}1,r} )) \in  \mathop{\sf S}^r({\bf M})^{{\leq}d} 
$. }\medskip 
\\
{\bf Proposition 1.17.}\\
{\it 1) There holds \medskip

 $\mathop{\sf A}^r({\bf M})^{{\leq}d}  =  \mathop{\rm Span}(\mathop{\sf 
A}^r_{\bf M}\, ({\bf M}^{{\leq}d} |_{{\times}r} )) =  \mathop{\rm 
Span}(\bigwedge  \, ({\bf M}^{{\leq}d} |_{{\times}r} ))$, \medskip
\\
i. e. $\mathop{\sf A}^r({\bf M})^{{\leq}d} $ is lineary generated by 
elements of the form  $\bigwedge  ({\bf m}_i|_{i{=}1,r} )$,  where  
${\bf m}_i{\in}{\bf M}^{{\leq}d} $ for all $i{=}1,r$. \medskip
\\
2) There holds \medskip

 $\mathop{\sf S}^r({\bf M})^{{\leq}d}  =  \mathop{\rm Span}(\mathop{\sf 
S}^r_{\bf M}\, ({\bf M}^{{\leq}d} |_{{\times}r} )) =  \mathop{\rm 
Span}(\prod \, ({\bf M}^{{\leq}d} |_{{\times}r} ))$, \medskip
\\
i. e. $\mathop{\sf S}^r({\bf M})^{{\leq}d} $ is lineary generated by 
elements of the form $\prod ({\bf m}_i|_{i{=}1,r} )$,  where  ${\bf 
m}_i{\in}{\bf M}^{{\leq}d} $ for all $i{=}1,r$. }\bigskip 
\\
Let ${\bf A}$  be a semigraded algebra with 
an unity.\medskip
\\
{\bf Proposition 1.18.} {\it $\mathop{\sf T}^r({\bf A})$ is 
a semigraded algebra with an unity, then $\mathop{\sf T}^r({\bf A})$ 
is a semigraded module and there holds 
\medskip

 $\mathop{\sf T}^r({\bf A})^{{\leq}d'} {\cdot}\mathop{\sf T}^r({\bf 
A})^{{\leq}d''}  \subseteq  \mathop{\sf T}^r({\bf A})^{{\leq}d'{+}d''} 
$, \medskip

 $1_{\mathop{\sf T}^r({\bf A})}  \in  \mathop{\sf T}^r({\bf 
A})^{{\leq}0} $. }\medskip 
\\
{\bf Proposition  1.19.}  {\it  $\mathop{\sf TS}^r({\bf A})$ 
 is a  semigraded   algebra   with an unity, i. e. 
$\mathop{\sf TS}^r({\bf A})$ is  semigraded  module 
 and  there holds \medskip

 $\mathop{\sf TS}^r({\bf A})^{{\leq}d'} {\cdot}\mathop{\sf 
TS}^r({\bf A})^{{\leq}d''}  \subseteq  \mathop{\sf TS}^r({\bf 
A})^{{\leq}d'{+}d''} $, \medskip

 $1_{\mathop{\sf TS}^r({\bf A})}  \in  \mathop{\sf TS}^r({\bf A})^{{\leq}0} $. }\\
\\
{\bf Proposition  1.20.}  {\it  $\mathop{\sf TA}^r({\bf A})$ 
 is a semigraded  module   over the semigraded 
algebra $\mathop{\sf TS}^r({\bf A})$, i. e. $\mathop{\sf TA}^r({\bf 
A})$ is a semigraded module, algebra over $\mathop{\sf 
TS}^r({\bf A})$, and there holds \medskip

 $\mathop{\sf TS}^r({\bf A})^{{\leq}d'} {\cdot}\mathop{\sf 
TA}^r({\bf A})^{{\leq}d''}  \subseteq  \mathop{\sf TA}^r({\bf 
A})^{{\leq}d'{+}d''} $. }\medskip 
\\
{\bf Proposition 1.21.} {\it $\mathop{\sf A}^r({\bf A})$ is 
a semigraded  unitary  module over the  semigraded 
  algebra   $\mathop{\sf TS}^r({\bf A})$, i. e. 
$\mathop{\sf A}^r({\bf A})$ 
is a semigraded module, unitary 
module over the algebra  $\mathop{\sf TS}^r({\bf A})$,  and  there holds \medskip

 $\mathop{\sf TS}^r({\bf A})^{{\leq}d'} {\cdot}\mathop{\sf 
A}^r({\bf A})^{{\leq}d''}  \subseteq  \mathop{\sf A}^r({\bf A})^{{\leq}d'{+}d''} 
$. }\medskip 
\\
{\bf Proposition  1.22.}   {\it   Let   embedding map 
  ${\bf 1}_{{\bf M}^{{\leq}d} {\subseteq}{\bf M}} $   is 
left invertible, i. e. there exists a linear map $\mathfrak{ 
p} ^d_{\bf M}{:}\ {\bf M} \rightarrow  {\bf M}^{{\leq}d} $  such 
  that \medskip

 $\mathfrak{ p} ^d_{\bf M}\, {\bf 1}_{{\bf M}^{{\leq}d} {\subseteq}{\bf 
M}}  =  {\bf 1}_{{\bf M}^{{\leq}d} } $. \medskip\\
Then \medskip\\
1) The map  $\mathop{\sf T}^r({\bf 1}_{{\bf M}^{{\leq}d} 
{\subseteq}{\bf M}} )$  induces the isomorphism   of modules $\mathop{\sf 
T}^r({\bf M}^{{\leq}d} ) \rightarrow  \mathop{\sf T}^r({\bf M})^{{\leq}d} 
$; \medskip\\
2) The map  $\mathop{\sf TS}^r({\bf 1}_{{\bf M}^{{\leq}d} 
{\subseteq}{\bf M}} )$  induces the isomorphism  of modules $\mathop{\sf 
T}^r({\bf M}^{{\leq}d} ) \rightarrow  \mathop{\sf T}^r({\bf M})^{{\leq}d} 
$;\medskip\\
3)  The map  $\mathop{\sf TA}^r({\bf 1}_{{\bf M}^{{\leq}d} 
{\subseteq}{\bf M}} )$  induces the  isomorphism  of modules $\mathop{\sf 
TA}^r({\bf M}^{{\leq}d} ) \rightarrow  \mathop{\sf TA}^r({\bf 
M})^{{\leq}d} $;\medskip\\
4) The map  $\mathop{\sf A}^r({\bf 1}_{{\bf M}^{{\leq}d} 
{\subseteq}{\bf M}} )$  induce the isomorphism  of modules $\mathop{\sf 
A}^r({\bf M}^{{\leq}d} ) \rightarrow  \mathop{\sf A}^r({\bf M})^{{\leq}d} 
$.}\hfil\vfill\eject
 
 \noindent Let $(x_i|_{i{=}1,r} )$  be a tuple of variables. 
\medskip\\
{\bf Notation 1.36.} Denote by $0(x_i|_{i{=}1,r} )$ the zero element 
of the ring ${\bf R}[x_i|_{i{=}1,r} ]$,  and $1(x_i|_{i{=}1,r} )$ 
the unit element of the ring ${\bf R}[x_i|_{i{=}1,r} ]$.\medskip 
\\
{\bf Definition 1.4.} Define the action of permutations ${\tau} 
{\in}!r$ on  ${\bf R}[x_i|_{i{=}1,r} ]$ by \medskip

 ${\tau} \, h(x_i|_{i{=}1,r} ) =  h(x_{{\tau} ^{{-}1} (i)} |_{i{=}1,r} )$. 
\medskip
\\
{\it 1.For any permutations  ${\tau} ',{\tau} ''{\in}!r$ 
 and any polynomial  $h(x_i|_{i{=}1,r} )  \in $ ${\bf R}[x_i|_{i{=}1,r} 
]$ there holds \medskip

 ${\tau} '\, {\tau} ''\, h(x_i|_{i{=}1,r} ) =  {\tau} '\, ({\tau} 
''\, h(x_i|_{i{=}1,r} ))$. \medskip\\
2. For any  polynomial $h(x_i|_{i{=}1,r} ) \in  {\bf R}[x_i|_{i{=}1,r} ]$ 
holds \medskip

 ${\varepsilon} _r\, h(x_i|_{i{=}1,r} ) =  h(x_i|_{i{=}1,r} 
)$. } \medskip
\\
{\bf Notation  1.37.}  Denote by  
$\mathop{\sf ts}^r({\bf R}[x_i|_{i{=}1,r} ])$  
the set of all polynomials $h(x_i|_{i{=}1,r} )$ 
in 
${\bf R}[x_i|_{i{=}1,r} ]$ 
such  that for any transposition 
${\tau} $ in $!r$ there holds \medskip

 ${\tau} \, h(x_i|_{i{=}1,r} ) =  h(x_i|_{i{=}1,r} )$. \medskip
\\
{\bf Notation 1.38.} Denote by $\mathop{\sf ta}^r({\bf R}[x_i|_{i{=}1,r} ])$ 
the set of all 
$h(x_i|_{i{=}1,r} ) \in {\bf R}[x_i|_{i{=}1,r} ]$ 
such  that for any transposition 
${\tau} $ in $!r$ there holds \medskip

 ${\tau} \, h(x_i|_{i{=}1,r} ) =  {-}h(x_i|_{i{=}1,r} )$. \medskip
\\
{\it There holds \medskip

 $\mathop{\sf ts}^r({\bf R}[x_i|_{i{=}1,r} ]){\cdot}\mathop{\sf 
ts}^r({\bf R}[x_i|_{i{=}1,r} ]) \subseteq  \mathop{\sf ts}^r({\bf 
R}[x_i|_{i{=}1,r} ])$, \medskip

 $\mathop{\sf ts}^r({\bf R}[x_i|_{i{=}1,r} ]){\cdot}\mathop{\sf 
ta}^r({\bf R}[x_i|_{i{=}1,r} ]) \subseteq  \mathop{\sf ta}^r({\bf 
R}[x_i|_{i{=}1,r} ])$, \medskip

 $1_{{\bf R}[x_i|_{i{=}1,r} ]}  \in  \mathop{\sf ts}^r({\bf 
R}[x_i|_{i{=}1,r} ])$. } \medskip\\
Hence, \medskip
\\
{\it 1) the product  in  ${\bf R}[x_i|_{i{=}1,r} ]$  induces 
a structure of an associative  and commutative  algebra on  
$\mathop{\sf ts}^r({\bf R}[x_i|_{i{=}1,r} ])$ 
with the unity $1_{{\bf R}[x_i|_{i{=}1,r} ]} $. \medskip
\\
2) the product in  ${\bf R}[x_i|_{i{=}1,r} ]$ induces 
structure of an unitary module on 
$\mathop{\sf ta}^r({\bf R}[x_i|_{i{=}1,r} ])$ 
over the algebra 
$\mathop{\sf ts}^r({\bf R}[x_i|_{i{=}1,r} ])$ 
with the unity 
$1_{{\bf R}[x_i|_{i{=}1,r} ]} $. }\bigskip
\\
Let $x$  be a variable. Since ${\bf R}[x]$  
is a commutative algebra over ${\bf R}$, then also  
${\bf R}[x]^{{\otimes}r} $ 
is a commutative algebra over ${\bf R}$ 
as tensor product of commutative algebras 
over ${\bf R}$. \medskip
\\
{\bf Proposition 1.23.} {\it There is an isomorphism of algebras over 
${\bf R}$  \medskip

 ${\nu} {:}\ {\bf R}[x_i|_{i{=}1,r} ] \rightarrow  \mathop{\sf 
T}^r({\bf R}[x])$, \medskip

 $x_i \mapsto  1^{{\otimes}(i{-}1)} {\otimes}x{\otimes}1^{{\otimes}(r{-}i)} 
$ for $i{=}1,r$, \medskip

 $1(x_i|_{i{=}1,r} ) \mapsto  1(x)^{{\otimes}r} $. \medskip
\\
This isomorphism induces \medskip\\
1) the isomorphisms of algebras \medskip

 ${\nu} _s{:}\ \mathop{\sf ts}^r({\bf R}[x_i|_{i{=}1,r} ]) 
 \rightarrow   \mathop{\sf TS}^r({\bf R}[x])$, \medskip\\
this means that there holds \medskip

 ${\nu} _s({\bf a}'(x_i|_{i{=}1,r} ){\cdot}{\bf a}''(x_i|_{i{=}1,r} 
)) =  {\nu} _s({\bf a}'(x_i|_{i{=}1,r} )){\cdot}{\nu} _s({\bf 
a}''(x_i|_{i{=}1,r} ))$  \medskip\\
for any ${\bf a}'(x_i|_{i{=}1,r} )$, ${\bf a}''(x_i|_{i{=}1,r} 
) \in  \mathop{\sf ts}^r({\bf R}[x_i|_{i{=}1,r} ])$; \medskip
\\
2) the isomorphism of modules over algebras \medskip

 ${\nu} _a{:}\ \mathop{\sf ta}^r({\bf R}[x_i|_{i{=}1,r} ]) 
 \rightarrow   \mathop{\sf TA}^r({\bf R}[x])$, ${\nu} _s{:}\ 
\mathop{\sf ts}^r({\bf R}[x_i|_{i{=}1,r} ]) \rightarrow  \mathop{\sf 
TS}^r({\bf R}[x])$, \medskip\\
this means that there holds  \medskip

 ${\nu} _a({\bf a}'(x_i|_{i{=}1,r} ){\cdot}{\bf m}''(x_i|_{i{=}1,r} 
)) =  {\nu} _s({\bf a}'(x_i|_{i{=}1,r} )){\cdot}{\nu} _a({\bf 
m}''(x_i|_{i{=}1,r} ))$  \medskip\\
for any ${\bf a}'(x_i|_{i{=}1,r} ) \in  \mathop{\sf ts}^r({\bf 
R}[x_i|_{i{=}1,r} ])$, ${\bf m}''(x_i|_{i{=}1,r} ) \in  \mathop{\sf 
ta}^r({\bf R}[x_i|_{i{=}1,r} ])$. }\\
\\
{\bf Convention 1.2.} Throughout this article put  
${\bf x}_i{=}1^{{\otimes}(i{-}1)} {\otimes}x{\otimes}1^{{\otimes}(r{-}i)} $  
for $i{=}1,r$. 
\hfil\vfill\eject
 \noindent {\bf Notation 1.39.} Denote by ${\bf R}[x]^{{\leq}d} 
$  the set of all polynomials in ${\bf R}[x]$  of the degree  ${\leq}d$. 
\medskip
\\
{\it  ${\bf R}[x]^{{\leq}d} $  is  submodule   of a module 
  ${\bf R}[x]$. The ring ${\bf R}[x]$   is a semigraded ring, with 
a semigraded  $({\bf R}[x]^{{\leq}d} |d{\in}{\bbb Z})$.} \medskip
\\
{\bf Notation  1.40.}  Denote by  ${\bf R}[(x_i)^{{\leq}d} |_{i{=}1,r} ]$ 
the set of all polynomials in ${\bf R}[x_i|_{i{=}1,r} ]$ 
of the degree in $x_i$ which is ${\leq}d$ for $i{=}1,r$.\medskip 
\\
By notations 1.31, 1.32, 1.33, 1.34, 1.35  there holds 
\medskip

 $\mathop{\sf T}^r({\bf R}[x])^{{\leq}d}  =  \mathop{\rm Span}(\mathop{\sf 
T}^r_{{\bf R}[x]} ({\bf R}[x^{{\leq}d} ]|_{{\times}r} ))$, \medskip

 $\mathop{\sf TS}^r({\bf R}[x])^{{\leq}d}  =  \mathop{\sf TS}^r({\bf 
R}[x]){\cap}\mathop{\sf T}^r({\bf R}[x])^{{\leq}d} $, \medskip

 $\mathop{\sf TA}^r({\bf R}[x])^{{\leq}d}  =  \mathop{\sf TA}^r({\bf 
R}[x]){\cap}\mathop{\sf T}^r({\bf R}[x])^{{\leq}d} $, \medskip

 $\mathop{\sf A}^r({\bf R}[x])^{{\leq}d}  =  (\mathop{\sf T}^r({\bf 
R}[x])^{{\leq}d} {+}\mathop{\sf ST}^r({\bf R}[x]))/\mathop{\sf 
ST}^r({\bf R}[x])$. \medskip
\\
{\bf Proposition 1.24.} {\it There holds \medskip

 $\mathop{\sf TS}^r({\bf R}[x])^{{\leq}d'} {\cdot}\mathop{\sf 
TS}^r({\bf R}[x])^{{\leq}d''}  \subseteq  \mathop{\sf TS}^r({\bf 
R}[x])^{{\leq}d'{+}d''} $; \medskip

 $\mathop{\sf TS}^r({\bf R}[x])^{{\leq}d'} {\cdot}\mathop{\sf 
TA}^r({\bf R}[x])^{{\leq}d''}  \subseteq  \mathop{\sf TA}^r({\bf 
R}[x])^{{\leq}d'{+}d''} $; \medskip

 $\mathop{\sf TS}^r({\bf R}[x])^{{\leq}d'} {\cdot}\mathop{\sf 
A}^r({\bf R}[x])^{{\leq}d''}  \subseteq  \mathop{\sf A}^r({\bf 
R}[x])^{{\leq}d'{+}d''} $; \medskip\\
i.    e.    $\mathop{\sf ts}^r({\bf R}[x_i|_{i{=}1,r} ])$  
  is a semigraded     algebra, $\mathop{\sf 
ta}^r({\bf R}[x_i|_{i{=}1,r} ])$ is a semigraded 
module  over the semigraded algebra $\mathop{\sf ts}^r({\bf 
R}[x_i|_{i{=}1,r} ])$. } \medskip\\
{\bf Notation 1.41.} Define \medskip

 $\mathop{\sf t}^r({\bf R}[x_i|_{i{=}1,r} ])^{{\leq}d}  =  
{\bf R}[(x_i)^{{\leq}d} |_{i{=}1,r} ]$, \medskip

 $\mathop{\sf ts}^r({\bf R}[x_i|_{i{=}1,r} ])^{{\leq}d}  = 
 \mathop{\sf ts}^r({\bf R}[x_i|_{i{=}1,r} ]){\cap}\mathop{\sf 
t}^r({\bf R}[x_i|_{i{=}1,r} ])^{{\leq}d} $, \medskip

 $\mathop{\sf ta}^r({\bf R}[x_i|_{i{=}1,r} ])^{{\leq}d}  = 
 \mathop{\sf ta}^r({\bf R}[x_i|_{i{=}1,r} ]){\cap}\mathop{\sf 
t}^r({\bf R}[x_i|_{i{=}1,r} ])^{{\leq}d} $. \medskip\\
{\bf Proposition 1.25.} {\it There holds \medskip

 $\mathop{\sf ts}^r({\bf R}[x_i|_{i{=}1,r} ])^{{\leq}d'} {\cdot}\mathop{\sf 
ts}^r({\bf R}[x_i|_{i{=}1,r} ])^{{\leq}d''}  \subseteq  \mathop{\sf 
ts}^r({\bf R}[x_i|_{i{=}1,r} ])^{{\leq}d'{+}d''} $; \medskip

 $\mathop{\sf ts}^r({\bf R}[x_i|_{i{=}1,r} ])^{{\leq}d'} {\cdot}\mathop{\sf 
ta}^r({\bf R}[x_i|_{i{=}1,r} ])^{{\leq}d''}  \subseteq  \mathop{\sf 
ta}^r({\bf R}[x_i|_{i{=}1,r} ])^{{\leq}d'{+}d''} $. \medskip
\\
i.    e.    $\mathop{\sf ts}^r({\bf R}[x_i|_{i{=}1,r} ])$  
  is a semigraded     algebra, $\mathop{\sf 
ta}^r({\bf R}[x_i|_{i{=}1,r} ])$ is a semigraded 
module  over the semigraded algebra $\mathop{\sf ts}^r({\bf 
R}[x_i|_{i{=}1,r} ])$. }\medskip 
\\ 
{\bf Proposition 1.26.} {\it There holds \medskip

 ${\nu} (\mathop{\sf t}^r({\bf R}[x_i|_{i{=}1,r} ])^{{\leq}d} 
) =  \mathop{\sf T}^r({\bf R}[x])^{{\leq}d} $, \medskip

 ${\nu} _s(\mathop{\sf ts}^r({\bf R}[x_i|_{i{=}1,r} ])^{{\leq}d} 
) =  \mathop{\sf TS}^r({\bf R}[x])^{{\leq}d} $, \medskip

 ${\nu} _a(\mathop{\sf ta}^r({\bf R}[x_i|_{i{=}1,r} ])^{{\leq}d} 
) =  \mathop{\sf TA}^r({\bf R}[x])^{{\leq}d} $. 
\medskip
\\
The isomorphism  of algebras  \medskip

 ${\nu} {:}\ {\bf R}[x_i|_{i{=}1,r} ] \rightarrow  \mathop{\sf T}^r({\bf R}[x])$, \medskip
\\
induces \medskip
\\
1) the isomorphism of semigraded of algebras \medskip

 ${\nu} _s{:}\ \mathop{\sf ts}^r({\bf R}[x_i|_{i{=}1,r} ]) 
 \rightarrow   \mathop{\sf TS}^r({\bf R}[x])$, \medskip
\\
2) the isomorphism of semigraded modules  over  semigraded 
 algebras \medskip

 ${\nu} _a{:}\ \mathop{\sf ta}^r({\bf R}[x_i|_{i{=}1,r} ]) 
\rightarrow  \mathop{\sf TA}^r({\bf R}[x])$, ${\nu} _s{:}\ \mathop{\sf 
ts}^r({\bf R}[x_i|_{i{=}1,r} ]) \rightarrow  \mathop{\sf TS}^r({\bf 
R}[x])$. }\medskip 
\\
{\bf Notation 1.42.} Denote by \medskip

 ${\sigma} (x_i|_{i{=}1,r} )(y) =  \prod ((y{-}x_i)|_{i{=}1,r} 
) =  \sum ({\sigma} _p(x_i|_{i{=}1,r} ){\cdot}y^{r{-}p} |_{p{=}0,r} 
)$. \medskip\\
Since ${\sigma} (x_i|_{i{=}1,r} )(y)$ symmetrically  
depends from  $(x_i|_{i{=}1,r} )$,  then  ${\sigma} _p(x_i|_{i{=}1,r} )$ 
is a symmetric polynomial in $(x_i|_{i{=}1,r} )$. Note 
 that  ${\sigma} (x_i|_{i{=}1,r} )(y)$ has in $x_i$ the degree 
${\leq}1$, hence, there holds \medskip
\\
{\bf Proposition 1.27.} \medskip

 ${\sigma} _p(x_i|_{i{=}1,r} ) \in  \mathop{\sf ta}^r({\bf 
R}[x_i|_{i{=}1,r} ])^{{\leq}1} $ for $p{=}0,r$.
\hfil\vfill\eject
 
 \noindent {\bf Proposition 1.28.} {\it Let $\mathfrak{ p} ^d_{{\bf R}[x]} 
{:}\ {\bf R}[x] \rightarrow   {\bf R}[x]^{{\leq}d} $   
be a linear  map such  that 
$x^{\delta}  \mapsto  ?({\delta} {\leq}d){\cdot}x^{\delta} $ 
for ${\delta} {\geq}0$. 
Then \medskip

 $\mathfrak{ p} ^d_{\bf M}\, {\bf 1}_{{\bf M}^{{\leq}d} {\subseteq}{\bf 
M}}  =  {\bf 1}_{{\bf M}^{{\leq}d} } $, \medskip\\
and \medskip
\\
1) The map $\mathop{\sf T}^r({\bf 1}_{{\bf R}[x]^{{\leq}d} 
{\subseteq}{\bf R}[x]} )$ induces the isomorphism of modules \medskip

 $\mathop{\sf T}^r({\bf R}[x]^{{\leq}d} ) \rightarrow  \mathop{\sf 
T}^r({\bf R}[x])^{{\leq}d} $. \medskip
\\
2) The map $\mathop{\sf TS}^r({\bf 1}_{{\bf R}[x]^{{\leq}d} 
{\subseteq}{\bf R}[x]} )$ induces the isomorphism of modules \medskip

 $\mathop{\sf T}^r({\bf R}[x]^{{\leq}d} ) \rightarrow  \mathop{\sf 
T}^r({\bf R}[x])^{{\leq}d} $. \medskip
\\
3) The map $\mathop{\sf TA}^r({\bf 1}_{{\bf R}[x]^{{\leq}d} 
{\subseteq}{\bf R}[x]} )$ induces the isomorphism of modules \medskip

 $\mathop{\sf TA}^r({\bf R}[x]^{{\leq}d} ) \rightarrow  \mathop{\sf 
TA}^r({\bf R}[x])^{{\leq}d} $.\medskip
\\
4) The map $\mathop{\sf A}^r({\bf 1}_{{\bf R}[x]^{{\leq}d} 
{\subseteq}{\bf R}[x]} )$ induces the isomorphism of modules \medskip

 $\mathop{\sf A}^r({\bf R}[x]^{{\leq}d} ) \rightarrow  \mathop{\sf 
A}^r({\bf R}[x])^{{\leq}d} $.} \medskip
\\
{\bf Proof.} The first statement proposition  obviously. 
 Then  1),  2),  3) proposition it follows from proposition 1.22 
 under substitution ${\bf M}^{{\leq}d}  \mapsto  {\bf R}[x]^{{\leq}d} 
$ for $d{\in}{\bbb Z}$, ${\bf M} \mapsto  {\bf R}[x]$.\bigskip 
\\
{\bf Definition 1.5.} Let $q'{\leq}p'$, $q''{\leq}p''$. 
A map ${\rho} {:}\ [q',p'] \rightarrow  [q'',p'']$ we call 
monotonous on  ${\Omega} {\subseteq}[q',p']$, if for any $l,k$ 
in ${\Omega} $ such  that  $l{<}k$  there holds ${\rho} (l){<}{\rho} (k)$. \medskip
\\
{\bf Notation  1.43.} Denote by ${\mathcal C} ^r_p$ 
the set of all ${\tau} {\in}!p$  such  which  are monotonous 
on  $[1,q]$ and on  $[q{+}1,p]$, where $0{\leq}r{\leq}p$, $q{=}p{-}r$. 
\medskip
\\
{\bf Notation 1.44.} Denote by ${\bf C}^r_p$  
the set of all maps 
${\chi} {:}\ [1,p]  \rightarrow   \{0,1\}$ 
such  that 
$\sum {\chi}  =  r$, 
where 
$0{\leq}r{\leq}p$. \medskip
\\
{\bf Notation 1.45.} Denote by 
$\lnot $ the map $\{0,1\} \rightarrow   \{0,1\}$  
such that  
${\lnot}1{=}0$, ${\lnot}0{=}1$. 
\medskip
\\
{\bf Notation 1.46.} Let ${\chi} {:}\ [1,p] \rightarrow 
 \{0,1\}$  be a map.  
Let  $h{=}(h_i|_{i{=}1,p} )$, denote by \medskip

 $h|_{{\times}{\chi} } \ {{:}{=}}\ (h_i|_{i{=}1,p} |_{{\times}{\chi} 
} )\ {{:}{=}}\ (h_i|_{{\times}{\chi} (i)} |_{i{=}1,p} )$. \medskip
\\
{\bf Proposition 1.29.} {\it Each ${\chi} {\in}{\bf C}^r_p$ 
corresponds unique ${\tau} {\in}{\mathcal C} ^r_p$ such 
 that  ${\chi} ({\tau} (l))  =   1$  for  all  $l{\in}[1,q]$, 
 ${\chi} ({\tau} (l))  =   0$  for  all $l{\in}[q{+}1,p]$, where 
$1{\leq}r{\leq}p$, $q{=}p{-}r$. There holds \medskip

 $(a_{{\tau} (i)} |_{i{=}1,q} ,b_{{\tau} (i)} |_{i{=}q{+}1,p} 
) =  (a_{l'} |_{{\times}{\lnot}{\chi} (l')} |_{l'{=}1,p} ,b_{l''} 
|_{{\times}{\chi} (l'')} |_{l''{=}1,p} )$  \medskip

 \qquad $=  (a_i|_{i{=}1,p} |_{{\times}{\lnot}{\chi} } ,b_i|_{i{=}1,p} 
|_{{\times}{\chi} } ) =  (a|_{{\times}{\lnot}{\chi} } ,b|_{{\times}{\chi} 
} )$. } \medskip\\
{\bf Notation 1.47.} In proposition 1.29 denote by \medskip

 $\mathop{\rm sgn}({\chi} ) =  \mathop{\rm sgn}({\tau} )$. 
\medskip\\
{\bf Definition 1.6.} Let $H$ is a module, ${\bf A}$ 
is an associative algebra with an unity, 
and be given a bilinear map 
${\bf A}{\times}\bigwedge  ^r(H) \rightarrow   \bigwedge  ^r(H)$  
of modules such  that $\bigwedge  ^r(H)$ 
is a module over ${\bf A}$, where $r{\geq}0$.\medskip
\\
Let $a^k_j{\in}{\bf A}$ for $j{=}1,p$, for $k{=}1,q$,  let 
 $h_j{\in}\bigwedge  ^1(H)$  for  $j{=}1,p$, where $q{=}p{-}r$. 
Denote by \bigskip

 ${\det}\left\| \begin{matrix}    & a^k_j|^{k{=}1,q} _{j{=}1,p} 
 \cr
  &   \cr
\wedge  & h_j|_{j{=}1,p} \end{matrix} \right\|  =  \sum (\mathop{\rm 
sgn}({\tau} ){\cdot}({\det}\left\| a^k_{{\tau} (l)} |^{k{=}1,q} 
_{l{=}1,q} \right\| {\cdot}\bigwedge (h_{{\tau} (l)} |_{l{=}q{+}1,p} 
)|{\tau} {\in}{\mathcal C} ^r_p)$. \bigskip\\
\\
{\it There holds \bigskip

 ${\det}\left\| \begin{matrix}    & a^k_j|^{k{=}1,q} _{j{=}1,p} 
 \cr
  &   \cr
\wedge  & h_j|_{j{=}1,p} \end{matrix} \right\|  =  \sum (\mathop{\rm 
sgn}({\chi} ){\cdot}({\det}\left\| a^k_i|^{k{=}1,q} _{i{=}1,p} 
|_{{\times}{\lnot}{\chi} } \right\| {\cdot}\bigwedge (h_i|_{i{=}1,p} 
|_{{\times}{\chi} } )|{\chi} {\in}{\bf C}^r_p)$. }\\
\hfil\vfill\eject
 
 \noindent {\bf Case 1.} Let $q{=}0$. Then $r{=}p{-}q{=}p{-}0{=}p$ 
and \bigskip

 ${\det}\left\| \begin{matrix}    & a^k_j|^{k{=}1,q} _{j{=}1,p} 
 \cr
  &   \cr
\wedge  & h_j|_{j{=}1,p} \end{matrix} \right\|  =  {\det}\left\| 
\begin{matrix}\wedge  & h_j|_{j{=}1,p} \end{matrix} \right\| 
$, \bigskip

 $\sum (\mathop{\rm sgn}({\tau} ){\cdot}({\det}\left\| a^k_{{\tau} 
(l)} |^{k{=}1,q} _{l{=}1,q} \right\| {\cdot}\bigwedge (h_{{\tau} 
(l)} |_{l{=}q{+}1,p} ))|{\tau} {\in}{\mathcal C} ^r_p)$  \medskip

 \qquad $=  (\mathop{\rm sgn}({\tau} ){\cdot}\bigwedge (h_{{\tau} 
(l)} |_{l{=}1,p} ))|{\tau} {\in}{\mathcal C} ^p_p) =  \bigwedge 
(h_l|_{l{=}1,p} )$. \bigskip\\
Here the condition ${\tau} {\in}{\mathcal C} ^p_p$ means that 
$t$ is monotonous on  $[1,0]$  and  on   $[1,p]$, such map 
is only  the identity  map ${\varepsilon} _p{:}\ 
[1,p] \rightarrow  [1,p]$, and then $\mathop{\rm sgn}({\tau} ) 
=  1$. Hence, \medskip

 ${\det}\left\| \begin{matrix}\wedge  & h_j|_{j{=}1,p} \end{matrix} 
\right\|  =  \bigwedge (h_l|_{l{=}1,p} )$.\\
\\
{\bf Case 2.} Let $q{=}p$. Then $r{=}p{-}q{=}p{-}p{=}0$ 
and  \bigskip

 ${\det}\left\| \begin{matrix}    & a^k_j|^{k{=}1,q} _{j{=}1,p} 
 \cr
  &   \cr
\wedge  & h_j|_{j{=}1,p} \end{matrix} \right\|  =  {\det}\left\| 
\begin{matrix}    & a^k_j|^{k{=}1,p} _{j{=}1,p}  \cr
  &   \cr
\wedge  & h_j|_{j{=}1,p} \end{matrix} \right\| $. \bigskip

 $\sum (\mathop{\rm sgn}({\tau} ){\cdot}({\det}\left\| a^k_{{\tau} 
(l)} |^{k{=}1,q} _{l{=}1,q} \right\| {\cdot}\bigwedge (h_{{\tau} 
(l)} |_{l{=}q{+}1,p} ))|{\tau} {\in}{\mathcal C} ^r_p)$  \bigskip

 \qquad $=  \sum (\mathop{\rm sgn}({\tau} ){\cdot}({\det}\left\| 
a^k_{{\tau} (l)} |^{k{=}1,p} _{l{=}1,p} \right\| {\cdot}\bigwedge 
(h_{{\tau} (l)} |_{l{=}p{+}1,p} ))|{\tau} {\in}{\mathcal C} ^0_p)$ 
 \bigskip

 \qquad $=  {\det}\left\| a^k_{{\tau} (l)} |^{k{=}1,p} _{l{=}1,p} 
\right\| $. \bigskip\\
Here the condition ${\tau} {\in}{\mathcal C} ^r_p$ means that 
$t$ is monotonous $[1,p]$ and on   $[p{+}1,p]$,  and such map 
is only the identity  map ${\varepsilon} _p{:}\ 
[1,p] \rightarrow  [1,p]$, and then $\mathop{\rm sgn}({\tau} ) 
=  1$.  There holds  $\bigwedge (h_{{\tau} (l)} |_{l{=}p{+}1,p} 
)  =   1$. Hence, \medskip

 ${\det}\left\| \begin{matrix}    & a^k_j|^{k{=}1,p} _{j{=}1,p} 
 \cr
  &   \cr
\wedge  & h_j|_{j{=}1,p} \end{matrix} \right\|  =  {\det}\left\| 
a^k_{{\tau} (l)} |^{k{=}1,p} _{l{=}1,p} \right\| $.\medskip 
\\\\
{\bf Lemma 1.1.} {\it Let $h_i(x){\in}{\bf R}[x]$ 
for  $i{=}1,r$.  
There holds \medskip

 $\mathop{\sf atr}^r(\bigwedge (h_i(x)|_{i{=}1,r} )) =  {\det}\| 
h_j({\bf x}_i)|^{i{=}1,r} _{j{=}1,r} \|  =  {\nu} _a({\det}\| 
h_j(x_i)|^{i{=}1,r} _{j{=}1,r} \| )$. } \medskip\\
{\bf Proof.} There holds \medskip

 $\mathop{\sf atr}^r(\bigwedge (h_i(x)|_{i{=}1,r} ))$  \medskip

 \qquad$=  \sum (\mathop{\rm sgn}({\tau} ){\cdot}({\tau} \, 
\bigotimes (h_i(x)|_{i{=}1,r} )))$  \medskip

 \qquad $=  \sum (\mathop{\rm sgn}({\tau} ){\cdot}\bigotimes 
(h_{{\tau} ^{{-}1} (i)} (x)|_{i{=}1,r} )|{\tau} {\in}!r)$  \medskip

 \qquad $=  \sum (\mathop{\rm sgn}({\tau} ){\cdot}\prod (h_{{\tau} 
^{{-}1} (i)} ({\bf x}_i)|_{i{=}1,r} )|{\tau} {\in}!r)$ \medskip

 \qquad $=  {\det}\| h_j({\bf x}_i)|^{i{=}1,r} _{j{=}1,r} \| $. \bigskip
\hfil\vfill\eject

\noindent 
{\bf Lemma 1.2.} {\it  Let  $a^k_j(x_i|_{i{=}1,r} ){\in}\mathop{\sf 
ts}^r({\bf R}[x_i|_{i{=}1,r} ])$  for  $j{=}1,p$,  for $k{=}1,q$; 
$h_j(x){\in}{\bf R}[x]$ for $j{=}1,p$; where $q{=}p{-}r$. Then \bigskip

 $\mathop{\sf atr}^r\, ({\det}\left\| \begin{matrix}    & a^k_j({\bf 
x}_i|_{i{=}1,r} )|^{k{=}1,q} _{j{=}1,p}  \cr
  &   \cr
\wedge  & h_j(x)\h     |_{j{=}1,p} \end{matrix} \right\| ) 
=  {\det}\left\| \begin{matrix}a^k_j({\bf x}_i|_{i{=}1,r} )|^{k{=}1,q} 
_{j{=}1,p}  \cr
  \cr
h_j({\bf x}_i)   \h |^{i{=}1,r} _{j{=}1,p} \end{matrix} \right\| 
 =  {\nu} _a({\det}\left\| \begin{matrix}a^k_j(x_i|_{i{=}1,r} 
)|^{k{=}1,q} _{j{=}1,p}  \cr
  \cr
h_j(x_i)   \h |^{i{=}1,r} _{j{=}1,p} \end{matrix} \right\| 
)$. }  \bigskip\\
{\bf Proof.}      Since      $a^k_j(x_i|_{i{=}1,r} 
){\in}\mathop{\sf ts}^r({\bf R}[x_i|_{i{=}1,r} ])$,      then $a^k_j({\bf 
x}_i|_{i{=}1,r} ){\in}\mathop{\sf TS}^r({\bf R}[x])$. There holds 
\bigskip

 $\mathop{\sf atr}^r\, ({\det}\left\| \begin{matrix}    & a^k_j({\bf 
x}_i|_{i{=}1,r} )|^{k{=}1,q} _{j{=}1,p}  \cr
  &   \cr
\wedge  & h_j(x)\h     |_{j{=}1,p} \end{matrix} \right\| )$ 
 \bigskip

 \qquad $=  \mathop{\sf atr}^r\, (\sum (\mathop{\rm sgn}({\tau} 
){\cdot}{\det}\left\| a^k_{{\tau} (l)} ({\bf x}_i|_{i{=}1,r} 
)|^{k{=}1,q} _{l{=}1,q} \right\| {\cdot}\bigwedge (h_{{\tau} 
(l)} (x)|_{l{=}q{+}1,p} )|{\tau} {\in}{\mathcal C} ^r_p))$  \bigskip

 \qquad $=  \sum (\mathop{\rm sgn}({\tau} ){\cdot}{\det}\left\| 
a^k_{{\tau} (l)} ({\bf x}_i|_{i{=}1,r} )|^{k{=}1,q} _{l{=}1,q} 
\right\| {\cdot}\mathop{\sf atr}^r\, (\bigwedge (h_{{\tau} (l)} 
(x)|_{l{=}q{+}1,p} ))|{\tau} {\in}{\mathcal C} ^r_p)$  \bigskip

 \qquad $=  \sum (\mathop{\rm sgn}({\tau} ){\cdot}{\det}\left\| 
a^k_{{\tau} (l)} ({\bf x}_i|_{i{=}1,r} )|^{k{=}1,q} _{l{=}1,q} 
\right\| {\cdot}{\det}\left\| h_{{\tau} (l)} ({\bf x}_i)|^{i{=}1,r} 
_{l{=}q{+}1,p} \right\| |{\tau} {\in}{\mathcal C} ^r_p)$  \bigskip

 \qquad $=  {\det}\left\| \begin{matrix}a^k_j({\bf x}_i|_{i{=}1,r} 
)|^{k{=}1,q} _{j{=}1,p}  \cr
  \cr
h_j({\bf x}_i)   \h |^{i{=}1,r} _{j{=}1,p} \end{matrix} \right\| $.\bigskip 
\\
\\
\\
{\bf Notation 1.48.} Denote by  $[x^{\alpha} ]_{*}$ 
the linear  map  ${\bf R}[x]  \rightarrow   {\bf R}$ such that \medskip

 $[x^{\alpha} ]_{*}.x^{\beta}  =  ?({\alpha} {=}{\beta} )$. \medskip
\\
{\bf Notation 1.49.} Denote by $x^{{\leq}d}  =  \| x^{d{-}{\delta} 
} |_{{\delta} {=}0,d} \| $, $[y^{{\leq}d} ]_{*} =  \| [y^{d{-}{\delta} 
} ]_{*}|^{{\delta} {=}0,{\delta} } \| $. \medskip\\
{\bf Notation 1.50.}  Let  $(x'_{i'} |_{i'{=}1,r'} )$, 
 $(x''_{i''} |_{i''{=}1,r''} )$ be tuples of variables. In general linear 
maps \medskip

 ${\bf R}[x''_{i''} |_{i''{=}1,r''} ] \rightarrow  {\bf R}[x'_{i'} 
|_{i'{=}1,r'} ]$  \medskip\\
we will denote by $a(x'_{i'} |_{i'{=}1,r'} ;x''_{i''} |_{i''{=}1,r''} )_{*}$, 
where $a$ is a symbol. \medskip
\\
{\bf Notation 1.51.} Let ${\bf A}$  be commutative algebra.  
Let  $x{=}(x_i|_{i{=}1,r} )$   be variables, 
$a{=}(a_i|_{i{=}1,r} )$  be  elements in ${\bf A}$.  
Denote by  ${\bf 1}(a_i|_{i{=}1,r} ;x_i|_{i{=}1,r} )_{*}$ 
the homomorphism of algebras ${\bf R}[x_i|_{i{=}1,r} ] \rightarrow 
 {\bf A}$ such  that  $x_i  \mapsto   a_i$  for  $i{=}1,r$. \medskip
\\
{\bf Proposition 1.30.} {\it Let  $x{=}(x_i|_{i{=}1,r} )$ 
  be  variables.  Then the map ${\bf 1}(x;x)_{*}$ is 
the identity map ${\bf R}[x] \rightarrow  {\bf R}[x]$.} 
\medskip
\\
{\bf Definition 1.7.} Let $x'{=}(x'_{i'} |_{i'{=}1,r'} )$, 
 $y'{=}(y_{j'} |_{j'{=}1,t'} )$,  $x''{=}(x''_{i''} |_{i''{=}1,r''} 
)$, $y''{=}(y_{j''} |_{j''{=}1,t''} )$,\\
$z{=}(z_k|_{k{=}1,s} )$  be tuples of variables. Let be given linear 
maps \medskip

 $a(x';y',z)_{*}{:}\ {\bf R}[y',z]\ (\ \simeq {\bf R}[y']{\otimes}{\bf 
R}[z]) \rightarrow  {\bf R}[x']$, \medskip

 $b(z,x'';y'')_{*}{:}\ {\bf R}[y''] \rightarrow  {\bf R}[z,x'']\ 
(\ \simeq {\bf R}[z]{\otimes}{\bf R}[x''])$. \medskip
\\
Denote by \medskip

 $a(x';y',z)_{*}\, b(z,x'';y'')_{*}$  \medskip

 \qquad $=  (a(x';y',z)_{*}{\otimes}{\bf 1}(x'';x'')_{*})\, 
({\bf 1}(y';y')_{*}{\otimes}b(z,x'';y'')_{*})$, \medskip
\\
and call composition of these maps. {\it This is a linear map } \medskip

 ${\bf R}[y',y'']\ (\ \simeq {\bf R}[y']{\otimes}{\bf R}[y'']) 
\rightarrow  {\bf R}[x',x'']\ (\ \simeq {\bf R}[x']{\otimes}{\bf R}[x''])$.\\
\hfil\vfill\eject
 
\section{Main results}

\noindent
\normalsize {\bf Transition 2.1.} {\it Let be given matrixs:  $A$ 
 of the size  $|^l_k$,  $B$  of the size $|^l_{d{+}1} $, $C$ of the size 
$|^l_m$ with elements in  ${\bf R}[x_i|_{i{=}1,r} ]$,  where  $l{+}r 
 =   k{+}m{+}d{+}1$.  Then \footnotesize\\

 ${\sigma} _0(x_i|_{i{=}1,r} )^{{\Delta} ''} {\cdot}{\sigma} 
_r(x_i|_{i{=}1,r} )^{{\Delta} '} {\cdot}{\det}\left\| \begin{matrix}A 
& B\h & C \cr
  &             &   \cr
  & x^{{\leq}d} _i|^{i{=}1,r}  &  \end{matrix} \right\|  = 
 ({-}1)^{r{\cdot}{\Delta} '} {\cdot}{\det}\left\| \begin{matrix}A 
& B\h & C \cr
  &                 &   \cr
  & x^{{\leq}d} _i{\cdot}x^{{\Delta} '} _i|^{i{=}1,r}  & 
 \end{matrix} \right\| $. } \normalsize \bigskip\\
{\bf Proof.} The matrix of the second determinant is obtained 
 from  matrix of the first determinant by follows transformations rows 
 \bigskip \footnotesize

 $x^{{\Delta} '} _i{\cdot}\left\| \begin{matrix}0|_{{\times}k} 
 & x^{{\leq}d} _i|^{i{=}1,r}  & 0|_{{\times}m} \end{matrix} \right\| 
 =  \left\| \begin{matrix}0|_{{\times}k}  & x^{{\leq}d} _i|^{i{=}1,r} 
{\cdot}x^{{\Delta} '} _i & 0|_{{\times}m} \end{matrix} \right\| 
$ for $i{=}1,r$. \bigskip \normalsize
\\
Then the determinant of the second matrix is equal to the determinant of the first 
matrix multiplied by the product \bigskip \footnotesize

 $\prod (x^{{\Delta} '} _i|_{i{=}1,r} ) =  \prod (x_i|_{i{=}1,r} 
)^{{\Delta} '}  =  ({-}1)^{{-}r{\cdot}{\Delta} '} {\cdot}\prod 
({-}x_i|_{[{=}1,r} )^{{\Delta} '} $  \bigskip

 \qquad $=  ({-}1)^{{-}r{\cdot}{\Delta} '} {\cdot}{\sigma} 
_0(x_i|_{i{=}1,r} )^{{\Delta} ''} {\cdot}{\sigma} _r(x_i|_{i{=}1,r} 
)^{{\Delta} '} $. \bigskip \normalsize\\
The last equality  holds since  ${\sigma} _0(x_i|_{i{=}1,r} 
)  =   1(x_i|_{i{=}1,r} )$.  This it implies equality being proved.
\\
\\
\normalsize {\bf Transition 2.2.} {\it Let be given vectors: $a(y)$ 
of the size $|_k$, $b(y)$ of the size $|_l$, $c(y)$ of the size $|_m$ 
of polynomials in ${\bf R}[x_i|_{i{=}1,r} ][y]$ of degrees ${\leq}{\delta} $, 
where $r{+}{\delta} {+}1  =   k{+}l{+}m$.  Then \footnotesize\\

 ${\det}\left\| \begin{matrix}[y^{{\leq}{\delta} } ]_{*}.| 
& a(y){\cdot}{\sigma} (x_i|_{i{=}1,r} )(y) & b(y) & c(y)\h \cr
  &                       &        &   \cr
           &                       &        & c(x_i)|^{i{=}1,r} 
\end{matrix} \right\| $\\
\\

 \qquad $=  {\det}\left\| \begin{matrix}[y^{{\leq}{\delta} 
} ]_{*}.| & a(y){\cdot}{\sigma} (x_i|_{i{=}1,r} )(y) & \,   b(y)\h 
& c(y) \cr
  &                       &   &   \cr
           &                       & {-}b(x_i)|^{i{=}1,r}  
&  \end{matrix} \right\| $. } \normalsize\\
\\
\\
{\bf Proof}. The matrix of the second determinant is obtained 
 from the matrix of  the second determinant by the following transformations 
its rows for $i{=}1,r$  \footnotesize\\

 $ \left\| \begin{matrix}0|_{{\times}k}  & 0|_{{\times}k}  
& c(x_i)\end{matrix} \right\|  -  \sum\limits^{ {\delta} } _{{\alpha} 
{=}0} x^{\alpha} _i{\cdot}(\left\| \begin{matrix}[y^{\alpha} 
]_{*}.| & a(y){\cdot}{\sigma} (x_i|_{i{=}1,r} )(y) & b(y) & c(y)\h\end{matrix} 
\right\| )$  \bigskip

 \qquad $=   \left\| \begin{matrix}0|_{{\times}k}  & 0|_{{\times}l} 
 & c(x_i)\end{matrix} \right\| $  \medskip

 \qquad \qquad $-  \left\| \begin{matrix}\sum\limits^{ {\delta} 
} _{{\alpha} {=}0} x^{\alpha} _i{\cdot}[y^{\alpha} ]_{*}.a(y){\cdot}{\sigma} 
(x_i|_{i{=}1,r} )(y) & \sum\limits^{ {\delta} } _{{\alpha} {=}0} 
x^{\alpha} _i{\cdot}[y^{\alpha} ]_{*}.b(y) & \sum\limits^{ {\delta} 
} _{{\alpha} {=}0} x^{\alpha} _i{\cdot}[y^{\alpha} ]_{*}.c(y)\end{matrix} 
\right\| $  \bigskip

 \qquad $=   \left\| \begin{matrix}0|_{{\times}k}  & 0|_{{\times}l} 
 & c(x_i)\end{matrix} \right\|  -  \left\| \begin{matrix}a(x_i){\cdot}{\sigma} 
(x_i|_{i{=}1,r} )(x_i) & b(x_i) & c(x_i)\end{matrix} \right\| 
$  \bigskip

 \qquad $=   \left\| \begin{matrix}0|_{{\times}k}  & 0|_{{\times}l} 
 & c(x_i)\end{matrix} \right\|  -  \left\| \begin{matrix}a(x_i){\cdot}0(x_i|_{i{=}1,r} 
) & b(x_i) & c(x_i)\end{matrix} \right\| $  \bigskip

 \qquad $=   \left\| \begin{matrix}0|_{{\times}k}  & 0|_{{\times}l} 
 & c(x_i)\end{matrix} \right\|  -  \left\| \begin{matrix}0|_{{\times}k} 
 & b(x_i) & c(x_i)\end{matrix} \right\| $  \bigskip

 \qquad $=   \left\| \begin{matrix}0|_{{\times}k} {-}0|_{{\times}k} 
 & 0|_{{\times}l} {-}b(x_i) & c(x_i){-}c(x_i)\end{matrix} \right\| 
$  \bigskip

 \qquad $=   \left\| \begin{matrix}0|_{{\times}k}  & {-}b(x_i) 
& 0|_{{\times}m} \end{matrix} \right\| $. \normalsize \bigskip
\\
Then determinants of the both  matrixs are equal.
\hfil\vfill\eject
 
 \noindent {\bf Transition 2.3.} {\it There holds \bigskip \footnotesize

 ${\det}\left\| \begin{matrix}[y^{{\leq}d{+}{\Delta} } ]_{*}.| 
& y^{{\leq}{\Delta} ''{-}1} {\cdot}y^{d{-}r{+}1{+}{\Delta} 
'} {\cdot}{\sigma} (x_i|_{i{=}1,r} )(y) & y^{{\leq}d} {\cdot}y^{{\Delta} 
'}  & y^{{\leq}{\Delta} '{-}1} {\cdot}{\sigma} (x_i|_{i{=}1,r} 
)(y)\end{matrix} \right\| $  \bigskip

 \quad $=  {\sigma} _0(x_i|_{i{=}1,r} )^{{\Delta} ''} {\cdot}{\sigma} 
_r(x_i|_{i{=}1,r} )^{{\Delta} '} $. } \bigskip \normalsize
\\
{\bf Proof.} \footnotesize\\

 ${\det}\left\| \begin{matrix}[y^{{\leq}d{+}{\Delta} } ]_{*}.| 
& y^{{\leq}{\Delta} ''{-}1} {\cdot}y^{d{-}r{+}1{+}{\Delta} 
'} {\cdot}{\sigma} (x_i|_{i{=}1,r} )(y) & y^{{\leq}d} {\cdot}y^{{\Delta} 
'}  & y^{{\leq}{\Delta} '{-}1} {\cdot}{\sigma} (x_i|_{i{=}1,r} 
)(y)\end{matrix} \right\| $  \bigskip\\

 \quad $=  {\det}\left\| \begin{matrix}[y^{{\leq}{\Delta} 
''{-}1} {\cdot}y^{d{+}{\Delta} '{+}1} ]_{*} \h & .| & y^{{\leq}{\Delta} 
''{-}1} {\cdot}y^{d{-}r{+}1{+}{\Delta} '} {\cdot}{\sigma} 
(x_i|_{i{=}1,r} )(y) & y^{{\leq}d} {\cdot}y^{{\Delta} '}  
& y^{{\leq}{\Delta} '{-}1} {\cdot}{\sigma} (x_i|_{i{=}1,r} 
)(y) \cr
[y^{{\leq}d} {\cdot}y^{{\Delta} '} ]_{*} \h & .| & y^{{\leq}{\Delta} 
''{-}1} {\cdot}y^{d{-}r{+}1{+}{\Delta} '} {\cdot}{\sigma} 
(x_i|_{i{=}1,r} )(y) & y^{{\leq}d} {\cdot}y^{{\Delta} '}  
& y^{{\leq}{\Delta} '{-}1} {\cdot}{\sigma} (x_i|_{i{=}1,r} 
)(y) \cr
[y^{{\leq}{\Delta} '{-}1} ]_{*} \h & .| & y^{{\leq}{\Delta} 
''{-}1} {\cdot}y^{d{-}r{+}1{+}{\Delta} '} {\cdot}{\sigma} 
(x_i|_{i{=}1,r} )(y) & y^{{\leq}d} {\cdot}y^{{\Delta} '}  
& y^{{\leq}{\Delta} '{-}1} {\cdot}{\sigma} (x_i|_{i{=}1,r} 
)(y)\end{matrix} \right\| $  \bigskip\\
\\
(Since elements of $y^{{\leq}d} {\cdot}y^{{\Delta} '} $ has in $y$  
the upper degrees ${\leq}d{+}{\Delta} '$, elements of $y^{{\leq}{\Delta} 
'{-}1} {\cdot}{\sigma} (x_i|_{i{=}1,r} )(y)$ has in $y$ the upper 
degrees ${\leq}{\Delta} '{+}r{-}1 {\leq}{\Delta} '{+}d$, 
elements of $y^{{\leq}{\Delta} ''{-}1} {\cdot}y^{d{+}{\Delta} '{+}1} 
$  has  in  $y$ the lower degree  ${\geq}d{+}{\Delta} '{+}1$,  then 
\bigskip

 $[y^{{\leq}{\Delta} ''{-}1} {\cdot}y^{d{+}{\Delta} '{+}1} 
]_{*}.|y^{{\leq}d} {\cdot}y^{{\Delta} '}  =  0|^{{\times}{\Delta} 
''} _{{\times}(d{+}1)} $, $[y^{{\leq}{\Delta} ''{-}1} {\cdot}y^{d{+}{\Delta} 
'{+}1} ]_{*}.|y^{{\leq}{\Delta} '{-}1} {\cdot}{\sigma} (x_i|_{i{=}1,r} 
)(y) =  0|^{{\times}{\Delta} ''} _{{\times}{\Delta} '} 
$. \bigskip
\\
Since elements of  $y^{{\leq}{\Delta} ''{-}1} {\cdot}y^{d{-}r{+}1{+}{\Delta} '} 
{\cdot}{\sigma} (x_i|_{i{=}1,r} )(y)$   has   in   $y$   
 the lower degrees ${\geq}d{-}r{+}1{+}{\Delta} '{\geq}{\Delta} 
'$, elements of $y^{{\leq}d} {\cdot}y^{{\Delta} '} $ has   in   $y$  
the lower   degrees ${\geq}{\Delta} '$, elements of $y^{{\leq}{\Delta} 
'{-}1} $ has in $y$ the upper degrees ${\leq}{\Delta} '{-}1$, 
then \bigskip

 $[y^{{\leq}{\Delta} '{-}1} ]_{*}.|y^{{\leq}{\Delta} 
''{-}1} {\cdot}y^{d{-}r{+}1{+}{\Delta} '} {\cdot}{\sigma} 
(x_i|_{i{=}1,r} )(y) =  0|^{{\times}{\Delta} '} _{{\times}{\Delta} 
''} $, $[y^{{\leq}{\Delta} '{-}1} ]_{*}.|y^{{\leq}d} {\cdot}y^{{\Delta} 
'}  =  0|^{{\times}{\Delta} '} _{{\times}(d{+}1)} $. ) \bigskip

 \qquad $=  {\det}\left\| \begin{matrix}[y^{{\leq}{\Delta} 
''{-}1} {\cdot}y^{d{+}{\Delta} '{+}1} ]_{*}\h & .| & y^{{\leq}{\Delta} 
''{-}1} {\cdot}y^{d{-}r{+}1{+}{\Delta} '} {\cdot}{\sigma} 
(x_i|_{i{=}1,r} )(y) & 0|^{{\times}{\Delta} ''} _{{\times}(d{+}1)} 
\h & 0|^{{\times}{\Delta} ''} _{{\times}{\Delta} '} \h 
\cr
[y^{{\leq}d} {\cdot}y^{{\Delta} '} ]_{*} \h & .| & y^{{\leq}{\Delta} 
''{-}1} {\cdot}y^{d{-}r{+}1{+}{\Delta} '} {\cdot}{\sigma} 
(x_i|_{i{=}1,r} )(y) & y^{{\leq}d} {\cdot}y^{{\Delta} '}  
& y^{{\leq}{\Delta} '{-}1} {\cdot}{\sigma} (x_i|_{i{=}1,r} 
)(y) \cr
[y^{{\leq}{\Delta} '{-}1} ]_{*} \h & .| & 0|^{{\times}{\Delta} 
'} _{{\times}{\Delta} ''} \h & 0|^{{\times}{\Delta} '} 
_{{\times}(d{+}1)} \h & y^{{\leq}{\Delta} '{-}1} {\cdot}{\sigma} 
(x_i|_{i{=}1,r} )(y)\end{matrix} \right\| $ \bigskip

 \qquad $=   {\det}\left\| \begin{matrix}[y^{{\leq}{\Delta} 
''{-}1} {\cdot}y^{d{+}{\Delta} '{+}1} ]_{*}.\h & .| & y^{{\leq}{\Delta} 
''{-}1} {\cdot}y^{d{-}r{+}1{+}{\Delta} '} {\cdot}{\sigma} 
(x_i|_{i{=}1,r} )(y) & 0|^{{\times}{\Delta} ''} _{{\times}(d{+}1)} 
\h & 0|^{{\times}{\Delta} ''} _{{\times}{\Delta} '} \h 
\cr
[y^{{\leq}d} {\cdot}y^{{\Delta} '} ]_{*}.\h & .| & 0|^{{\times}(d{+}1)} 
_{{\times}{\Delta} ''} \h & y^{{\leq}d} {\cdot}y^{{\Delta} 
'}  & 0|^{{\times}(d{+}1)} _{{\times}{\Delta} '} \h \cr
[y^{{\leq}{\Delta} '{-}1} ]_{*}.\h & .| & 0|^{{\times}{\Delta} 
'} _{{\times}{\Delta} ''} \h & 0|^{{\times}{\Delta} '} 
_{{\times}(d{+}1)} \h & y^{{\leq}{\Delta} '{-}1} {\cdot}{\sigma} 
(x_i|_{i{=}1,r} )(y)\end{matrix} \right\| $ \bigskip

 \qquad $=  {\det}\left\| \begin{matrix}[y^{{\leq}{\Delta} 
''{-}1} ]_{*} \h & .| & y^{{\leq}{\Delta} ''{-}1} {\cdot}y^{{-}r} 
{\cdot}{\sigma} (x_i|_{i{=}1,r} )(y) & 0|^{{\times}{\Delta} 
''} _{{\times}(d{+}1)} \h & 0|^{{\times}{\Delta} ''} _{{\times}{\Delta} 
'} \h \cr
[y^{{\leq}d} ]_{*} \h & .| & 0|^{{\times}(d{+}1)} _{{\times}{\Delta} 
''} \h & y^{{\leq}d} {\cdot}y^{{\Delta} '}  & 0|^{{\times}(d{+}1)} 
_{{\times}{\Delta} '} \h \cr
[y^{{\leq}{\Delta} '{-}1} ]_{*} \h & .| & 0|^{{\times}{\Delta} 
'} _{{\times}{\Delta} ''} \h & 0|^{{\times}{\Delta} '} 
_{{\times}(d{+}1)} \h & y^{{\leq}{\Delta} '{-}1} {\cdot}{\sigma} 
(x_i|_{i{=}1,r} )(y)\end{matrix} \right\| $  \bigskip

 \qquad $=  {\det}\left\| \begin{matrix}[y^{{\leq}{\Delta} 
''{-}1} ]_{*}.| & y^{{\leq}{\Delta} ''{-}1} {\cdot}y^{{-}r} 
{\cdot}{\sigma} (x_i|_{i{=}1,r} )(y)\end{matrix} \right\| $  
\medskip

 \qquad \qquad \qquad ${\cdot}{\det}\left\| \begin{matrix}[y^{{\leq}d} 
]_{*}.| & y^{{\leq}d} \end{matrix} \right\| $  \medskip

 \qquad \qquad \qquad \qquad ${\cdot}{\det}\left\| \begin{matrix}[y^{{\leq}{\Delta} 
'{-}1} ]_{*}.| & y^{{\leq}{\Delta} '{-}1} {\cdot}{\sigma} 
(x_i|_{i{=}1,r} )(y)\end{matrix} \right\| $  \bigskip\\
(The upper degree of $y^{\alpha} {\cdot}y^{{-}r} {\cdot}{\sigma} 
(x_i|_{i{=}1,r} )(y)$ in $y$ is ${\leq}{\alpha} {-}r{+}r 
=   {\alpha} $,  then  if ${\beta} {>}{\alpha} $ there holds 
\medskip

 $[y^{\beta} ]_{*}.y^{\alpha} {\cdot}y^{{-}r} {\cdot}{\sigma} 
(x_i|_{i{=}1,r} )(y)  =    0$, \medskip\\
if ${\beta} {=}{\alpha} $ there holds \medskip

 $[y^{\beta} ]_{*}.y^{\alpha} {\cdot}y^{{-}r} {\cdot}{\sigma} 
(x_i|_{i{=}1,r} )(y) =  [y^r]_{*}.{\sigma} (x_i|_{i{=}1,r} )(y) 
=  {\sigma} _r(x_i|_{i{=}1,r} )$. \medskip\\
Hence, the matrix of the first determinant is triangular 
with  elements ${\sigma} _r(x_i|_{i{=}1,r} )$ on  diagonal. \medskip
\\
The lower degree of 
$y^{\alpha} {\cdot}{\sigma} (x_i|_{i{=}1,r} )(y)$ in $y$ is ${\geq}{\alpha} $, 
then if ${\beta} {<}{\alpha} $ there holds \medskip

 $[y^{\beta} ]_{*}.y^{\alpha} {\cdot}{\sigma} (x_i|_{i{=}1,r} 
)(y) =  0$, \medskip\\
if ${\beta} {=}{\alpha} $ there holds \medskip

 $[y^{\beta} ]_{*}.y^{\alpha} {\cdot}{\sigma} (x_i|_{i{=}1,r} 
)(y) =  [y^0]_{*}.{\sigma} (x_i|_{i{=}1,r} )(y) =  {\sigma} _0(x_i|_{i{=}1,r} 
)(y)$. \medskip\\
Hence, the matrix of the third determinant is triangular 
with elements ${\sigma} _r(x_i|_{i{=}1,r} )$ on  diagonal. \medskip
\\
There holds $[y^{\beta} ]_{*}.y^{\alpha}  =  ?({\alpha} {=}{\beta} 
)$.  Hence,  matrix of  the second  determinant is the unit matrix.) \bigskip

 \qquad $=  {\sigma} _0(x_i|_{i{=}1,r} )^{{\Delta} ''} {\cdot}1^{d{+}1} 
{\cdot}{\sigma} _r(x_i|_{i{=}1,r} )^{{\Delta} '}  =  {\sigma} 
_0(x_i|_{i{=}1,r} )^{{\Delta} ''} {\cdot}{\sigma} _r(x_i|_{i{=}1,r} 
)^{{\Delta} '} $. \normalsize\\
\\
{\bf Convention 2.1.} In proofs of theorems 2.1 and 2.3 by 
${\buildrel{ 1}\over  = }$  we   shall  mean transition 2.1, 
by ${\buildrel{ 2}\over  = }$ we  shall mean transition 2.2, 
by ${\buildrel{ 3}\over  = }$  we   shall  mean transition 2.3.
\hfil\vfill\eject

\noindent 
{\bf Lemma 2.1.} {\it Let $r{\geq}1$. \medskip\\
1. $\bigwedge (x^{{\leq}d} |_{{\times}{\chi} } ) \in  \bigwedge 
 ^r({\bf R}[x])^{{\leq}d} $ for ${\chi} {\in}{\bf C}^r_{d{+}1} 
$,  elements  $(\bigwedge (x^{{\leq}d} |_{{\times}{\chi} } )|{\chi} 
{\in}{\bf C}^r_{d{+}1} )$  are lineary  independent,  $\bigwedge 
 ^r({\bf R}[x])^{{\leq}d} $  is a  free  module  with 
  a basis  $(\bigwedge (x^{{\leq}d} |_{{\times}{\chi} } )|{\chi} 
{\in}{\bf C}^r_{d{+}1} )$, where $d{\geq}r{-}1$. \medskip\\
2. $\bigwedge (x^{{\leq}r{-}1} ) \in  \bigwedge  ^r({\bf R}[x])^{{\leq}r{-}1} 
$, $\bigwedge  ^r({\bf R}[x])^{{\leq}r{-}1} $ is an one generated 
 free module with a basis  $(\bigwedge (x^{{\leq}r{-}1} ))$.}\medskip\\
{\bf Proof 1.} Since ${\chi} {\in}{\bf C}^r_{d{+}1} 
$, then $\sum {\chi}  =  r$.  Since  elements  in  $(x^{{\leq}d} 
|_{{\times}{\chi} } )$ belongs to ${\bf R}[x]^{{\leq}d} $, the  
number of their is equal to $\sum {\chi} ' =  r$, then by of proposition  1.16 
$\bigwedge (x^{{\leq}d} |_{{\times}{\chi} '} ) \in 
 \mathop{\sf A}^r({\bf R}[x])^{{\leq}d}  =  \bigwedge  ^r({\bf R}[x])^{{\leq}d} $. \medskip
\\
Since  ${\bf R}[x]^{{\leq}d} $   is a  free   module   with   a basis    
$(x^{{\leq}d} )$,  then $(\bigwedge (x^{{\leq}d} 
|_{{\times}{\chi} } )|{\chi} {\in}{\bf C}^r_{d{+}1} )$ are 
lineary independent  in  $\bigwedge  ^r({\bf R}[x]^{{\leq}d} 
)$,  $\bigwedge  ^r({\bf R}[x]^{{\leq}d} )$ is a free 
module with the basis  
$(\bigwedge (x^{{\leq}d} |_{{\times}{\chi} } )|{\chi} {\in}{\bf C}^r_{d{+}1} )$.  
Since  by of 3 of proposition 1.28 there is an isomorphism of modules  
$\bigwedge  ^r({\bf R}[x])^{{\leq}d}   \simeq  \bigwedge  ^r({\bf R}[x]^{{\leq}d} )$,  
then elements 
 $(\bigwedge (x^{{\leq}d} |_{{\times}{\chi} } )|{\chi} {\in}{\bf 
C}^r_{d{+}1} )$  
are  lineary  independent  in  $\bigwedge  ^r({\bf R}[x])^{{\leq}d} $, 
$\bigwedge  ^r({\bf R}[x])^{{\leq}d} $ 
is a free module with a basis  $(\bigwedge (x^{{\leq}d} 
|_{{\times}{\chi} } )|{\chi} {\in}{\bf C}^r_{d{+}1} )$. \medskip
\\
{\bf Proof 2.} Substituting into 1 of lemma $d =  r{-}1$ we have 
${\chi}  \in  {\bf C}^r_{d{+}1}  =  {\bf C}^r_r$,  then  ${\bf 
C}^r_r$ consists only of one maps $[1,r] \rightarrow 
 \{0,1\}$ such that ${\chi} _i{=}1$ for  all $i{=}1,r$, for 
such ${\chi} $ there holds $\bigwedge (x^{{\leq}d} |_{\chi} 
) =  \bigwedge (x^{{\leq}r{-}1} )$. Then from 1 of lemma it follows 
$\bigwedge (x^{{\leq}r{-}1} ) \in  \bigwedge  ^r({\bf R}[x])^{{\leq}r{-}1} 
$, $\bigwedge  ^r({\bf R}[x])^{{\leq}r{-}1}  =  \bigwedge  ^r({\bf 
R}[x])^{{\leq}d} $ is an one generated free module 
with the basis  $(\bigwedge (x^{{\leq}r{-}1} )) =  (\bigwedge (x^{{\leq}d} 
|_{{\times}{\chi} } )|{\chi} {\in}{\bf C}^r_{d{+}1} )$.\medskip 
\\ 
{\bf Lemma 2.2.} {\it Let $r{\geq}1$, $d{=}r{-}1$, 
${\Delta} {\geq}0$. Let \bigskip \footnotesize

 $S(x_i|_{i{=}1,r} ) =  \sum\limits_{q{\in}Q} {\det}\left\| 
\begin{matrix}[y^{{\leq}d{+}{\Delta} } ]_{*}.| & y^{{\leq}{\Delta} 
{-}1} {\cdot}{\sigma} (x_i|_{i{=}1,r} )(y) & h^q_{{\leq}r} (y)\end{matrix} 
\right\| $  \bigskip \normalsize
\\
for some $h^q_i(x) \in  {\bf R}[x]^{{\leq}d{+}{\Delta} } $ 
for $i{=}1,r$, for $q{\in}Q$. Then \bigskip
\\
1) $S(x_i|_{i{=}1,r} )$ is polynomially expressed via  
$({\sigma} _p(x_i|_{i{=}1,r} )|_{p{=}1,r} )$  
polynomial of the degree ${\leq}{\Delta} $;\medskip
\\
2) $S(x_i|_{i{=}1,r} ){\in}(\mathop{\sf ts}^r({\bf R}[x_i|_{i{=}1,r} ])^{{\leq}1} )^{\Delta} $; \medskip
\\
3) $S(x_i|_{i{=}1,r} ){\in}\mathop{\sf ts}^r({\bf R}[x_i|_{i{=}1,r} ])^{{\leq}{\Delta} } $. } \bigskip
\\
{\bf Proof.} For  $q{\in}Q$ put \bigskip 
\footnotesize

 $S^q(x_i|_{i{=}1,r} ) =  {\det}\left\| \begin{matrix}[y^{{\leq}d{+}{\Delta} 
} ]_{*}.| & y^{{\leq}{\Delta} {-}1} {\cdot}{\sigma} (x_i|_{i{=}1,r} 
)(y) & h_{{\leq}r} (y)\end{matrix} \right\| $. \normalsize \bigskip
\\
Then \bigskip \footnotesize

 $S(x_i|_{i{=}1,r} ) =  \sum\limits_{q{\in}Q} S^q(x_i|_{i{=}1,r} 
)$. \normalsize \bigskip
\\
Columns  of determinant  are   
$[y^{{\leq}d{+}{\Delta} } ]_{*}.y^{{\leq}{\Delta} {-}1} 
{\cdot}{\sigma} (x_i|_{i{=}1,r} )(y)$,   
elements which is elements in $({\sigma} _p(x_i|_{i{=}1,r} )|_{p{=}0,r} )$, 
or zeroes, the number of such columns is  ${\Delta} $.  
Also  columns  of determinant are  columns   
$[y^{{\leq}d{+}{\Delta} } ]_{*}.h^q_{{\leq}r} (y)$, 
elements  which  are  elements in ${\bf R}$.  
Hence, the determinant is polynomially expressed via elements  
in  $({\sigma} _p(x_i|_{i{=}1,r} )|_{p{=}0,r} )$  
by a homogeneouse polynomial of the degree  ${=}{\Delta} $,  
since it polylineary  depends from own  columns.  
Then 
$S^q(x_i|_{i{=}1,r} )$  
is polynomially  expressed via   
$({\sigma}_p(x_i|_{i{=}1,r} )|_{p{=}0,r} )$   
by a homogeneouse polynomial of the degree ${=}{\Delta} $, 
then $S(x_i|_{i{=}1,r} )$ as the sum of $S^q(x_i|_{i{=}1,r} )$ 
over $q{\in}Q$ is polynomially expressed via 
$({\sigma} _p(x_i|_{i{=}1,r} )|_{p{=}0,r} )$ 
by a homogeneouse polynomial of degree ${=}{\Delta} $. 
\medskip
\\
Since ${\sigma} _0(x_i|_{i{=}1,r} ){=}1(x_i|_{i{=}1,r} )$, 
then  $S(x_i|_{i{=}1,r} )$ is  polynomially  expressed via elements 
$({\sigma} _p(x_i|_{i{=}1,r} )|_{p{=}1,r} )$ by a polynomial of the degree 
${\leq}{\Delta} $. 
Since  elements in  $({\sigma} _p(x_i|_{i{=}1,r} )|_{p{=}0,r} )$  
belongs to  $\mathop{\sf ts}^r({\bf R}[x_i|_{i{=}1,r} ])^{{\leq}1} $,  
and $S(x_i|_{i{=}1,r} )$ is expressed via  elements  
$({\sigma} _p(x_i|_{i{=}1,r} )|_{p{=}0,r} )$  
by a homogeneouse  polynomial of the degree ${=}{\Delta} $, 
then $S(x_i|_{i{=}1,r} ) \in  
(\mathop{\sf ts}^r({\bf R}[x_i|_{i{=}1,r} ])^{{\leq}1} )^{\Delta}$. 
Since by of proposition 1.25 \bigskip \footnotesize

 $(\mathop{\sf ts}^r({\bf R}[x_i|_{i{=}1,r} ])^{{\leq}1} )^{\Delta} 
 \subseteq  \mathop{\sf ts}^r({\bf R}[x_i|_{i{=}1,r} ])^{{\leq}{\Delta} 
} $, \normalsize \bigskip\\
then $S(x_i|_{i{=}1,r} ) \in  \mathop{\sf ts}^r({\bf R}[x_i|_{i{=}1,r} 
])^{{\leq}{\Delta} } $.  
\hfil\vfill\eject
 
 \noindent {\bf Lemma 2.3.} {\it Let $r{\geq}1$, $d{=}r{-}1$, 
${\Delta} {\geq}0$. 
\bigskip\\
1. If \bigskip \footnotesize

 $S(x_i|_{i{=}1,r} ) \in  \mathop{\sf ts}^r({\bf R}[x_i|_{i{=}1,r} 
])^{{\leq}{\Delta} } $, \normalsize \bigskip\\
then \bigskip \footnotesize

 $S({\bf x}_i|_{i{=}1,r} ){\cdot}\bigwedge (x^{{\leq}r{-}1} 
) \in  \bigwedge  ^r({\bf R}[x])^{{\leq}d{+}{\Delta} } $, 
\normalsize \bigskip\\
and for some $h^q_i(x) \in  {\bf R}[x]^{{\leq}d{+}{\Delta} 
} $ for $i{=}1,r$, for $q{\in}Q$ there holds \bigskip \footnotesize

 $S({\bf x}_i|_{i{=}1,r} ){\cdot}\bigwedge (x^{{\leq}r{-}1} 
) =  \sum\limits_{q{\in}Q} \bigwedge (h^q_{{\leq}r} (x))$. \normalsize 
\bigskip\\
2. Let $h^q_i(x) \in  {\bf R}[x]^{{\leq}d{+}{\Delta} } 
$ for $i{=}1,r$, for $q{\in}Q$, let \bigskip \footnotesize

 $S(x_i|_{i{=}1,r} ) \in  \mathop{\sf ts}^r({\bf R}[x_i|_{i{=}1,r} 
])$, \normalsize \bigskip\\
then equalities \bigskip \footnotesize

 $S({\bf x}_i|_{i{=}1,r} ){\cdot}\bigwedge (x^{{\leq}r{-}1} 
) =  \sum\limits_{q{\in}Q} \bigwedge (h^q_{{\leq}r} (x))$, \bigskip

 $S(x_i|_{i{=}1,r} ){\cdot}{\det}\left\| x^{{\leq}r{-}1} _i|^{i{=}1,r} 
\right\|  =  \sum\limits_{q{\in}Q} {\det}\left\| h^q_{{\leq}r} 
(x_i)|^{i{=}1,r} \right\| $, \normalsize \bigskip\\
are equivalent. \bigskip\\
3. Let \bigskip \footnotesize

 $S(x_i|_{i{=}1,r} ) \in  \mathop{\sf ts}^r({\bf R}[x_i|_{i{=}1,r} 
])$, \normalsize \bigskip\\
then equalities \bigskip \footnotesize

 $S({\bf x}_i|_{i{=}1,r} ){\cdot}\bigwedge (x^{{\leq}r{-}1} 
) =  0_{\bigwedge  ^r({\bf R}[x])} $, \bigskip

 $S(x_i|_{i{=}1,r} ){\cdot}{\det}\left\| x^{{\leq}r{-}1} _i|^{i{=}1,r} 
\right\|  =  0(x_i|_{i{=}1,r} )$, \normalsize \bigskip\\
are equivalent. \medskip\\
4. Let \bigskip \footnotesize 

 $S(x_i|_{i{=}1,r} ) \in  \mathop{\sf ts}^r({\bf R}[x_i|_{i{=}1,r} 
])$, \normalsize \bigskip\\
then statements \bigskip \footnotesize

 $S(x_i|_{i{=}1,r} ) \in  \mathop{\sf ts}^r({\bf R}[x_i|_{i{=}1,r} 
])^{{\leq}{\Delta} } $, \bigskip

 $S({\bf x}_i|_{i{=}1,r} ) \in  \mathop{\sf TS}^r({\bf R}[x])^{{\leq}{\Delta} 
} $, \bigskip \normalsize\\
are equivalent. }\bigskip
\\
{\bf Proof 1.} By of 2 of lemma 2.1 $\bigwedge 
(x^{{\leq}r{-}1} ) \in  \bigwedge  ^r({\bf R}[x])^{{\leq}r{-}1} 
$. \medskip\\
Since $S(x_i|_{i{=}1,r} ) \in  \mathop{\sf ts}^r({\bf R}[x_i|_{i{=}1,r} 
])^{{\leq}{\Delta} } $, then by of proposition 1.26 \bigskip 
\footnotesize

 $S({\bf x}_i|_{i{=}1,r} ) =  {\nu} _s(S(x_i|_{i{=}1,r} )) 
\in  {\nu} _s(\mathop{\sf ts}^r({\bf R}[x_i|_{i{=}1,r} ])^{{\leq}{\Delta} 
} ) =  \mathop{\sf TS}^r({\bf R}[x])^{{\leq}{\Delta} } $. 
\bigskip \normalsize\\
Then by of proposition 1.20 \footnotesize \bigskip

 $S({\bf x}_i|_{i{=}1,r} ){\cdot}\bigwedge (x^{{\leq}r{-}1} 
) \in  \mathop{\sf TS}^r({\bf R}[x])^{{\leq}{\Delta} } {\cdot}\mathop{\sf 
A}^r({\bf R}[x])^{{\leq}d}  \subseteq  \mathop{\sf A}^r({\bf 
R}[x])^{{\leq}d{+}{\Delta} }  =  \bigwedge  ^r({\bf R}[x])^{{\leq}d{+}{\Delta} 
} $. \bigskip \normalsize\\
Then by of proposition 1.17 \footnotesize \bigskip

 $S({\bf x}_i|_{i{=}1,r} ){\cdot}\bigwedge (x^{{\leq}r{-}1} 
) =  \sum\limits_{q{\in}Q} \bigwedge (h^q_{{\leq}r} (x))$, \bigskip 
\normalsize\\
for some $h^q_i(x) \in  {\bf R}[x]^{{\leq}d{+}{\Delta} 
} $ for $i{=}1,r$, for $q{\in}Q$. 
\hfil\vfill\eject
\noindent 
{\bf Proof 2,3.} There there holds \footnotesize \bigskip

 $\mathop{\sf atr}(S({\bf x}_i|_{i{=}1,r} ){\cdot}\bigwedge 
(x^{{\leq}r{-}1} )) =  S({\bf x}_i|_{i{=}1,r} ){\cdot}\mathop{\sf 
atr}(\bigwedge (x^{{\leq}r{-}1} ))$  \medskip

 \qquad $=  S({\bf x}_i|_{i{=}1,r} ){\cdot}{\det}\left\| {\bf 
x}^{{\leq}r{-}1} _i|^{i{=}1,r} \right\|  =  {\nu} _a(S(x_i|_{i{=}1,r} 
)){\cdot}{\nu} _s({\det}\left\| x^{{\leq}r{-}1} _i|^{i{=}1,r} 
\right\| )$  \medskip

 \qquad $=  {\nu} _a(S(x_i|_{i{=}1,r} ){\cdot}{\det}\left\| 
x^{{\leq}r{-}1} _i|^{i{=}1,r} \right\| )$, \normalsize \bigskip
\\
since the first equality  there holds by of proposition 1.13, 
the second equality  holds by of lemma 1.1,  the fourth  equality 
   holds  by of  2) of proposition 1.23. \medskip\\
There holds \footnotesize \bigskip

 $\mathop{\sf atr}(\sum\limits_{q{\in}Q} \bigwedge (h^q_{{\leq}r} 
(x))) =  \sum\limits_{q{\in}Q} \mathop{\sf atr}(\bigwedge (h^q_{{\leq}r} 
(x))) =  \sum\limits_{q{\in}Q} {\det}\left\| h^q_{{\leq}r} ({\bf 
x}_i)|^{i{=}1,r} \right\| $  \medskip

 \qquad $=  {\nu} _a(\sum\limits_{q{\in}Q} {\det}\left\| h^q_{{\leq}r} 
(x_i)|^{i{=}1,r} \right\| )$, \normalsize \bigskip\\
since the second equality  holds by of lemma 1.1. \medskip
\\
There holds \footnotesize \bigskip

 $\mathop{\sf atr}(0_{\mathop{\sf A}^r({\bf R}[x])} ) =  0_{\mathop{\sf 
TA}^r({\bf R}[x])}  =  {\nu} _a(0_{\mathop{\sf ta}^r({\bf R}[x_i|_{i{=}1,r} 
])} ) =  {\nu} _a(0(x_i|_{i{=}1,r} ))$. \normalsize \bigskip
\\
Since by of proposition 1.5 and proposition 1.23 maps 
\bigskip \footnotesize

 $\mathop{\sf atr}^r_{{\bf R}[x]} {:}\ \mathop{\sf A}^r({\bf 
R}[x]) \rightarrow  \mathop{\sf TA}^r({\bf R}[x])$, ${\nu} _a{:}\ 
\mathop{\sf ta}^r({\bf R}[x_i|_{i{=}1,r} ]) \rightarrow  \mathop{\sf 
TA}^r({\bf R}[x])$, \bigskip\\
\normalsize are invertible, then equalities \footnotesize 
\bigskip

 $S({\bf x}_i|_{i{=}1,r} ){\cdot}\bigwedge (x^{{\leq}r{-}1} 
) =  \sum\limits_{q{\in}Q} \bigwedge (h^q_{{\leq}r} (x))$, \bigskip

 $S(x_i|_{i{=}1,r} ){\cdot}{\det}\left\| x^{{\leq}r{-}1} _i|^{i{=}1,r} 
\right\|  =  \sum\limits_{q{\in}Q} {\det}\left\| h^q_{{\leq}r} 
(x_i)|^{i{=}1,r} \right\| $, \normalsize \bigskip
\\
are equivalent, and equalities \footnotesize \bigskip

 $S({\bf x}_i|_{i{=}1,r} ){\cdot}\bigwedge (x^{{\leq}r{-}1} 
) =  0_{\bigwedge  ^r({\bf R}[x])} $, \bigskip

 $S(x_i|_{i{=}1,r} ){\cdot}{\det}\left\| x^{{\leq}r{-}1} _i|^{i{=}1,r} 
\right\|  =  0(x_i|_{i{=}1,r} )$. \normalsize \bigskip\\
are equivalent.\medskip
\\
{\bf Proof 4.} Since \bigskip \footnotesize

 ${\nu} _s(S(x_i|_{i{=}1,r} )) =  S({\bf x}_i|_{i{=}1,r} )$, 
\bigskip \normalsize\\
and by of proposition 1.26 \bigskip \footnotesize

 ${\nu} _s(\mathop{\sf ts}^r({\bf R}[x_i|_{i{=}1,r} ])^{{\leq}{\Delta} 
} ) =  \mathop{\sf TS}^r({\bf R}[x])^{{\leq}{\Delta} } $, 
\bigskip \normalsize\\
then \bigskip \footnotesize

 $S(x_i|_{i{=}1,r} ) \in  \mathop{\sf ts}^r({\bf R}[x_i|_{i{=}1,r} 
])^{{\leq}{\Delta} } $, \bigskip

 $S({\bf x}_i|_{i{=}1,r} ) \in  \mathop{\sf TS}^r({\bf R}[x])^{{\leq}{\Delta} 
} $, \bigskip \normalsize\\
are equivalent. Since there holds the first, then there holds and the second.

\hfil\vfill\eject
\noindent 
 \noindent {\bf Theorem 2.1.} {\it Let $r{\geq}1$, $d{=}r{-}1$, 
${\Delta} {\geq}0$. \bigskip\\
1. Let $h_i(x) \in  {\bf R}[x]^{{\leq}d{+}{\Delta} } $ 
for $i{=}1,r$. Let \bigskip \footnotesize

 $S(x_i|_{i{=}1,r} ) =  {\det}\left\| \begin{matrix}[y^{{\leq}d{+}{\Delta} 
} ]_{*}.| & y^{{\leq}{\Delta} {-}1} {\cdot}{\sigma} (x_i|_{i{=}1,r} 
)(y) & h_{{\leq}r} (y)\end{matrix} \right\| $, \normalsize \bigskip
\\
then \bigskip \footnotesize

 $S(x_i|_{i{=}1,r} ) \in  \mathop{\sf ts}^r({\bf R}[x_i|_{i{=}1,r} 
])^{{\leq}{\Delta} } $, \bigskip

 ${\det}\left\| h_{{\leq}r} (x_i)|^{i{=}1,r} \right\|  =  S(x_i|_{i{=}1,r} 
){\cdot}{\det}\left\| x^{{\leq}r{-}1} _i|^{i{=}1,r} \right\| 
$, \normalsize \bigskip\\
 that is equivalent to \bigskip \footnotesize

 $S({\bf x}_i|_{i{=}1,r} ) \in  \mathop{\sf TS}^r({\bf R}[x])^{{\leq}{\Delta} 
} $, \bigskip

 $\bigwedge (h_{{\leq}r} (x)) =  S({\bf x}_i|_{i{=}1,r} ){\cdot}\bigwedge 
(x^{{\leq}r{-}1} )$. \normalsize \bigskip\\
2. There holds \bigskip \footnotesize

 $\bigwedge  ^r({\bf R}[x])^{{\leq}r{-}1{+}{\Delta} }  \subseteq 
 \mathop{\sf TS}^r({\bf R}[x])^{{\leq}{\Delta} } {\cdot}\bigwedge 
 ^r({\bf R}[x])^{{\leq}r{-}1} $, \normalsize \bigskip\\
3.  $\bigwedge  ^r({\bf R}[x])$  is a free  one generated 
 module  over   $\mathop{\sf TS}^r({\bf R}[x])$ by element $\bigwedge 
(x^{{\leq}r{-}1} )$. } \bigskip
\\
{\bf Proof 1.} Let ${\Delta} {=}{\Delta} 
'{+}{\Delta} ''$, where ${\Delta} '{\geq}0$, ${\Delta} 
''{\geq}0$. There holds \footnotesize \bigskip

 ${\sigma} _0(x_i|_{i{=}1,r} )^{{\Delta} ''} {\cdot}{\sigma} 
_r(x_i|_{i{=}1,r} )^{{\Delta} '} {\cdot}S(x_i|_{i{=}1,r} ){\cdot}{\det}\left\| 
x^{{\leq}r{-}1} _i|^{i{=}1,r} \right\| $  \bigskip

 \quad $=  {\sigma} _0(x_i|_{i{=}1,r} )^{{\Delta} ''} {\cdot}{\sigma} 
_r(x_i|_{i{=}1,r} )^{{\Delta} '} {\cdot}{\det}\left\| \begin{matrix}[y^{{\leq}d{+}{\Delta} 
} ]_{*}.| & y^{{\leq}{\Delta} {-}1} {\cdot}{\sigma} (x_i|_{i{=}1,r} 
)(y) & h_{{\leq}r} (y)\end{matrix} \right\| {\cdot}{\det}\left\| 
x^{{\leq}r{-}1} _i|^{i{=}1,r} \right\| $  \bigskip

 \quad $\buildrel{ 1}\over  = ({-}1)^{r{\cdot}{\Delta} '} 
{\cdot}{\det}\left\| \begin{matrix}[y^{{\leq}d{+}{\Delta} 
} ]_{*}.| & y^{{\leq}{\Delta} {-}1} {\cdot}{\sigma} (x_i|_{i{=}1,r} 
)(y) & h_{{\leq}r} (y)\end{matrix} \right\| {\cdot}{\det}\left\| 
x^{{\leq}r{-}1} _i{\cdot}x^{{\Delta} '} _i|^{i{=}1,r} \right\| 
$  \bigskip\\

 \quad $=  ({-}1)^{r{\cdot}{\Delta} '} {\cdot}{\det}\left\| 
\begin{matrix}[y^{{\leq}d{+}{\Delta} } ]_{*}.| & y^{{\leq}{\Delta} 
{-}1} {\cdot}{\sigma} (x_i|_{i{=}1,r} )(y) & h_{{\leq}r} (y) 
& y^{{\leq}r{-}1} {\cdot}y^{{\Delta} '} \h \cr
  &                        &          &   \cr
             &                        &          & x^{{\leq}r{-}1} 
_i{\cdot}x^{{\Delta} '} _i|^{i{=}1,r} \end{matrix} \right\| 
$  \bigskip\\

 \quad $\buildrel{ 2}\over  = ({-}1)^{r{\cdot}{\Delta} '} 
{\cdot}{\det}\left\| \begin{matrix}[y^{{\leq}d{+}{\Delta} 
} ]_{*}.| & y^{{\leq}{\Delta} {-}1} {\cdot}{\sigma} (x_i|_{i{=}1,r} 
)(y) & \ \ \, h_{{\leq}r} (y)\h & y^{{\leq}r{-}1} {\cdot}y^{{\Delta} 
'}  \cr
  &                        &                   &   \cr
             &                        & {-}h_{{\leq}r} (x_i)|^{i{=}1,r} 
 &  \end{matrix} \right\| $  \bigskip\\

 \quad $=  ({-}1)^r$  \medskip

 \qquad ${\cdot}{\det}\left\| \begin{matrix}[y^{{\leq}d{+}{\Delta} 
} ]_{*}.| & y^{{\leq}{\Delta} ''{-}1} {\cdot}y^{{\Delta} 
'} {\cdot}{\sigma} (x_i|_{i{=}1,r} )(y) & y^{{\leq}r{-}1} {\cdot}y^{{\Delta} 
'}  & y^{{\leq}{\Delta} '{-}1} {\cdot}{\sigma} (x_i|_{i{=}1,r} 
)(y) & \ \ \, h_{{\leq}r} (y)\h \cr
  &                             &   &                      
   &                \cr
             &                             &   &           
              & {-}h_{{\leq}r} (x_i)|^{i{=}1,r} \end{matrix} 
\right\| $  \bigskip\\

 \quad $=  ({-}1)^r{\cdot}{\det}\left\| \begin{matrix}[y^{{\leq}d{+}{\Delta} 
} ]_{*}.| & y^{{\leq}{\Delta} ''{-}1} {\cdot}y^{{\Delta} 
'} {\cdot}{\sigma} (x_i|_{i{=}1,r} )(y) & y^{{\leq}r{-}1} {\cdot}y^{{\Delta} 
'}  & y^{{\leq}{\Delta} '{-}1} {\cdot}{\sigma} (x_i|_{i{=}1,r} 
)(y)\end{matrix} \right\| $  \medskip

 \qquad \qquad \qquad \qquad \qquad \qquad \qquad \qquad \qquad 
\qquad \qquad \qquad  \qquad  \qquad  \qquad \qquad \qquad \qquad 
${\cdot}{\det}\left\| {-}h_{{\leq}r} (x_i)|^{i{=}1,r} \right\| 
$  \bigskip

 \quad $\buildrel{ 3}\over  = ({-}1)^r{\cdot}{\sigma} _0(x_i|_{i{=}1,r} 
)^{{\Delta} ''} {\cdot}{\sigma} _1(x_i|_{i{=}1,r} )^{{\Delta} 
'} {\cdot}{\det}\left\| {-}h_{{\leq}r} (x_i)|^{i{=}1,r} \right\| 
$  \bigskip

 \quad $=  {\sigma} _0(x_i|_{i{=}1,r} )^{{\Delta} ''} {\cdot}{\sigma} 
_1(x_i|_{i{=}1,r} )^{{\Delta} '} {\cdot}{\det}\left\| h_{{\leq}r} 
(x_i)|^{i{=}1,r} \right\| $. \normalsize \bigskip\\
Dividing both parts of the obtained equality by the common multiple 
we obtain \footnotesize \bigskip

 $S(x_i|_{i{=}1,r} ){\cdot}{\det}\left\| x^{{\leq}r{-}1} _i|^{i{=}1,r} 
\right\|  =  {\det}\left\| h_{{\leq}r} (x_i)|^{i{=}1,r} \right\| 
$. \normalsize \bigskip
\hfil\vfill\eject
\noindent 
By of 3) of lemma 2.2 \bigskip \footnotesize

 $\mathop{\sf ts}^r({\bf R}[x_i|_{i{=}1,r} ])^{{\leq}{\Delta} 
} $. \normalsize \bigskip\\
By of 4) of lemma 2.3 statements \bigskip \footnotesize

 $S(x_i|_{i{=}1,r} ) \in  \mathop{\sf ts}^r({\bf R}[x_i|_{i{=}1,r} 
])^{{\leq}{\Delta} } $, \bigskip

 $S({\bf x}_i|_{i{=}1,r} ) \in  \mathop{\sf TS}^r({\bf R}[x])^{{\leq}{\Delta} 
} $, \bigskip \normalsize\\
are equivalent. Since there holds  the first  statement, 
 then  there holds and the second statement. \medskip\\
By of 2 of lemma 2.3 equality \bigskip \footnotesize

 $\bigwedge (h_{{\leq}r} (x)) =  S({\bf x}_i|_{i{=}1,r} ){\cdot}\bigwedge 
(x^{{\leq}r{-}1} )$, \bigskip

 ${\det}\left\| h_{{\leq}r} (x_i)|^{i{=}1,r} \right\|  =  S(x_i|_{i{=}1,r} 
){\cdot}{\det}\left\| x^{{\leq}r{-}1} _i|^{i{=}1,r} \right\| 
$,  \bigskip \normalsize\\
are equivalent. Since there holds the second equality, 
then there holds the first equality.
\\
{\bf Proof  2.}  By of  2  of lemma  2.1 
 $\bigwedge  ^r({\bf R}[x])^{{\leq}r{-}1} {\in}\bigwedge (x^{{\leq}r{-}1} )$.  
Since $h_i(x){\in}{\bf R}[x]^{{\leq}d{+}{\Delta} } $ 
for $i{=}1,r$, then by of proposition 1.16 $\bigwedge (h_{{\leq}r} 
(x)) \in  \mathop{\sf A}^r(x)^{{\leq}d{+}{\Delta} } $. By of 
proposition  1.17  $\bigwedge  ^r({\bf R}[x])^{{\leq}d{+}{\Delta} } $  
is lineary  generated by  elements   of the form $\bigwedge (h_{{\leq}r} 
(x))$, where $h_i(x) \in  {\bf R}[x]^{{\leq}d{+}{\Delta} } $ 
for $i{=}1,r$. By  proposition  1.26  \bigskip \footnotesize

 ${\nu} _s(\mathop{\sf ts}^r({\bf R}[x_i|_{i{=}1,r} ])^{{\leq}{\Delta} 
} ) =  \mathop{\sf TS}^r({\bf R}[x])^{{\leq}{\Delta} } $, 
\bigskip \normalsize
\\
then any element in $\mathop{\sf TS}^r({\bf R}[x])^{{\leq}{\Delta} 
} $ is of the form $S({\bf x}_i|_{i{=}1,r} )  =   {\nu} _s(S(x_i|_{i{=}1,r} 
))$, where $S(x_i|_{i{=}1,r} ) \in  \mathop{\sf ts}^r({\bf R}[x_i|_{i{=}1,r} 
])^{{\leq}{\Delta} } $. Then from equality \bigskip \normalsize

 $\bigwedge (h_{{\leq}r} (x)) =  S({\bf x}_i|_{i{=}1,r} ){\cdot}\bigwedge 
(x^{{\leq}r{-}1} )$, \normalsize \bigskip\\
it follows \bigskip \footnotesize

 $\bigwedge  ^r({\bf R}[x])^{{\leq}r{-}1{+}{\Delta} }  \subseteq 
 \mathop{\sf TS}^r({\bf R}[x])^{{\leq}{\Delta} } {\cdot}\bigwedge 
 ^r({\bf R}[x])^{{\leq}r{-}1} $. \normalsize \bigskip
\\
{\bf Proof 3.} Since by of 2 of lemma 2.1 $\bigwedge 
 ^r({\bf R}[x])^{{\leq}r{-}1} $ is lineary generated by element 
$\bigwedge (x^{{\leq}r{-}1} )$, then by of 2  of theorem  $\bigwedge 
 ^r({\bf R}[x])$  is an  one generated module over 
$\mathop{\sf TS}^r({\bf R}[x])$ element $\bigwedge (x^{{\leq}r{-}1} )$. 
By of proposition  1.26  \bigskip \footnotesize

 ${\nu} _s(\mathop{\sf ts}^r({\bf R}[x_i|_{i{=}1,r} ])^{{\leq}{\Delta} 
} ) =  \mathop{\sf TS}^r({\bf R}[x])^{{\leq}{\Delta} } $, 
\bigskip \normalsize
\\
then any element in $\mathop{\sf TS}^r({\bf R}[x])^{{\leq}{\Delta} } $ 
is of the form $S({\bf x}_i|_{i{=}1,r} )  =   {\nu} _s(S(x_i|_{i{=}1,r} ))$, 
where \bigskip \footnotesize

 $S(x_i|_{i{=}1,r} ) \in  \mathop{\sf ts}^r({\bf R}[x_i|_{i{=}1,r} ])$. \normalsize \bigskip
\\
Let there holds \bigskip \footnotesize

 $S({\bf x}_i|_{i{=}1,r} ){\cdot}\bigwedge (x^{{\leq}r{-}1} 
) =  0_{\bigwedge  ^r({\bf R}[x])} $, \normalsize \bigskip
\\
By of 3 of lemma 2.3 equalities \bigskip \footnotesize

 $S({\bf x}_i|_{i{=}1,r} ){\cdot}\bigwedge (x^{{\leq}r{-}1} 
) =  0_{\bigwedge  ^r({\bf R}[x])} $, \bigskip

 $S(x_i|_{i{=}1,r} ){\cdot}{\det}\left\| x^{{\leq}r{-}1} _i|^{i{=}1,r} 
\right\|  =  0(x_i|_{i{=}1,r} )$, \normalsize \bigskip\\
are equivalent, and since there holds  the first  equality, 
 then  there holds the second equality. Then \bigskip \footnotesize

 $S(x_i|_{i{=}1,r} ) =  0(x_i|_{i{=}1,r} )$, \normalsize \bigskip
\\
hence, \bigskip \footnotesize

 $S({\bf x}_i|_{i{=}1,r} ) =  0({\bf x}_i|_{i{=}1,r} ) =  0_{\mathop{\sf 
TS}^r({\bf R}[x])} $. \normalsize \bigskip\\
Hence,  $\bigwedge  ^r({\bf R}[x])$  
is a free  one generated  module  over 
$\mathop{\sf TS}^r({\bf R}[x])$ 
by element $\bigwedge (x^{{\leq}r{-}1} )$.\\
\hfil\vfill\eject
 
 \noindent {\bf Theorem 2.2.} {\it Let $r{\geq}1$, $d{=}r{-}1$, 
${\Delta} {\geq}0$. 
\bigskip\\
1. Let $h^q_i(x) \in  {\bf R}[x]^{{\leq}d{+}{\Delta} } 
$ for $i{=}1,r$, for $q{\in}Q$; let \bigskip \footnotesize

 $S(x_i|_{i{=}1,r} ) \in  \mathop{\sf ts}^r({\bf R}[x_i|_{i{=}1,r} 
])$. \normalsize \bigskip\\
If \bigskip \footnotesize

 $S({\bf x}_i|_{i{=}1,r} ){\cdot}\bigwedge (x^{{\leq}r{-}1} 
) =  \sum\limits_{q{\in}Q} \bigwedge (h^q_{{\leq}r} (x))$, \normalsize 
\bigskip\\
 that is equivalent to \bigskip \footnotesize

 $S(x_i|_{i{=}1,r} ){\cdot}{\det}\left\| x^{{\leq}r{-}1} _i|^{i{=}1,r} 
\right\|  =  \sum\limits_{q{\in}Q} {\det}\left\| h^q_{{\leq}r} 
(x_i)|^{i{=}1,r} \right\| $, \normalsize \bigskip\\
then \bigskip \footnotesize

 $S(x_i|_{i{=}1,r} ) =  \sum\limits_{q{\in}Q} {\det}\left\| 
\begin{matrix}[y^{{\leq}d{+}{\Delta} } ]_{*}.| & y^{{\leq}{\Delta} 
{-}1} {\cdot}{\sigma} (x_i|_{i{=}1,r} )(y) & h^q_{{\leq}r} (y)\end{matrix} 
\right\| $, \bigskip

 $S(x_i|_{i{=}1,r} ) \in  \mathop{\sf ts}^r({\bf R}[x_i|_{i{=}1,r} 
])^{{\leq}{\Delta} } $. \normalsize \bigskip 
\\
2. There holds \bigskip \footnotesize

 $\mathop{\sf TS}^r({\bf R}[x])^{{\leq}{\Delta} }  \subseteq 
 (\mathop{\sf TS}^r({\bf R}[x])^{{\leq}1} )^{\Delta} $, \normalsize 
\bigskip
\\
this is equivalent to \bigskip \footnotesize

 $\mathop{\sf ts}^r({\bf R}[x_i|_{i{=}1,r} ])^{{\leq}{\Delta} 
}  \subseteq  (\mathop{\sf ts}^r({\bf R}[x_i|_{i{=}1,r} ])^{{\leq}1} 
)^{\Delta} $. \normalsize \bigskip
\\
3. $\mathop{\sf ts}^r({\bf R}[x_i|_{i{=}1,r} ])^{{\leq}{\Delta} } $ 
is polynomially generated by elements 
$({\sigma} _p(x_i|_{i{=}1,r} )|_{p{=}1,r} )$ 
by polynomials of the degree ${\leq}{\Delta} $. } 
\bigskip
\\
{\bf Proof 1.} By of 2 of lemma 2.3 equalities \bigskip 
\footnotesize

 $S({\bf x}_i|_{i{=}1,r} ){\cdot}\bigwedge (x^{{\leq}r{-}1} 
) =  \sum\limits_{q{\in}Q} \bigwedge (h^q_{{\leq}r} (x))$, \bigskip

 $S(x_i|_{i{=}1,r} ){\cdot}{\det}\left\| x^{{\leq}r{-}1} _i|^{i{=}1,r} 
\right\|  =  \sum\limits_{q{\in}Q} {\det}\left\| h^q_{{\leq}r} 
(x_i)|^{i{=}1,r} \right\| $, \normalsize \bigskip\\
are equivalent. Since there holds the first equality, 
then there holds the second equality. \bigskip\\
Let $q{\in}Q$. Put \bigskip \footnotesize

 $S^q(x_i|_{i{=}1,r} ) =  {\det}\left\| \begin{matrix}[y^{{\leq}d{+}{\Delta} 
} ]_{*}.| & y^{{\leq}{\Delta} {-}1} {\cdot}{\sigma} (x_i|_{i{=}1,r} 
)(y) & h^q_{{\leq}r} (y)\end{matrix} \right\| $. \normalsize 
\bigskip\\
Since $h^q_i(x) \in  {\bf R}[x]^{{\leq}d{+}{\Delta} } 
$  for  $i{=}1,r$,  $r{\geq}1$,  $d{=}r{-}1$,  ${\Delta} {\geq}0$, 
then for  $S^q(x_i|_{i{=}1,r} )$ it is satisfies the conditions of theorem 
 2.1.  Then  by   of theorem   2.1   \bigskip \footnotesize

 $S^q(x_i|_{i{=}1,r} ){\cdot}{\det}\left\| x^{{\leq}r{-}1} 
_i|^{i{=}1,r} \right\|  =  {\det}\left\| h^q_{{\leq}r} (x_i)|^{i{=}1,r} 
\right\| $. \normalsize \bigskip\\
Then \bigskip \footnotesize

 $S(x_i|_{i{=}1,r} ){\cdot}{\det}\left\| x^{{\leq}r{-}1} _i|^{i{=}1,r} 
\right\|  =  \sum\limits_{q{\in}Q} {\det}\left\| h^q_{{\leq}r} 
(x_i)|^{i{=}1,r} \right\|  =  (\sum\limits_{q{\in}Q} S^q(x_i|_{i{=}1,r} 
)){\cdot}{\det}\left\| x^{{\leq}r{-}1} _i|^{i{=}1,r} \right\| 
$. \normalsize \bigskip\\
Dividing both parts of the obtained equality by the common multiple 
we obtain \bigskip \footnotesize

 $S(x_i|_{i{=}1,r} ) =  \sum\limits_{q{\in}Q} S^q(x_i|_{i{=}1,r} 
) =  \sum\limits_{q{\in}Q} {\det}\left\| \begin{matrix}[y^{{\leq}d{+}{\Delta} 
} ]_{*}.| & y^{{\leq}{\Delta} {-}1} {\cdot}{\sigma} (x_i|_{i{=}1,r} 
)(y) & h^q_{{\leq}r} (y)\end{matrix} \right\| $. \normalsize 
\bigskip\\
Then from 3 of lemma 2.2 it follows \bigskip \footnotesize

 $S(x_i|_{i{=}1,r} ) \in  \mathop{\sf ts}^r({\bf R}[x_i|_{i{=}1,r} 
])^{{\leq}{\Delta} } $. \normalsize \bigskip
\hfil\vfill\eject
 \noindent 
{\bf Proof 2,3.} Let there holds \bigskip \footnotesize

 $S(x_i|_{i{=}1,r} ) \in  \mathop{\sf ts}^r({\bf R}[x_i|_{i{=}1,r} 
])^{{\leq}{\Delta} } $. \bigskip \normalsize
\\
Then by of 1 of lemma 2.3 for some 
$h^q_i(x) \in  {\bf R}[x]^{{\leq}d{+}{\Delta} } $  
for  $i{=}1,r$, for  $q{\in}Q$ 
there holds \bigskip\footnotesize

 $S({\bf x}_i|_{i{=}1,r} ){\cdot}\bigwedge (x^{{\leq}r{-}1} ) 
=  \sum\limits_{q{\in}Q} \bigwedge (h^q_{{\leq}r} (x))$. \normalsize \bigskip\\
Then by of 1 of theorem \bigskip \footnotesize

 $S(x_i|_{i{=}1,r} ) =  \sum\limits_{q{\in}Q} {\det}\left\| 
\begin{matrix}[y^{{\leq}d{+}{\Delta} } ]_{*}.| & y^{{\leq}{\Delta} 
{-}1} {\cdot}{\sigma} (x_i|_{i{=}1,r} )(y) & h^q_{{\leq}r} (y)\end{matrix} 
\right\| $. \normalsize \bigskip\\
Then by of 2 of lemma 2.2 \bigskip \footnotesize

 $S(x_i|_{i{=}1,r} )  \in  
(\mathop{\sf ts}^r({\bf R}[x_i|_{i{=}1,r} ])^{{\leq}1} )^{\Delta} $, \normalsize 
\bigskip
\\
by of  1   of lemma   2.2   $S(x_i|_{i{=}1,r} )$ is  polynomially 
  expressed via $({\sigma} _p(x_i|_{i{=}1,r} )|_{p{=}1,r} )$ 
by polynomial  of the degree  ${\leq}{\Delta} $. Since $S(x_i|_{i{=}1,r} 
)$ is  arbitrary   element   in   $\mathop{\sf ts}^r({\bf 
R}[x_i|_{i{=}1,r} ])^{{\leq}{\Delta} } $,   then   \bigskip 
\footnotesize

 $\mathop{\sf ts}^r({\bf R}[x_i|_{i{=}1,r} ])^{{\leq}{\Delta} 
}  \subseteq  (\mathop{\sf ts}^r({\bf R}[x_i|_{i{=}1,r} ])^{{\leq}1} 
)^{\Delta} $, \bigskip \normalsize\\
and      $\mathop{\sf ts}^r({\bf R}[x_i|_{i{=}1,r} ])^{{\leq}{\Delta} } $ 
is polynomially generated by elements 
$({\sigma} _p(x_i|_{i{=}1,r} )|_{p{=}1,r} )$ 
by polynomials of the degree ${\leq}{\Delta} $. \bigskip
\\
Since by of proposition 1.26 the map \footnotesize 
\bigskip

 ${\nu} _s{:}\ \mathop{\sf ts}^r({\bf R}[x_i|_{i{=}1,r} ]) 
\rightarrow  \mathop{\sf TS}^r({\bf R}[x])$, \bigskip \normalsize
\\
is an isomorphism of an semigraded of rings, then inclusions 
\bigskip \footnotesize

 $\mathop{\sf TS}^r({\bf R}[x])^{{\leq}{\Delta} }  \subseteq 
 (\mathop{\sf TS}^r({\bf R}[x])^{{\leq}1} )^{\Delta} $, \bigskip

 $\mathop{\sf ts}^r({\bf R}[x_i|_{i{=}1,r} ])^{{\leq}{\Delta} 
}  \subseteq  (\mathop{\sf ts}^r({\bf R}[x_i|_{i{=}1,r} ])^{{\leq}1} 
)^{\Delta} $. \normalsize \bigskip
\\
are equivalent. Since there holds the second inclusion, 
then there holds the first inclusion.\\
\hfil\vfill\eject
 
 \noindent \normalsize {\bf Theorem 2.3.} {\it Let $r{\geq}1$, 
$d{\geq}r{-}1$, ${\Delta} {=}{\Delta} '{+}{\Delta} ''$, 
${\Delta} '{\geq}0$, ${\Delta} ''{\geq}0$. 
\bigskip\\
1. Let $h_i(x)  \in   {\bf R}[x]^{{\leq}d{+}{\Delta} } 
$  for  $i{=}1,r$. Then \bigskip \footnotesize

 $\bigwedge (h_{{\leq}r} (x)) =  ({-}1)^{r{\cdot}{\Delta} 
'} {\cdot}({-}1)^r$  \bigskip

 \qquad ${\cdot}{\det}\left\| \begin{matrix}  & [y^{{\leq}d{+}{\Delta} 
} ]_{*}.| & y^{{\leq}{\Delta} ''{-}1} {\cdot}y^{d{-}r{+}1{+}{\Delta} 
'} {\cdot}{\sigma} ({\bf x}_i|_{i{=}1,r} )(y) & y^{{\leq}d} {\cdot}y^{{\Delta} 
'}  & y^{{\leq}{\Delta} '{-}1} {\cdot}{\sigma} ({\bf x}_i|_{i{=}1,r} 
)(y) & h_{{\leq}r} (y) \cr
  &              &                                   &     
      &                         &   \cr
\wedge  &              &                                   
& x^{{\leq}d} \h &                         &  \end{matrix} \right\| 
$. \bigskip\\
\\
\normalsize This equality is equivalent to  \bigskip \footnotesize

 ${\det}\left\| h_{{\leq}r} (x_i)|^{i{=}1,r} \right\|  =  ({-}1)^{r{\cdot}{\Delta} 
'} {\cdot}({-}1)^r$  \bigskip

 \qquad  ${\cdot}{\det}\left\| \begin{matrix}[y^{{\leq}d{+}{\Delta} 
} ]_{*}.| & y^{{\leq}{\Delta} ''{-}1} {\cdot}y^{d{-}r{+}1{+}{\Delta} 
'} {\cdot}{\sigma} (x_i|_{i{=}1,r} )(y) & y^{{\leq}d} {\cdot}y^{{\Delta} 
'} \h & y^{{\leq}{\Delta} '{-}1} {\cdot}{\sigma} ({\bf x}_i|_{i{=}1,r} 
)(y) & h_{{\leq}r} (y) \cr
  &                                   &             &      
                   &   \cr
             &                                   & x^{{\leq}d} 
_i|^{i{=}1,r}  &                         &  \end{matrix} \right\| 
$. \bigskip\\
\normalsize\\
2. There holds \bigskip \footnotesize 

 $\bigwedge  ^r({\bf R}[x])^{{\leq}d{+}{\Delta} }  \subseteq 
 \mathop{\sf TS}^r({\bf R}[x])^{{\leq}{\Delta} } {\cdot}\bigwedge 
 ^r({\bf R}[x])^{{\leq}d} $. } \normalsize \bigskip\\
{\bf Proof 1.} There holds \bigskip \footnotesize

 $({-}1)^{{-}r{\cdot}{\Delta} '} {\cdot}{\sigma} _0(x_i|_{i{=}1,r} 
)^{{\Delta} ''} {\cdot}{\sigma} _r(x_i|_{i{=}1,r} )^{{\Delta} 
'} \cdot $  \medskip\\

 \quad ${\det}\left\| \begin{matrix}[y^{{\leq}d{+}{\Delta} 
} ]_{*}.| & y^{{\leq}{\Delta} ''{-}1} {\cdot}y^{d{-}r{+}1{+}{\Delta} 
'} {\cdot}{\sigma} (x_i|_{i{=}1,r} )(y) & y^{{\leq}d} {\cdot}y^{{\Delta} 
'} \h & y^{{\leq}{\Delta} '{-}1} {\cdot}{\sigma} (x_i|_{i{=}1,r} 
)(y) & h_{{\leq}r} (y) \cr
  &                                   &             &      
                   &    \cr
             &                                   & x^{{\leq}d} 
_i|^{i{=}1,r}  &                         &   \end{matrix} \right\| 
$\\
\\

 \quad $\buildrel{ 1}\over  = {\det}\left\| \begin{matrix}[y^{{\leq}d{+}{\Delta} 
} ]_{*}.| & y^{{\leq}{\Delta} ''{-}1} {\cdot}y^{d{-}r{+}1{+}{\Delta} 
'} {\cdot}{\sigma} (x_i|_{i{=}1,r} )(y) & y^{{\leq}d} {\cdot}y^{{\Delta} 
'} \h & y^{{\leq}{\Delta} '{-}1} {\cdot}{\sigma} (x_i|_{i{=}1,r} 
)(y) & h_{{\leq}r} (y) \cr
  &                                   &                 &  
                       &   \cr
             &                                   & x^{{\leq}d} 
_i{\cdot}x^{{\Delta} '} _i|^{i{=}1,r}  &                  
       &  \end{matrix} \right\| $\\
\\

 \quad $\buildrel{ 2}\over  = {\det}\left\| \begin{matrix}[y^{{\leq}d{+}{\Delta} 
} ]_{*}.| & y^{{\leq}{\Delta} ''{-}1} {\cdot}y^{d{-}r{+}1{+}{\Delta} 
'} {\cdot}{\sigma} (x_i|_{i{=}1,r} )(y) & y^{{\leq}d} {\cdot}y^{{\Delta} 
'}  & y^{{\leq}{\Delta} '{-}1} {\cdot}{\sigma} (x_i|_{i{=}1,r} 
)(y) & \ \ \, h_{{\leq}r} (y)\h \cr
  &                                   &           &        
                 &         \cr
             &                                   &         
  &                         & {-}h_{{\leq}r} (x_i)|^{i{=}1,r} 
\end{matrix} \right\| $\\
\\

 \quad $=  {\det}\left\| \begin{matrix}[y^{{\leq}d{+}{\Delta} 
} ]_{*}.| & y^{{\leq}{\Delta} ''{-}1} {\cdot}y^{d{-}r{+}1{+}{\Delta} 
'} {\cdot}{\sigma} (x_i|_{i{=}1,r} )(y) & y^{{\leq}d} {\cdot}y^{{\Delta} 
'}  & y^{{\leq}{\Delta} '{-}1} {\cdot}{\sigma} (x_i|_{i{=}1,r} 
)(y)\end{matrix} \right\| $  \medskip

 \qquad \qquad \qquad \qquad \qquad \qquad \qquad \qquad \qquad 
\qquad \qquad \qquad  \qquad  \qquad  \qquad  \qquad \qquad \qquad 
${\cdot}{\det}\left\| {-}h_{{\leq}r} (x_i)|^{i{=}1,r} \right\| 
$\\

 \quad $\buildrel{ 3}\over  = {\sigma} _0(x_i|_{i{=}1,r} )^{{\Delta} 
''} {\cdot}{\sigma} _r(x_i|_{i{=}1,r} )^{{\Delta} '} {\cdot}{\det}\left\| 
{-}h_{{\leq}r} (x_i)|^{i{=}1,r} \right\| $\\

 \quad $=  {\sigma} _0(x_i|_{i{=}1,r} )^{{\Delta} ''} {\cdot}{\sigma} 
_r(x_i|_{i{=}1,r} )^{{\Delta} '} {\cdot}({-}1)^r{\cdot}{\det}\left\| 
h_{{\leq}r} (x_i)|^{i{=}1,r} \right\| $. \normalsize\\
\\
Dividing both parts of the obtained equality by the common multiple 
 and replacing both parts of the equality we obtain  \bigskip \footnotesize

 ${\det}\left\| h_{{\leq}r} (x_i)|^{i{=}1,r} \right\|  =  ({-}1)^{r{\cdot}{\Delta} 
'} {\cdot}({-}1)^r$\\

 \quad ${\cdot}{\det}\left\| \begin{matrix}[y^{{\leq}d{+}{\Delta} 
} ]_{*}.| & y^{{\leq}{\Delta} ''{-}1} {\cdot}y^{d{-}r{+}1{+}{\Delta} 
'} {\cdot}{\sigma} (x_i|_{i{=}1,r} )(y) & y^{{\leq}d} {\cdot}y^{{\Delta} 
'} \h & y^{{\leq}{\Delta} '{-}1} {\cdot}{\sigma} (x_i|_{i{=}1,r} 
)(y) & h_{{\leq}r} (y) \cr
  &                                   &             &      
                   &   \cr
             &                                   & x^{{\leq}d} 
_i|^{i{=}1,r}  &                         &  \end{matrix} \right\| 
$. \normalsize\\
\hfil\vfill\eject
\noindent 
By of lemma 1.2 and of lemma 1.1 there holds \footnotesize 
\bigskip\\
$\mathop{\sf atr}^r_{{\bf R}[x]} ({\det}\left\| \begin{matrix} 
 & [y^{{\leq}d{+}{\Delta} } ]_{*}.| & y^{{\leq}{\Delta} 
''{-}1} {\cdot}y^{d{-}r{+}1{+}{\Delta} '} {\cdot}{\sigma} 
({\bf x}_i|_{i{=}1,r} )(y) & y^{{\leq}d} {\cdot}y^{{\Delta} 
'}  & y^{{\leq}{\Delta} '{-}1} {\cdot}{\sigma} ({\bf x}_i|_{i{=}1,r} 
)(y) & h_{{\leq}r} (y) \cr
  &              &                                   &     
      &                         &   \cr
\wedge  &              &                                   
& x^{{\leq}d} \h &                         &  \end{matrix} \right\| 
)$  \bigskip

 $=  {\nu} _a\, ({\det}\left\| \begin{matrix}[y^{{\leq}d{+}{\Delta} 
} ]_{*}.| & y^{{\leq}{\Delta} ''{-}1} {\cdot}y^{d{-}r{+}1{+}{\Delta} 
'} {\cdot}{\sigma} (x_i|_{i{=}1,r} )(y) & y^{{\leq}d} {\cdot}y^{{\Delta} 
'} \h & y^{{\leq}{\Delta} '{-}1} {\cdot}{\sigma} ({\bf x}_i|_{i{=}1,r} 
)(y) & h_{{\leq}r} (y) \cr
  &                                   &             &      
                   &   \cr
             &                                   & x^{{\leq}d} 
_i|^{i{=}1,r}  &                         &  \end{matrix} \right\| 
)$, \bigskip

 $\mathop{\sf atr}^r_{{\bf R}[x]} (\bigwedge (h_{{\leq}r} (x))) 
=  {\nu} _a\, ({\det}\left\| h_{{\leq}r} (x_i)|^{i{=}1,r} \right\| 
)$. \normalsize \bigskip
\\
Since by of proposition 1.5 and proposition 1.23 maps 
\bigskip \footnotesize 

 $\mathop{\sf atr}^r_{{\bf R}[x]} {:}\ \mathop{\sf A}^r({\bf 
R}[x]) \rightarrow  \mathop{\sf TA}^r({\bf R}[x])$, ${\nu} _a{:}\ 
\mathop{\sf ta}^r({\bf R}[x|_{i{=}1,r} ]) \rightarrow  \mathop{\sf TA}^r({\bf R}[x])$, 
\bigskip
\\
\normalsize are invertible, then equalities \bigskip \footnotesize

 $\bigwedge (h_{{\leq}r} (x)) =  ({-}1)^{r{\cdot}{\Delta} 
'} {\cdot}({-}1)^r$  \bigskip

 \qquad ${\cdot}{\det}\left\| \begin{matrix}  & [y^{{\leq}d{+}{\Delta} 
} ]_{*}.| & y^{{\leq}{\Delta} ''{-}1} {\cdot}y^{d{-}r{+}1{+}{\Delta} 
'} {\cdot}{\sigma} ({\bf x}_i|_{i{=}1,r} )(y) & y^{{\leq}d} {\cdot}y^{{\Delta} 
'}  & y^{{\leq}{\Delta} '{-}1} {\cdot}{\sigma} ({\bf x}_i|_{i{=}1,r} 
)(y) & h_{{\leq}r} (y) \cr
  &              &                                   &     
      &                         &   \cr
\wedge  &              &                                   
& x^{{\leq}d} \h &                         &  \end{matrix} \right\| 
$, \bigskip\\
\normalsize and \bigskip \footnotesize

 ${\det}\left\| h_{{\leq}r} (x_i)|^{i{=}1,r} \right\|  =  ({-}1)^{r{\cdot}{\Delta} 
'} {\cdot}({-}1)^r$  \bigskip

 \qquad ${\cdot}{\det}\left\| \begin{matrix}[y^{{\leq}d{+}{\Delta} 
} ]_{*}.| & y^{{\leq}{\Delta} ''{-}1} {\cdot}y^{d{-}r{+}1{+}{\Delta} 
'} {\cdot}{\sigma} (x_i|_{i{=}1,r} )(y) & y^{{\leq}d} {\cdot}y^{{\Delta} 
'} \h & y^{{\leq}{\Delta} '{-}1} {\cdot}{\sigma} ({\bf x}_i|_{i{=}1,r} 
)(y) & h_{{\leq}r} (y) \cr
  &                                   &             &      
                   &   \cr
             &                                   & x^{{\leq}d} 
_i|^{i{=}1,r}  &                         &  \end{matrix} \right\| 
$  \bigskip\\
\normalsize are equivalent. Since there holds 
 the second  equality, then there holds the first equality.\\
\hfil\vfill\eject
\noindent 

{\bf Proof 2.} By of 1 of theorem \footnotesize
\\

 $({-}1)^{{-}r{\cdot}{\Delta} '} {\cdot}({-}1)^{{-}r} {\cdot}\bigwedge 
(h_{{\leq}r} (x))$  \bigskip\\

 \qquad $=  {\det}\left\| \begin{matrix}  & [y^{{\leq}d{+}{\Delta} 
} ]_{*}.| & y^{{\leq}{\Delta} ''{-}1} {\cdot}y^{d{-}r{+}1{+}{\Delta} 
'} {\cdot}{\sigma} ({\bf x}_i|_{i{=}1,r} )(y) & y^{{\leq}d} {\cdot}y^{{\Delta} 
'}  & y^{{\leq}{\Delta} '{-}1} {\cdot}{\sigma} ({\bf x}_i|_{i{=}1,r} 
)(y) & h_{{\leq}r} (y) \cr
  &              &                                   &     
      &                         &   \cr
\wedge  &              &                                   
& x^{{\leq}d} \h &                         &  \end{matrix} \right\| 
$  \bigskip

 \qquad $=  {\det}\left\| \begin{matrix}  & [y^{{\leq}d{+}{\Delta} 
} ]_{*}.| & A''({\bf x}_i|_{i{=}1,r} )(y) & y^{{\leq}d} {\cdot}y^{{\Delta} 
'}  & A'({\bf x}_i|_{i{=}1,r} )(y) & h_{{\leq}r} (y) \cr
  &              &                   &           &         
          &   \cr
\wedge  &              &                   & x^{{\leq}d} \h 
&                   &  \end{matrix} \right\| $  \bigskip\\
(by of definition 1.6, with using of notations  
1.43, 1.44,  1.45, 1.46, 1.47, definition 1.5) \bigskip

 \qquad $=  \sum (\mathop{\rm sgn}({\chi} ){\cdot}{\det}\left\| 
\begin{matrix}[y^{{\leq}d{+}{\Delta} } ]_{*}.| & A''({\bf 
x}_i|_{i{=}1,r} )(y) & y^{{\leq}d} {\cdot}y^{{\Delta} '}  
& A'({\bf x}_i|_{i{=}1,r} )(y) & h_{{\leq}r} (y)\end{matrix} 
\right\| |_{{\lnot}{\chi} } $  \bigskip

 \qquad \qquad ${\cdot}\bigwedge (\left( \begin{matrix}0(x)|_{{\times}{\Delta} 
''}  & x^{{\leq}d}  & 0(x)|_{{\times}{\Delta} '}  & 0(x)|_{{\times}r} 
\end{matrix} \right) |_{\chi} )|{\chi} {\in}{\bf C}^r_{d{+}1{+}{\Delta} 
{+}r} )$  \bigskip\\
(Nonzero summands can be only for these 
${\chi} {\in}{\bf C}^r_{d{+}1{+}{\Delta} {+}r} $ 
for which there holds \medskip

 ${\chi} _i{=}0$ for all $i{\in}[1,{\Delta} '']$ and for 
all $i{\in}[{\Delta} ''{+}d{+}1,{\Delta} ''{+}d{+}1{+}{\Delta} 
'{+}r]$, \medskip
\\
denote by ${\mathcal B} $ the set of all such ${\chi} $. Then 
for any ${\chi} {\in}{\mathcal B} $  \medskip

 ${\chi} _i =  {\chi} '_{i{-}{\Delta} ''} $ for all $i{\in}[{\Delta} 
''{+}1,{\Delta} ''{+}d{+}1]$, \medskip\\
where ${\chi} '{\in}{\bf C}^r_{d{+}1} $. Denote by ${\varphi} $ 
the map ${\bf C}^r_{d{+}1}  \rightarrow  {\mathcal B} $, 
${\chi} ' \mapsto  {\chi} $. One can see,   that the map 
${\varphi} $ is invertible.) \bigskip

 \qquad $=  \sum (\mathop{\rm sgn}({\varphi} ({\chi} ')){\cdot}{\det}\left\| 
\begin{matrix}[y^{{\leq}d{+}{\Delta} } ]_{*}.| & A''({\bf 
x}_i|_{i{=}1,r} )(y) & (y^{{\leq}d} {\cdot}y^{{\Delta} '} 
)|_{{\times}{\lnot}{\chi} '}  & A'({\bf x}_i|_{i{=}1,r} )(y) 
& h_{{\leq}r} (y)\end{matrix} \right\| $ \bigskip

 \qquad \qquad ${\cdot}\bigwedge (x^{{\leq}d} |_{{\times}{\chi} 
'} )|{\chi} '{\in}{\bf C}^r_{d{+}1} )$  \bigskip\\
(Since  elements in $[y^{{\leq}d{+}{\Delta} } ]_{*}.A''({\bf 
x}_i|_{i{=}1,r} )(y)$,   $[y^{{\leq}d{+}{\Delta} } ]_{*}.A'({\bf 
x}_i|_{i{=}1,r} )(y)$ belongs to $\mathop{\sf TS}^r({\bf R}[x])^{{\leq}1} 
$, $A''$ include ${\Delta} ''$ columns, $A'$ include ${\Delta} 
'$ columns, on  the rest columns $[y^{{\leq}d{+}{\Delta} 
} ]_{*}.y^{{\leq}d} {\cdot}y^{{\Delta} '} $  
there are elements  in  ${\bf R}$,  
then the determinant 
${\in}(\mathop{\sf TS}^r({\bf R}[x])^{{\leq}1} )^{{\Delta} '{+}{\Delta} ''}  
=  (\mathop{\sf TS}^r({\bf R}[x])^{{\leq}1} )^{\Delta}   
=   \mathop{\sf TS}^r({\bf R}[x])^{{\leq}{\Delta} } $,  
the last  equality holds by of proposition 1.24. 
Since ${\chi} '{\in}{\bf C}^r_{d{+}1} $, then by of 1 of lemma 
2.1 $\bigwedge (x^{{\leq}d} |_{{\times}{\chi} '} ) \in  \bigwedge 
 ^r({\bf R}[x])^{{\leq}d} $.) \bigskip

 \qquad $\in  \mathop{\sf TS}^r({\bf R}[x])^{{\leq}{\Delta} 
} {\cdot}\bigwedge  ^r({\bf R}[x])^{{\leq}d} $. \normalsize \bigskip
\\
Hence, \bigskip \footnotesize

 $({-}1)^{{-}r{\cdot}{\Delta} '} {\cdot}({-}1)^{{-}r} {\cdot}\bigwedge 
(h_{{\leq}r} (x)) \in  \mathop{\sf TS}^r({\bf R}[x])^{{\leq}{\Delta} 
} {\cdot}\bigwedge  ^r({\bf R}[x])^{{\leq}d} $. \normalsize \bigskip
\\
By of proposition 1.17 $\bigwedge  ^r({\bf R}[x])^{{\leq}d{+}{\Delta} 
} $ is lineary  generated by  elements   of the form $\bigwedge (h_{{\leq}r} (x))$, 
where $h_i(x) \in  {\bf R}[x]^{{\leq}d{+}{\Delta} } $ 
for $i{=}1,r$. Then \bigskip \footnotesize

 $\bigwedge  ^r({\bf R}[x])^{{\leq}d{+}{\Delta} }  \subseteq 
 \mathop{\sf TS}^r({\bf R}[x])^{{\leq}{\Delta} } {\cdot}\bigwedge 
 ^r({\bf R}[x])^{{\leq}d} $. \normalsize\\
\hfil\vfill\eject
  
\section{Some application}

\noindent
\normalsize {\bf Definition  3.1.}  Divided difference operator   we call   
and denote by $\nabla (x_1,x_2;x)_{*}$ 
the map  ${\bf R}[x] \rightarrow  {\bf R}[x_1,x_2]$ such 
 that \medskip

 $(x_1{-}x_2){\cdot}(\nabla (x_1,x_2;x)_{*}\, F(x)) =  F(x_1){-}F(x_2)$.
\\
\\
{\bf Property 3.1.} {\it $\nabla (x_1,x_2;x)_{*}$ exists 
and is  unique such map, it is a linear map.} \bigskip
\\
{\bf Property 3.2.} {\it Let ${\bf x}_1 =   x{\otimes}1$, 
 ${\bf x}_2  =   1{\otimes}x$.  Then  $\nabla ({\bf x}_1,{\bf 
x}_2;x)_{*}$  is coproduct ${\bf R}[x] \rightarrow 
 {\bf R}[x]{\otimes}{\bf R}[x]$. This coproduct is 
cocommutative  and coassociative.\medskip\\
Cocommutativity of $\nabla ({\bf x}_1,{\bf x}_2;x)_{*}$ it follows 
from property \medskip

 $\nabla (x_1,x_2;x)_{*} =  \nabla (x_2,x_1;x)_{*}$. \medskip
\\
Coassociativity of $\nabla ({\bf x}_1,{\bf x}_2;x)_{*}$ it follows 
from property \medskip

 $\nabla (x_1,x_2;x)_{*}\, \nabla (x,x_3;x)_{*} =  \nabla (x_2,x_3;x)_{*}\, 
\nabla (x_1,x;x)_{*}$. } \bigskip
\\
{\bf Definition 3.2.} Denote by \bigskip

 $\nabla (x_1;x)_{*} =  {\bf 1}(x_1;x)_{*}$, \medskip

 $\nabla (x_i|_{i{=}1,r} ;x)_{*} =  \nabla (x_i|_{i{=}1,r{-}1} 
;x)_{*}\, \nabla (x,x_r;x)_{*}$ for $r{\geq}3$. \bigskip\\
{\it The equality \bigskip

 $\nabla (x_i|_{i{=}1,r} ;x)_{*} =  \nabla (x_i|_{i{=}1,r{-}1} 
;x)_{*}\, \nabla (x,x_r;x)_{*}$  \bigskip
\\
holds for all $r{\geq}2$. } \bigskip
\\
{\bf Property 3.3.} {\it Let $F(x)$  be a polynomial in $x$, then $\nabla 
(x_i|_{i{=}1,r} ;x)_{*}\, F(x)$  is symmetric polynomial 
in $(x_i|_{i{=}1,r} )$, i. e. belongs to  $\mathop{\sf ts}^r({\bf 
R}[x_i|_{i{=}1,s} ])$.} \bigskip\\
Property 3.3 it follows from property 3.2 of cocommutativity  and 
 coassociativity of the coproduct $\nabla $.\\
\hfil\vfill\eject
 
 \noindent {\bf Lemma 3.1.} {\it Let $F(x)$ be a polynomial in $x$, $r{\geq}1$. There holds \bigskip

 1. $(\nabla (x_i|_{i{=}1,r} ;x)_{*}\, F(x)){\cdot}{\det}\left\| 
x^{{\leq}r{-}1} _i|^{i{=}1,r} \right\|  =  {\det}\left\| \begin{matrix}F(x_i)|^{i{=}1,r} 
 & x^{{\leq}r{-}2} _i|^{i{=}1,r} \end{matrix} \right\| $. \bigskip

 2. ${\det}\left\| x^{{\leq}r{-}1} _i|^{i{=}1,r} \right\|  
=  \prod ((x_i{-}x_k)|_{i{=}1,k{-}1} |_{k{=}2,r} )$. }  \bigskip
\\
{\bf Proof 2.} Equality being proved represent in the form $E(r)$. 
Let us assume,  that for some $r{\geq}2$ there holds $E(r{-}1)$. There holds \bigskip

 ${\det}\left\| x^{{\leq}r{-}1} _i|^{i{=}1,r} \right\|  =  
\left\| \begin{matrix}(x_i{-}x_r){\cdot}x^{{\leq}r{-}2} _i|^{i{=}1,r} 
 & x^0_i|^{i{=}1,r} \end{matrix} \right\| $  \bigskip

 \qquad $=  {\det}\left\| \begin{matrix}(x_i{-}x_r){\cdot}x^{{\leq}r{-}2} 
_i|^{i{=}1,r{-}1}  & x^0_i|^{i{=}1,r{-}1}  \cr
(x_r{-}x_r){\cdot}x^{{\leq}r{-}2} _r\h & x^0_r\h\end{matrix} 
\right\| $  \bigskip

 \qquad $=  {\det}\left\| \begin{matrix}(x_i{-}x_r){\cdot}x^{{\leq}r{-}2} 
_i|^{i{=}1,r{-}1}  & 1|^{{\times}(r{-}1)}  \cr
0|_{{\times}(r{-}1)} \h & 1\h\end{matrix} \right\| $  \bigskip

 \qquad $=  {\det}\left\| (x_i{-}x_r){\cdot}x^{{\leq}r{-}2} 
_i|^{i{=}1,r{-}1} \right\| $  \bigskip

 \qquad $=  \prod ((x_i{-}x_r)|_{i{=}1,r{-}1} ){\cdot}{\det}\left\| 
x^{{\leq}r{-}2} _i|^{i{=}1,r{-}1} \right\| $. \bigskip
\\
The first transition is obtained by following transformations of columns \bigskip

 ${\delta} {=}r{-}1,1{:}\ \left\| x^{\delta} _i|^{i{=}1,r} 
\right\|  \mapsto  \left\| x^{\delta} _i|^{i{=}1,r} \right\| 
{-}x_r{\cdot}\left\| x^{{\delta} {-}1} _i|^{i{=}1,r} \right\| 
$  \bigskip

 \qquad $=  \left\| x_i{\cdot}x^{{\delta} {-}1} _i|^{i{=}1,r} 
\right\| {-}\left\| x_r{\cdot}x^{{\delta} {-}1} _i|^{i{=}1,r} 
\right\|  =  \left\| (x_i{-}x_r){\cdot}x^{{\delta} {-}1} _i|^{i{=}1,r} 
\right\| $. \bigskip
\\
By of $E(r{-}1)$ there holds \bigskip

 $\prod ((x_i{-}x_r)|_{i{=}1,r{-}1} ){\cdot}{\det}\left\| x^{{\leq}r{-}2} 
_i|^{i{=}1,r{-}1} \right\| $  \bigskip

 \qquad $=  \prod ((x_i{-}x_r)|_{i{=}1,r{-}1} ){\cdot}\prod 
((x_i{-}x_k)|_{i{=}1,k{-}1} |_{k{=}2,r{-}1} )$  \bigskip

 \qquad $=  \prod ((x_i{-}x_k)|_{i{=}1,k{-}1} |_{k{=}2,r} )$. 
\bigskip
\\
Thus there holds $E(r)$, i. e. \bigskip

 ${\det}\left\| x^{{\leq}r{-}1} _i|^{i{=}1,r} \right\|  =  
\prod ((x_i{-}x_k)|_{i{=}1,k{-}1} |_{k{=}2,r} )$. \bigskip\\
There holds $E(2)$, since \bigskip

 ${\det}\left\| \begin{matrix}x^1_1 & x^0_1 \cr
x^1_2 & x^0_2\end{matrix} \right\|  =  {\det}\left\| \begin{matrix}x_1 
& 1 \cr
x_2 & 1\end{matrix} \right\|  =  (x_1{-}x_2) =  \prod ((x_i{-}x_k)|_{i{=}1,k{-}1} 
|_{k{=}2,2} )$. \bigskip\\
There holds $E(1)$, since \bigskip

 ${\det}\left\| x^0_1\right\|  =  1 =  \prod ((x_i{-}x_k)|_{i{=}1,k{-}1} 
|_{k{=}2,1} )$. \bigskip\\
Then for all $r{\geq}1$ there holds $E(r)$.\\
\hfil\vfill\eject
 
 \noindent {\bf Proof 1.} Equality being proved represent  
in the form $E(r)$. Let us assume,  that for some $r{\geq}2$ there holds $E(r{-}1)$. 
There holds \bigskip

 ${\det}\left\| \begin{matrix}F(x_i)|^{i{=}1,r}  & x^{{\leq}r{-}2} 
_i|^{i{=}1,r} \end{matrix} \right\| $  \bigskip

 \qquad $=  {\det}\left\| \begin{matrix}(x_i{-}x_r){\cdot}\nabla 
F(x_i,x_r)|^{i{=}1,r}  & (x_i{-}x_r){\cdot}x^{{\leq}r{-}3} _i|^{i{=}1,r} 
 & x^0_i|^{i{=}1,r} \end{matrix} \right\| $  \bigskip

 \qquad $=  {\det}\left\| \begin{matrix}(x_i{-}x_r){\cdot}\nabla 
F(x_i,x_r)|^{i{=}1,r{-}1}  & (x_i{-}x_r){\cdot}x^{{\leq}r{-}3} 
_i|^{i{=}1,r{-}1}  & x^0_i|^{i{=}1,r{-}1}  \cr
  &   &   \cr
(x_r{-}x_r){\cdot}\nabla F(x_r,x_r)\h & (x_r{-}x_r){\cdot}x^{{\leq}r{-}3} 
_r\h & x^0_r\h\end{matrix} \right\| $  \bigskip

 \qquad $=  {\det}\left\| \begin{matrix}(x_i{-}x_r){\cdot}\nabla 
F(x_i,x_r)|^{i{=}1,r{-}1}  & (x_i{-}x_r){\cdot}x^{{\leq}r{-}3} 
_i|^{i{=}1,r{-}1}  & 1|^{{\times}(r{-}1)}  \cr
  &   &   \cr
0\h & 0|_{{\times}(r{-}2)} \h & 1\h\end{matrix} \right\| $ 
 \bigskip

 \qquad $=  {\det}\left\| \begin{matrix}(x_i{-}x_r){\cdot}\nabla 
F(x_i,x_r)|^{i{=}1,r{-}1}  & (x_i{-}x_r){\cdot}x^{{\leq}r{-}3} 
_i|^{i{=}1,r{-}1} \end{matrix} \right\| $  \bigskip

 \qquad $=  \prod ((x_i{-}x_r)|_{i{=}1,r{-}1} ){\cdot}{\det}\left\| 
\begin{matrix}\nabla F(x_i,x_r)|^{i{=}1,r{-}1}  & x^{{\leq}r{-}3} 
_i|^{i{=}1,r{-}1} \end{matrix} \right\| $. \bigskip
\\
The first transition is obtained by foloows transformations of columns \bigskip

 ${\delta} {=}r{-}2,1{:}\ \left\| x^{\delta} _i|^{i{=}1,r} 
\right\|  \mapsto  \left\| x^{\delta} _i|^{i{=}1,r} \right\| 
{-}x_r{\cdot}\left\| x^{{\delta} {-}1} _i|^{i{=}1,r} \right\| 
$  \bigskip

 \qquad $=  \left\| x_i{\cdot}x^{{\delta} {-}1} _i|^{i{=}1,r} 
\right\| {-}\left\| x_r{\cdot}x^{{\delta} {-}1} _i|^{i{=}1,r} 
\right\|  =  \left\| (x_i{-}x_r){\cdot}x^{{\delta} {-}1} _i|^{i{=}1,r} 
\right\| $, \bigskip

 $\left\| F(x_i)|^{i{=}1,r} \right\|  \mapsto  \left\| F(x_i)|^{i{=}1,r} 
\right\| {-}F(x_r){\cdot}\left\| x^0_i|^{i{=}1,r} \right\|  = 
 \left\| F(x_i){-}F(x_r)|^{i{=}1,r} \right\| $  \bigskip

 \qquad $=  \left\| (x_i{-}x_r){\cdot}\nabla F(x_i,x_r)|^{i{=}1,r} 
\right\| $. \bigskip
\\
By of $E(r{-}1)$, the first equality obtained in proof of 
 2  of lemma, and definition 3.2 there holds \bigskip

 $\prod ((x_i{-}x_r)|_{i{=}1,r{-}1} ){\cdot}{\det}\left\| \begin{matrix}\nabla 
F(x_i,x_r)|^{i{=}1,r{-}1}  & x^{{\leq}r{-}3} _i|^{i{=}1,r{-}1} 
\end{matrix} \right\| $  \bigskip

 \qquad $=  \prod ((x_i{-}x_r)|_{i{=}1,r{-}1} ){\cdot}{\det}\left\| 
x^{{\leq}r{-}2} _i|^{i{=}1,r{-}1} \right\| {\cdot}(\nabla (x_i|_{i{=}1,r{-}1} 
;x)_{*}\, \nabla F(x,x_r))$  \bigskip

 \qquad $=  {\det}\left\| x^{{\leq}r{-}1} _i|^{i{=}1,r} \right\| 
{\cdot}(\nabla (x_i|_{i{=}1,r{-}1} ;x)_{*}\, \nabla (x,x_r;x)_{*}\, 
F(x))$  \bigskip

 \qquad $=  {\det}\left\| x^{{\leq}r{-}1} _i|^{i{=}1,r} \right\| 
{\cdot}(\nabla (x_i|_{i{=}1,r} ;x)_{*}\, F(x))$. \bigskip
\\
Thus there holds $E(r)$, i. e. \bigskip

 ${\det}\left\| x^{{\leq}r{-}1} _i|^{i{=}1,r} \right\| {\cdot}(\nabla 
(x_i|_{i{=}1,r} ;x)_{*}\, F(x)) =  {\det}\left\| \begin{matrix}F(x_i)|^{i{=}1,r} 
 & x^{{\leq}r{-}2} _i|^{i{=}1,r} \end{matrix} \right\| $. \bigskip
\\
\\
There holds $E(2)$, since \bigskip

 $\nabla (x_1,x_2;x)_{*}\, F(x) =  \frac{ F(x_1){-}F(x_2)} 
{    x_1{-}x_2}  =  \frac{ {\det}\left\| \begin{matrix}F(x_1) 
& 1 \cr
F(x_2) & 1\end{matrix} \right\| } { {\det}\left\| \begin{matrix}x_1 
& 1 \cr
x_2 & 1\end{matrix} \right\| }  =  \frac{ {\det}\left\| \begin{matrix}F(x_1) 
& x^0_1 \cr
F(x_2) & x^0_2\end{matrix} \right\| } { {\det}\left\| \begin{matrix}x^1_1 
& x^0_1 \cr
x^1_2 & x^0_2\end{matrix} \right\| } $. \bigskip
\\
There holds $E(1)$, since \bigskip

 $\nabla (x_1;x)_{*}\, F(x) =  F(x_1) =  \frac{ {\det}\| F(x_1)\| } 
{ {\det}\| x^0_1\| } $. \bigskip
\\
Then for all $r{\geq}1$ there holds $E(r)$.\\
\hfil\vfill\eject
\noindent 
{\bf Corollary 3.1.} {\it Let $F(x)$  be a polynomial in $x$ of the degree 
${\leq}d$, $d{\geq}r{-}1$, $r{\geq}1$. Then \bigskip

 $\nabla (x_i|_{i{=}1,r} ;x)_{*}\, F(x) =  {\det}\left\| \begin{matrix}[y^{{\leq}d} 
]_{*}.| & y^{{\leq}d{-}r} {\cdot}{\sigma} (x_i|_{i{=}1,r} )(y) 
& F(y) & y^{{\leq}r{-}2} \end{matrix} \right\| $, \bigskip

 $\nabla (x_i|_{i{=}1,r} ;x)_{*}\, F(x) \in  \mathop{\sf ts}^r({\bf 
R}[x_i|_{i{=}1,s} ])^{{\leq}d{-}r{+}1} $, \bigskip

 $F(x){\wedge}\bigwedge (x^{{\leq}r{-}2} ) =  (\nabla ({\bf 
x}_i|_{i{=}1,r} ;x)_{*}\, F(x)){\cdot}\bigwedge (x^{{\leq}r{-}1} 
)$.} \ \bigskip
\\
{\bf Proof.} Let's make the substitution 
${\Delta} \mapsto  d{-}r{+}1$, $r \mapsto  r$ in 2 of theorem 2.2, then \bigskip

 $h_{{\leq}r} (x) \mapsto  (F(x),x^{{\leq}r{-}2} )$, \bigskip

 $S(x_i|_{i{=}1,r} ) \mapsto  \nabla (x_i|_{i{=}1,r} ;x)_{*}\, 
F(x)$. \bigskip
\\
It is satisfies  conditions of 1 of theorem 2.2, since $r \mapsto  r{\geq}1$, 
${\Delta}  \mapsto  d{-}r{+}1{\geq}0$, \bigskip

 ($h_i(x){\in}{\bf R}[x]^{{\leq}r{-}1{+}{\Delta} } $ for 
$i{=}1,r$) $\mapsto $ ($F(x){\in}{\bf R}[x]^{{\leq}d} $, $x^{{\leq}{\delta} 
}  \in  {\bf R}[x]^{{\leq}d} $ for  ${\delta} {=}0,r{-}2$), \bigskip
\\
 that holds, since $F(x)$ is a polynomial of the degree ${\leq}d$ 
and ${\delta} {\leq}r{-}2{\leq}r{-}1{\leq}d$; \bigskip

 $S(x_i|_{i{=}1,r} ) \mapsto  \nabla (x_i|_{i{=}1,r} ;x)_{*}\, 
F(x) \in  \mathop{\sf ts}^r({\bf R}[x_i|_{i{=}1,s} ])$, \bigskip
\\
 that holds according to property 3.3: \bigskip

 $S(x_i|_{i{=}1,r} ){\cdot}{\det}\left\| x^{{\leq}r{-}1} _i|^{i{=}1,r} 
\right\|  =  {\det}\| h_{{\leq}r{-}1} (x_i)|^{i{=}1,r} \| $  
\bigskip

 \qquad $\mapsto  (\nabla (x_i|_{i{=}1,r} ;x)_{*}\, F(x)){\cdot}{\det}\left\| 
x^{{\leq}r{-}1} _i|^{i{=}1,r} \right\|  =  {\det}\left\| \begin{matrix}F(x_i)|^{i{=}1,r} 
 & x^{{\leq}r{-}2} _i|^{i{=}1,r} \end{matrix} \right\| $, \bigskip
\\
 that holds by of {\it 1} of lemma 3.1, and  that is equivalent to
\bigskip

 $(S({\bf x}_i|_{i{=}1,r} ){\cdot}\bigwedge (x^{{\leq}r{-}1} 
) =  \bigwedge (h_{{\leq}r} (x)))$  \bigskip

 \qquad $\mapsto  ((\nabla ({\bf x}_i|_{i{=}1,r} ;x)_{*}\, 
F(x)){\cdot}\bigwedge (x^{{\leq}r{-}1} ) =  F(x){\wedge}\bigwedge 
(x^{{\leq}r{-}2} ))$. \bigskip\\
Then by of 1 of theorem 2.2 \bigskip

 $(S(x_i|_{i{=}1,r} ) =  {\det}\left\| \begin{matrix}[y^{{\leq}d{+}{\Delta} 
} ]_{*}.| & y^{{\leq}{\Delta} {-}1} {\cdot}{\sigma} (x_i|_{i{=}1,r} 
)(y) & h_{{\leq}r} (y)\end{matrix} \right\| )$  \bigskip

 \qquad $\mapsto  (\nabla (x_i|_{i{=}1,r} ;x)_{*}\, F(x) = 
 {\det}\left\| \begin{matrix}[y^{{\leq}d} ]_{*}.| & y^{{\leq}d{-}r} 
{\cdot}{\sigma} (x_i|_{i{=}1,r} )(y) & F(y) & y^{{\leq}r{-}2} 
\end{matrix} \right\| )$, \bigskip

 $(S(x_i|_{i{=}1,r} ) \in  \mathop{\sf ts}^r({\bf R}[x_i|_{i{=}1,s} 
])^{{\leq}{\Delta} } ) \mapsto  (\nabla (x_i|_{i{=}1,r} ;x)_{*}\, 
F(x) \in  \mathop{\sf ts}^r({\bf R}[x_i|_{i{=}1,s} ])^{{\leq}d{-}r{+}1} )$.
\bigskip 
\\ 
{\bf Corollary 3.2.} {\it Let $F(x)$  be a polynomial in $x$
of the degree ${\leq}d$, $d{\geq}r{-}1$, $r{\geq}1$. If \bigskip

 $f(x) =  \prod ((x{-}{\lambda} _i)|_{i{=}1,r} )$, \bigskip
\\
where ${\lambda} _i{\in}{\bf R}$ for $i{=}1,r$, then \bigskip

 $\nabla ({\lambda} _i|_{i{=}1,r} ;x)_{*}\, F(x) =  {\det}\left\| 
\begin{matrix}[y^{{\leq}d} ]_{*}.| & y^{{\leq}d{-}r} {\cdot}f(y) 
& F(y) & y^{{\leq}r{-}2} \end{matrix} \right\| $. } \bigskip
\\
{\bf Proof.} There holds \bigskip

 ${\sigma} (x_i|_{i{=}1,r} )(y) =  \prod ((y{-}x_i)|_{i{=}1,r} 
)$. \bigskip\\
Then \bigskip

 ${\sigma} ({\lambda} _i|_{i{=}1,r} )(y) =  \prod ((y{-}{\lambda} 
_i)|_{i{=}1,r} ) =  f(y)$, \bigskip\\
hence, \bigskip

 $\nabla ({\lambda} _i|_{i{=}1,r} ;x)_{*}\, F(x) =  {\det}\left\| 
\begin{matrix}[y^{{\leq}d} ]_{*}.| & y^{{\leq}d{-}r} {\cdot}{\sigma} 
({\lambda} _i|_{i{=}1,r} )(y) & F(y) & y^{{\leq}r{-}2} \end{matrix} 
\right\| $  \bigskip

 \qquad $=  {\det}\left\| \begin{matrix}[y^{{\leq}d} ]_{*}.| 
& y^{{\leq}d{-}r} {\cdot}f(y) & F(y) & y^{{\leq}r{-}2} \end{matrix} 
\right\| $.
\hfil\vfill\eject
 
\section{Examples}

\noindent
{\bf Example 4.1.} \medskip

 $(a_i|_{i{=}2,5} ) =  (a_2,a_3,a_4,a_5)$, \medskip

 $\| a_i|_{i{=}2,5} \|  =  \| \begin{matrix}a_2 & a_3 & a_4 
& a_5\end{matrix} \| $, \medskip

 $\| a^j|^{j{=}4,6} \|  =  \left\| \begin{matrix}a^4 \cr
a^5 \cr
a^6\end{matrix} \right\| $, \medskip

 $\| a^j_i|^{j{=}4,6} _{i{=}2,5} \|  =  \left\| \begin{matrix}a^4_2 
& a^4_3 & a^4_4 & a^4_5 \cr
a^5_2 & a^5_3 & a^5_4 & a^5_5 \cr
a^6_2 & a^6_3 & a^6_4 & a^6_5\end{matrix} \right\| $, \medskip

 $(a_{{\leq}4} ) =  (a_1,a_2,a_3,a_4)$, \medskip

 $\| a_{{\leq}4} \|  =  \| \begin{matrix}a_1 & a_2 & a_3 & 
a_4\end{matrix} \| $, \medskip

 $\| a^j_{{\leq}4} |^{j{=}4,6} \|  =  \left\| \begin{matrix}a^4_1 
& a^4_2 & a^4_3 & a^4_4 \cr
a^5_1 & a^5_2 & a^5_3 & a^5_4 \cr
a^6_1 & a^6_2 & a^6_3 & a^6_4\end{matrix} \right\| $, \medskip

 $(a|_{{\times}4} ) =  (a,a,a,a)$, \medskip

 $\| a|_{{\times}4} \|  =  \| \begin{matrix}a & a & a & a\end{matrix} 
\| $, \medskip

 $\| a|^{{\times}3} \|  =  \left\| \begin{matrix}a \cr
a \cr
a\end{matrix} \right\| $, \medskip

 $\| a|^{{\times}3} _{{\times}4} \|  =  \left\| \begin{matrix}a 
& a & a & a \cr
a & a & a & a \cr
a & a & a & a\end{matrix} \right\| $.\\
\\
\\
{\bf Example 4.2.} Let \medskip

 $a^j =  \| \begin{matrix}f^j & g^j & h^j\end{matrix} \| $ 
for $i{=}3,5$, \medskip\\
then \medskip

 $\| a^j|^{j{=}3,5} \|  =  \left\| \begin{matrix}f^3 & g^3 
& h^3 \cr
f^4 & g^4 & h^4 \cr
f^5 & g^5 & h^5\end{matrix} \right\| $.\\
\\
\\
{\bf Example 4.3.} $?(3{=}3) =  1$, $?(3{=}2) =  0$, $?(3{<}2) 
=  0$, $?(2{\leq}3) =  1$.\bigskip 
\\
\\
{\bf Example 4.4.} \bigskip

 $\sum (a_1,a_2,a_3,a_4) =  a_1{+}a_2{+}a_3{+}a_4$; \bigskip

 $\prod (a_i|_{i{=}3,6} ) =  \prod (a_3,a_4,a_5,a_6) =  a_3{\cdot}a_4{\cdot}a_5{\cdot}a_6$; 
\bigskip

 $\sum (a|_{{\times}4} ) =  a{+}a{+}a{+}a$; \bigskip

 $\bigwedge (a_{{\leq}4} ) =  \bigwedge (a_1,a_2,a_3,a_4) = 
 a_1{\wedge}a_2{\wedge}a_3{\wedge}a_4$.\\

\hfil\vfill\eject
\noindent 
{\bf Example 4.5.} Let $\top $ be a binary operation, \bigskip

 $\| l\|  =  \left\| \begin{matrix}l^4 \cr 
l^5 \cr
l^6\end{matrix} \right\| $, $\| a\|  =  
\| \begin{matrix}a_2 & a_3 & a_4 & a_5\end{matrix} \| $, \bigskip
\\
then \bigskip

 $\| \begin{matrix}l & \top | & a_2 & a_3 & a_4 & a_5\end{matrix} 
\|  =  \left\| \begin{matrix}l^4\top a_2 & l^4\top a_3 & l^4\top 
a_4 \cr
l^5\top a_2 & l^5\top a_3 & l^5\top a_4 \cr
l^6\top a_2 & l^6\top a_3 & l^6\top a_4\end{matrix} \right\| $.
\bigskip
\\
\\
{\bf Example 4.6.} Let $\top $ be a binary operation,  
\bigskip

 $\| a\|  =  \| \begin{matrix}a_2 & a_3 & a_4 & a_5\end{matrix} 
\| $, \bigskip
\\
then \bigskip

 $\| a\| \top b =  \| \begin{matrix}a_2 & a_3 & a_4 & a_5\end{matrix} 
\| \top b =  \| \begin{matrix}a_2\top b & a_3\top b & a_4\top b & a_5\top b\end{matrix} \| $.\\
\bigskip
\\
\\
{\bf Example 4.7.} Let ${\chi} $ be a map $[1,5] \rightarrow 
 \{0,1\}$ such that  ${\chi} =  (0,0,1,0,1)$. There holds \bigskip

 $\sum {\chi}  =  \sum (0,0,1,0,1) =  0{+}0{+}1{+}0{+}1 =  
2$. \bigskip\\
Hence, ${\chi} {\in}{\bf C}^2_5$. There holds \bigskip

 ${\lnot}{\chi}  =  {\lnot}(0,0,1,0,1) =  ({\lnot}0,{\lnot}0,{\lnot}1,{\lnot}0,{\lnot}1) 
=  (1,1,0,1,0)$. \bigskip
\\
Let $a =  (a_i|_{i{=}1,5} ) =  (a_1,a_2,a_3,a_4,a_5)$, $b 
=  (b_i|_{i{=}1,5} ) =  (b_1,b_2,b_3,b_4,b_5)$. There holds \bigskip

 $a|_{{\times}{\lnot}{\chi} }  =  (a_1,a_2,a_3,a_4,a_5)|_{{\times}(1,1,0,1,0)} 
 =  (a_1|_{{\times}1} ,a_2|_{{\times}1} ,a_3|_{{\times}0} ,a_4|_{{\times}1} 
,a_5|_{{\times}0} ) =  (a_1,a_2,a_4)$, \medskip

 $b|_{{\times}{\chi} }  =  (b_1,b_2,b_3,b_4,b_5)|_{{\times}(0,0,1,0,1)} 
 =  (b_1|_{{\times}0} ,b_2|_{{\times}0} ,b_3|_{{\times}1} ,b_4|_{{\times}0} 
,b_5|_{{\times}1} ) =  (b_3,b_5)$. \bigskip\\
Then \bigskip

 $(a|_{{\times}{\lnot}{\chi} } ,b|_{{\times}{\chi} } ) =  (a_1,a_2,a_4,b_3,b_5) 
=  (a_{{\tau} (1)} ,a_{{\tau} (2)} ,a_{{\tau} (3)} ,b_{{\tau} 
(4)} ,b_{{\tau} (5)} )$, \bigskip\\
where ${\tau} {=}(1,2,4,3,5)$. There holds ${\tau} {\in}!5$, 
the map  ${\tau} $  is  monotonous  on  $[1,3]$, and 
monotonous on  $[4,5]$, hence, ${\tau} {\in}{\mathcal 
C} ^2_5$. There holds \bigskip

 ${\chi} ({\tau} ) =  {\chi} (1,2,4,3,5) =  ({\chi} (1),{\chi} 
(2),{\chi} (4),{\chi} (3),{\chi} (5)) =  (0,0,0,1,1)$.\\
\hfil\vfill\eject
\noindent 
{\bf Example 4.8.}\\

 ${\det}\left\| \begin{matrix}  & a^1_1 & a^1_2 & a^1_3 & a^1_4 
\cr
  & a^2_1 & a^2_2 & a^2_3 & a^2_4 \cr
\wedge  & h_1 & h_2 & h_3 & h_4\end{matrix} \right\| $\\
\\

 \qquad $=  \mathop{\rm sgn}(1,2,3,4){\cdot}{\det}\left\| \begin{matrix}a^1_1 
& a^1_2 \cr
a^2_1 & a^2_2\end{matrix} \right\| {\cdot}(h_3{\wedge}h_4){+}\mathop{\rm 
sgn}(1,3,2,4){\cdot}{\det}\left\| \begin{matrix}a^1_1 & a^1_3 
\cr
a^2_1 & a^2_3\end{matrix} \right\| {\cdot}(h_2{\wedge}h_4)$
\\

 \qquad \qquad ${+}\mathop{\rm sgn}(2,3,1,4){\cdot}{\det}\left\| 
\begin{matrix}a^1_2 & a^1_3 \cr
a^2_2 & a^2_3\end{matrix} \right\| {\cdot}(h_1{\wedge}h_4){+}\mathop{\rm 
sgn}(1,4,2,3){\cdot}{\det}\left\| \begin{matrix}a^1_1 & a^1_4 
\cr
a^2_1 & a^2_4\end{matrix} \right\| {\cdot}(h_2{\wedge}h_3)$
\\

 \qquad \qquad ${+}\mathop{\rm sgn}(2,4,1,3){\cdot}{\det}\left\| 
\begin{matrix}a^1_2 & a^1_4 \cr
a^2_2 & a^2_4\end{matrix} \right\| {\cdot}(h_1{\wedge}h_3){+}{+}\mathop{\rm 
sgn}(3,4,1,2){\cdot}{\det}\left\| \begin{matrix}a^1_3 & a^1_4 
\cr
a^2_3 & a^2_4\end{matrix} \right\| {\cdot}(h_1{\wedge}h_2)$.
\\
\\

 \qquad $=  \mathop{\rm sgn}(0,0,1,1){\cdot}{\det}\left\| \begin{matrix}a^1_1 
& a^1_2 \cr
a^2_1 & a^2_2\end{matrix} \right\| {\cdot}(h_3{\wedge}h_4){+}\mathop{\rm 
sgn}(0,1,0,1){\cdot}{\det}\left\| \begin{matrix}a^1_1 & a^1_3 
\cr
a^2_1 & a^2_3\end{matrix} \right\| {\cdot}(h_2{\wedge}h_4)$
\\

 \qquad \qquad ${+}\mathop{\rm sgn}(1,0,0,1){\cdot}{\det}\left\| 
\begin{matrix}a^1_2 & a^1_3 \cr
a^2_2 & a^2_3\end{matrix} \right\| {\cdot}(h_1{\wedge}h_4){+}\mathop{\rm 
sgn}(0,1,1,0){\cdot}{\det}\left\| \begin{matrix}a^1_1 & a^1_4 
\cr
a^2_1 & a^2_4\end{matrix} \right\| {\cdot}(h_2{\wedge}h_3)$
\\

 \qquad \qquad ${+}\mathop{\rm sgn}(1,0,0,1){\cdot}{\det}\left\| 
\begin{matrix}a^1_2 & a^1_4 \cr
a^2_2 & a^2_4\end{matrix} \right\| {\cdot}(h_1{\wedge}h_3){+}\mathop{\rm 
sgn}(1,1,0,0){\cdot}{\det}\left\| \begin{matrix}a^1_3 & a^1_4 
\cr
a^2_3 & a^2_4\end{matrix} \right\| {\cdot}(h_1{\wedge}h_2)$.
\\
\\

 \qquad $=  {\det}\left\| \begin{matrix}a^1_1 & a^1_2 \cr
a^2_1 & a^2_2\end{matrix} \right\| {\cdot}(h_3{\wedge}h_4){-}{\det}\left\| 
\begin{matrix}a^1_1 & a^1_3 \cr
a^2_1 & a^2_3\end{matrix} \right\| {\cdot}(h_2{\wedge}h_4){+}{\det}\left\| 
\begin{matrix}a^1_2 & a^1_3 \cr
a^2_2 & a^2_3\end{matrix} \right\| {\cdot}(h_1{\wedge}h_4)$
\\

 \qquad \qquad ${+}{\det}\left\| \begin{matrix}a^1_1 & a^1_4 
\cr
a^2_1 & a^2_4\end{matrix} \right\| {\cdot}(h_2{\wedge}h_3){-}{\det}\left\| 
\begin{matrix}a^1_2 & a^1_4 \cr
a^2_2 & a^2_4\end{matrix} \right\| {\cdot}(h_1{\wedge}h_3){+}{\det}\left\| 
\begin{matrix}a^1_3 & a^1_4 \cr
a^2_3 & a^2_4\end{matrix} \right\| {\cdot}(h_1{\wedge}h_2)$.
\\
\\
{\bf Example 4.9.}\\

 ${\det}\| \begin{matrix}\wedge  & h_1 & h_2 & h_3\end{matrix} 
\|  =  (h_1{\wedge}h_2{\wedge}h_3)$,\\
\\

 ${\det}\left\| \begin{matrix}    & a^1_1 & a^1_2 & a^1_3 \cr
\wedge  & h_1 & h_2 & h_3\end{matrix} \right\|  =  {\det}\| 
a^1_1\| {\cdot}(h_2{\wedge}h_3){-}{\det}\| a^1_2\| {\cdot}(h_1{\wedge}h_3){+}{\det}\| 
a^1_3\| {\cdot}(h_1{\wedge}h_2)$,\\
\\

 ${\det}\left\| \begin{matrix}  & a^1_1 & a^1_2 & a^1_3 \cr
  & a^2_1 & a^2_2 & a^2_3 \cr
\wedge  & h_1 & h_2 & h_3\end{matrix} \right\|  =  {\det}\left\| 
\begin{matrix}a^1_1 & a^1_2 \cr
a^2_1 & a^2_2\end{matrix} \right\| {\cdot}(h_3){-}{\det}\left\| 
\begin{matrix}a^1_1 & a^1_3 \cr
a^2_1 & a^2_3\end{matrix} \right\| {\cdot}(h_2){+}{\det}\left\| 
\begin{matrix}a^1_2 & a^1_3 \cr
a^2_2 & a^2_3\end{matrix} \right\| {\cdot}(h_1)$,\\
\\

 ${\det}\left\| \begin{matrix}  & a^1_1 & a^1_2 & a^1_3 \cr
  & a^2_1 & a^2_2 & a^2_3 \cr
  & a^3_1 & a^3_2 & a^3_3 \cr
\wedge  & h_1 & h_2 & h_3\end{matrix} \right\|  =  {\det}\left\| 
\begin{matrix}a^1_1 & a^1_2 & a^1_3 \cr
a^2_1 & a^2_2 & a^2_3 \cr
a^3_1 & a^3_2 & a^3_3\end{matrix} \right\| $,\\
\\

 ${\det}\left\| \begin{matrix}  & a^1_1 & a^1_2 & a^1_3 \cr
  & a^2_1 & a^2_2 & a^2_3 \cr
  & a^3_1 & a^3_2 & a^3_3 \cr
  & a^4_1 & a^4_2 & a^4_3 \cr
\wedge  & h_1 & h_2 & h_3\end{matrix} \right\|  =  0$.
\hfil\vfill\eject

 \noindent {\bf Convention 4.1.} Let $r{\geq}1$. In following examples 
put \medskip

 $\overline {{\sigma} }_p =  {\sigma} _p(x_i|_{i{=}1,r} )$ 
for $p{=}0,r$,

 $\widetilde {{\sigma} }_p =  {\sigma} _p({\bf x}_i|_{i{=}1,r} 
)$ for $p{=}0,r$.\\
\\
{\bf Example 4.10.} \bigskip \footnotesize

 $x^{{\leq}0}  =  \| x^0\| $, $x^{{\leq}1}  =  \| \begin{matrix}x^1 
& x^0\end{matrix} \| $, $x^{{\leq}2}  =  \| \begin{matrix}x^2 
& x^1 & x^0\end{matrix} \| $, $x^{{\leq}3}  =  \| \begin{matrix}x^3 
& x^2 & x^1 & x^0\end{matrix} \| $. \normalsize \bigskip\\
{\bf Example 4.10.} \bigskip \footnotesize

 $\| [y^{{\leq}0} ]_{*}\|  =  \| [y^0]_{*}\| $, $\| [y^{{\leq}1} 
]_{*}\|  =  \left\| \begin{matrix}[y^1]_{*} \cr
[y^0]_{*}\end{matrix} \right\| $, $\| [y^{{\leq}2} ]_{*}\| 
 =  \left\| \begin{matrix}[y^2]_{*} \cr
[y^1]_{*} \cr
[y^0]_{*}\end{matrix} \right\| $, $\| [y^{{\leq}3} ]_{*}\| 
 =  \left\| \begin{matrix}[y^3]_{*} \cr
[y^2]_{*} \cr
[y^1]_{*} \cr
[y^0]_{*}\end{matrix} \right\| $. \normalsize \bigskip\\
{\bf Example 4.12.} Let ${\Delta} {=}2$, $r{=}3$, $d{=}r{-}1$. 
Let \bigskip \footnotesize

 $H(x) =  \sum\limits^{d{+}{\Delta} } _{{\delta} {=}0} H_{\delta} 
{\cdot}x^{\delta}  =  H_4{\cdot}x^4{+}H_3{\cdot}x^3{+}H_2{\cdot}x^2{+}H_1{\cdot}x^1{+}H_0{\cdot}x^0$. 
\normalsize \bigskip\\
Then \bigskip \footnotesize

 $[y^4]_{*}.H(y){=}H_4$, $[y^3]_{*}.H(y){=}H_3$, $[y^2]_{*}.H(y){=}H_2$, 
$[y^1]_{*}.H(y){=}H_1$, $[y^0]_{*}.H(y){=}H_0$; \bigskip

 $\| [y^{{\leq}d{+}{\Delta} } ]_{*}.|H(y)\|  =  \| [y^{{\leq}4} 
]_{*}.|H(y)\|  =  \left\| \begin{matrix}[y^4]_{*}.H(y) \cr
[y^3]_{*}.H(y) \cr
[y^2]_{*}.H(y) \cr
[y^1]_{*}.H(y) \cr
[y^0]_{*}.H(y)\end{matrix} \right\|  =  \left\| \begin{matrix}H_4 
\cr
H_3 \cr
H_2 \cr
H_1 \cr
H_0\end{matrix} \right\| $. \normalsize \bigskip\\
{\bf Example 4.13.} Let ${\Delta} {=}2$, $r{=}3$, $d{=}r{-}1$. 
Let \bigskip \footnotesize

 $H(x) =  \sum\limits^{ d} _{{\delta} {=}0} H_{\delta} {\cdot}x^{\delta} 
 =  H_3{\cdot}x^3{+}H_2{\cdot}x^2{+}H_1{\cdot}x^1{+}H_0{\cdot}x^0$. 
\normalsize \bigskip\\
Then \bigskip \footnotesize

 $[y^4]_{*}.y^1{\cdot}H(y){=}H_3$, $[y^3]_{*}.y^1{\cdot}H(y){=}H_2$, 
$[y^2]_{*}.y^1{\cdot}H(y){=}H_1$, $[y^1]_{*}.y^1{\cdot}H(y){=}H_0$, 
$[y^0]_{*}.y^1{\cdot}H(y){=}0$; \bigskip

 $[y^4]_{*}.y^0{\cdot}H(y){=}0$, $[y^3]_{*}.y^0{\cdot}H(y){=}H_3$, 
$[y^2]_{*}.y^0{\cdot}H(y){=}H_2$, $[y^1]_{*}.y^0{\cdot}H(y){=}H_1$, 
$[y^0]_{*}.y^0{\cdot}H(y){=}H_0$; \bigskip\\
\\

 $\| [y^{{\leq}d{+}{\Delta} } ]_{*}.|y^{{\leq}{\Delta} 
{-}1} {\cdot}H(y)\|  =  \| [y^{{\leq}4} ]_{*}.|y^{{\leq}1} {\cdot}H(y)\| 
 =  \left\| \begin{matrix}[y^4]_{*}.y^1{\cdot}H(y) & [y^4]_{*}.y^0{\cdot}H(y) 
\cr
[y^3]_{*}.y^1{\cdot}H(y) & [y^3]_{*}.y^0{\cdot}H(y) \cr
[y^2]_{*}.y^1{\cdot}H(y) & [y^2]_{*}.y^0{\cdot}H(y) \cr
[y^1]_{*}.y^1{\cdot}H(y) & [y^1]_{*}.y^0{\cdot}H(y) \cr
[y^0]_{*}.y^1{\cdot}H(y) & [y^0]_{*}.y^0{\cdot}H(y)\end{matrix} 
\right\|  =  \left\| \begin{matrix}H_3 &   \cr
H_2 & H_3 \cr
H_1 & H_2 \cr
H_0 & H_1 \cr
     & H_0\end{matrix} \right\| $. \normalsize\\
\hfil\vfill\eject
 
 \noindent {\bf Example 4.14.} Let ${\Delta} {=}2$, $r{=}3$, 
$d{=}r{-}1$. Let \bigskip \footnotesize

 $h_i(x) =  \sum\limits^{d{+}{\Delta} } _{{\delta} {=}0} 
h_{i,{\delta} } {\cdot}x^{\delta}  =  h_{i,4} {\cdot}x^4{+}h_{i,3} 
{\cdot}x^3{+}h_{i,2} {\cdot}x^2{+}h_{i,1} {\cdot}x^1{+}h_{i,0} 
{\cdot}x^0$ for $i{=}1,r$. \normalsize \bigskip\\
Then \bigskip \footnotesize

 $\| [y^{{\leq}d{+}{\Delta} } ]_{*}.|h_{{\leq}r} (y)\|  
=  \| \begin{matrix}[y^{{\leq}d{+}{\Delta} } ]_{*}.h_1(y) 
& [y^{{\leq}d{+}{\Delta} } ]_{*}.h_2(y) & [y^{{\leq}d{+}{\Delta} 
} ]_{*}.h_3(y)\end{matrix} \|  =  \left\| \begin{matrix}h_{1,4} 
 & h_{2,4}  & h_{3,4}  \cr
h_{1,3}  & h_{2,3}  & h_{3,3}  \cr
h_{1,2}  & h_{2,2}  & h_{3,2}  \cr
h_{1,1}  & h_{2,1}  & h_{3,1}  \cr
h_{1,0}  & h_{2,0}  & h_{3,0} \end{matrix} \right\| $. \normalsize 
\bigskip\\
{\bf Example 4.15.} Let $r{=}3$, \bigskip \footnotesize

 $\mathfrak{ s} (x) =  \sum\limits^{ r} _{{\delta} {=}0} \mathfrak{ 
s} _{r{-}{\delta} } {\cdot}x^{\delta}  =  \mathfrak{ s} _0{\cdot}x^3{+}\mathfrak{ 
s} _1{\cdot}x^2{+}\mathfrak{ s} _2{\cdot}x^1{+}\mathfrak{ s} 
_3{\cdot}x^0$. \normalsize \bigskip\\
If ${\Delta} {=}2$, $d{=}r{-}1$, then \bigskip \footnotesize

 $\| [y^{{\leq}d{+}{\Delta} } ]_{*}.|y^{{\leq}{\Delta} 
{-}1} {\cdot}\mathfrak{ s} (y)\|  =  \| [y^{{\leq}4} ]_{*}.|y^{{\leq}1} 
{\cdot}\mathfrak{ s} (y)\|  =  \| \begin{matrix}[y^{{\leq}4} 
]_{*}.| & y^1{\cdot}\mathfrak{ s} (y) & y^0{\cdot}\mathfrak{ 
s} (y)\end{matrix} \| $  \bigskip

 \qquad $=  \left\| \begin{matrix}[y^4]_{*}.y^1{\cdot}\mathfrak{ 
s} (y) & [y^4]_{*}.y^0{\cdot}\mathfrak{ s} (y) \cr
[y^3]_{*}.y^1{\cdot}\mathfrak{ s} (y) & [y^3]_{*}.y^0{\cdot}\mathfrak{ 
s} (y) \cr
[y^2]_{*}.y^1{\cdot}\mathfrak{ s} (y) & [y^2]_{*}.y^0{\cdot}\mathfrak{ 
s} (y) \cr
[y^1]_{*}.y^1{\cdot}\mathfrak{ s} (y) & [y^1]_{*}.y^0{\cdot}\mathfrak{ 
s} (y) \cr
[y^0]_{*}.y^1{\cdot}\mathfrak{ s} (y) & [y^0]_{*}.y^0{\cdot}\mathfrak{ 
s} (y)\end{matrix} \right\|  =  \left\| \begin{matrix}\mathfrak{ 
s} _0 &   \cr
\mathfrak{ s} _1 & \mathfrak{ s} _0 \cr
\mathfrak{ s} _2 & \mathfrak{ s} _1 \cr
\mathfrak{ s} _3 & \mathfrak{ s} _2 \cr
     & \mathfrak{ s} _3\end{matrix} \right\| $. \normalsize 
\bigskip\\
If ${\Delta} {=}3$, $d{=}3$, ${\Delta} '{=}1$, ${\Delta} 
''{=}2$, then \bigskip \footnotesize

 $\| [y^{{\leq}d{+}{\Delta} } ]_{*}.|y^{{\leq}{\Delta} 
''{-}1} {\cdot}y^{d{-}r{+}1{+}{\Delta} '} {\cdot}\mathfrak{ 
s} (y)\|  =  \| [y^{{\leq}6} ]_{*}.|y^{{\leq}1} {\cdot}y^2{\cdot}\mathfrak{ 
s} (y)\|  =  \| \begin{matrix}[y^{{\leq}6} ]_{*}.| & y^1{\cdot}y^2{\cdot}\mathfrak{ 
s} (y) & y^0{\cdot}y^2{\cdot}\mathfrak{ s} (y)\end{matrix} \| 
$  \bigskip

 \qquad $=  \| \begin{matrix}[y^{{\leq}6} ]_{*}.| & y^3{\cdot}\mathfrak{ 
s} (y) & y^2{\cdot}\mathfrak{ s} (y)\end{matrix} \|  =  \left\| 
\begin{matrix}[y^6]_{*}.y^3{\cdot}\mathfrak{ s} (y) & [y^6]_{*}.y^2{\cdot}\mathfrak{ 
s} (y) \cr
[y^5]_{*}.y^3{\cdot}\mathfrak{ s} (y) & [y^5]_{*}.y^2{\cdot}\mathfrak{ 
s} (y) \cr
[y^4]_{*}.y^3{\cdot}\mathfrak{ s} (y) & [y^4]_{*}.y^2{\cdot}\mathfrak{ 
s} (y) \cr
[y^3]_{*}.y^3{\cdot}\mathfrak{ s} (y) & [y^3]_{*}.y^2{\cdot}\mathfrak{ 
s} (y) \cr
[y^2]_{*}.y^3{\cdot}\mathfrak{ s} (y) & [y^2]_{*}.y^2{\cdot}\mathfrak{ 
s} (y) \cr
[y^1]_{*}.y^3{\cdot}\mathfrak{ s} (y) & [y^1]_{*}.y^2{\cdot}\mathfrak{ 
s} (y) \cr
[y^0]_{*}.y^3{\cdot}\mathfrak{ s} (y) & [y^0]_{*}.y^2{\cdot}\mathfrak{ 
s} (y)\end{matrix} \right\|  =  \left\| \begin{matrix}\mathfrak{ 
s} _0 &   \cr
\mathfrak{ s} _1 & \mathfrak{ s} _0 \cr
\mathfrak{ s} _2 & \mathfrak{ s} _1 \cr
\mathfrak{ s} _3 & \mathfrak{ s} _2 \cr
     & \mathfrak{ s} _3 \cr
0 & 0 \cr
0 & 0\end{matrix} \right\| $, \bigskip\\

 $\| [y^{{\leq}d{+}{\Delta} } ]_{*}.|y^{{\leq}{\Delta} 
'{-}1} {\cdot}\mathfrak{ s} (y)\|  =  \| [y^{{\leq}6} ]_{*}.|y^{{\leq}0} 
{\cdot}\mathfrak{ s} (y)\|  =  \| [y^{{\leq}6} ]_{*}.|\mathfrak{ 
s} (y)\|  =  \left\| \begin{matrix}0\h \cr
0\h \cr
0\h \cr
\mathfrak{ s} _0 \cr
\mathfrak{ s} _1 \cr
\mathfrak{ s} _2 \cr
\mathfrak{ s} _3\end{matrix} \right\| $. \normalsize \bigskip
\\
{\bf Example 4.16.} Let ${\Delta} {=}2$, $r{=}3$, $d{=}r{-}1$, 
${\Delta} '{=}1$. Then \bigskip \footnotesize

 $\| [y^{{\leq}d{+}{\Delta} } ]_{*}.|y^{{\leq}d} {\cdot}y^{{\Delta} 
'} \|  =  \| [y^{{\leq}4} ]_{*}.|y^{{\leq}2} {\cdot}y^1\|  = 
 \| \begin{matrix}[y^{{\leq}4} ]_{*}.| & y^2{\cdot}y^1 & y^1{\cdot}y^1 
& y^0{\cdot}y^1\end{matrix} \|  =  \| \begin{matrix}[y^{{\leq}4} 
]_{*}.| & y^3 & y^2 & y^1\end{matrix} \| $  \bigskip

 \qquad $=  \left\| \begin{matrix}[y^4]_{*}.y^3 & [y^4]_{*}.y^2 
& [y^4]_{*}.y^1 \cr
[y^3]_{*}.y^3 & [y^3]_{*}.y^2 & [y^3]_{*}.y^1 \cr
[y^2]_{*}.y^3 & [y^2]_{*}.y^2 & [y^2]_{*}.y^1 \cr
[y^1]_{*}.y^3 & [y^1]_{*}.y^2 & [y^1]_{*}.y^1 \cr
[y^0]_{*}.y^3 & [y^0]_{*}.y^2 & [y^0]_{*}.y^1\end{matrix} \right\| 
 =  \left\| \begin{matrix}    &     &   \cr
1 &     &   \cr
    & 1 &   \cr
    &     & 1 \cr
    &     &  \end{matrix} \right\| $. \normalsize \bigskip
\\
\hfil\vfill\eject
 
 \noindent {\bf Example 4.17.} Let $r{=}3$. \medskip\\
If $d{=}r{-}1$, then \bigskip \footnotesize

 $\| x^{{\leq}d} _i|^{i{=}1,r} \|  =  \| x^{{\leq}2} _i|^{i{=}1,3} 
\|  =  \left\| \begin{matrix}x^{{\leq}2} _1 \cr
x^{{\leq}2} _2 \cr
x^{{\leq}2} _3\end{matrix} \right\|  =  \left\| \begin{matrix}x^2_1 
& x^1_1 & x^0_1 \cr
x^2_2 & x^1_2 & x^0_2 \cr
x^2_3 & x^1_3 & x^0_3\end{matrix} \right\| $. \bigskip

 $\bigwedge (x^{{\leq}r{-}1} ) =  \bigwedge (x^{{\leq}2} ) 
=  \bigwedge (x^2,x^1,x^0) =  x^2{\wedge}x^1{\wedge}x^0$. \normalsize 
\bigskip\\
If $d{=}3$, then \bigskip \footnotesize

 $\| x^{{\leq}d} _i|^{i{=}1,r} \|  =  \| x^{{\leq}3} _i|^{i{=}1,3} 
\|  =  \left\| \begin{matrix}x^{{\leq}3} _1 \cr
x^{{\leq}3} _2 \cr
x^{{\leq}3} _3\end{matrix} \right\|  =  \left\| \begin{matrix}x^3_1 
& x^2_1 & x^1_1 & x^0_1 \cr
x^3_2 & x^2_2 & x^1_2 & x^0_2 \cr
x^3_3 & x^2_3 & x^1_3 & x^0_3\end{matrix} \right\| $. \normalsize 
\bigskip\\
\\
{\bf Example 4.18.} Let $r{=}3$, $d{=}r{-}1$, ${\Delta} 
'{=}1$. \medskip\\
If $d{=}r{-}1$, then \bigskip \footnotesize

 $\| x^{{\leq}d} _i{\cdot}x^{{\Delta} '} _i|^{i{=}1,r} \| 
 =  \| x^{{\leq}2} _i{\cdot}x^1_i|^{i{=}1,3} \|  =  \left\| \begin{matrix}x^{{\leq}2} 
_1{\cdot}x^1_1 \cr
x^{{\leq}2} _2{\cdot}x^1_2 \cr
x^{{\leq}2} _3{\cdot}x^1_3\end{matrix} \right\|  =  \left\| 
\begin{matrix}x^2_i{\cdot}x^1_i & x^1_i{\cdot}x^1_i & x^0_i{\cdot}x^1_i 
\cr
x^2_i{\cdot}x^1_i & x^1_i{\cdot}x^1_i & x^0_i{\cdot}x^1_i \cr
x^2_i{\cdot}x^1_i & x^1_i{\cdot}x^1_i & x^0_i{\cdot}x^1_i\end{matrix} 
\right\|  =  \left\| \begin{matrix}x^3_1 & x^2_1 & x^1_1 \cr
x^3_2 & x^2_2 & x^1_2 \cr
x^3_3 & x^2_3 & x^1_3\end{matrix} \right\| $. \normalsize \bigskip
\\
\\
If $d{=}3$, then \bigskip \footnotesize

 $\| x^{{\leq}d} _i{\cdot}x^{{\Delta} '} _i|^{i{=}1,r} \| 
 =  \| x^{{\leq}3} _i{\cdot}x^1_i|^{i{=}1,3} \|  =  \left\| \begin{matrix}x^{{\leq}3} 
_1{\cdot}x^1_1 \cr
x^{{\leq}3} _2{\cdot}x^1_2 \cr
x^{{\leq}3} _3{\cdot}x^1_3\end{matrix} \right\|  =  \left\| 
\begin{matrix}x^3_i{\cdot}x^1_i & x^2_i{\cdot}x^1_i & x^1_i{\cdot}x^1_i 
& x^0_i{\cdot}x^1_i \cr
x^3_i{\cdot}x^1_i & x^2_i{\cdot}x^1_i & x^1_i{\cdot}x^1_i & 
x^0_i{\cdot}x^1_i \cr
x^3_i{\cdot}x^1_i & x^2_i{\cdot}x^1_i & x^1_i{\cdot}x^1_i & 
x^0_i{\cdot}x^1_i\end{matrix} \right\| $  \bigskip

 \qquad $=  \left\| \begin{matrix}x^4_1 & x^3_1 & x^2_1 & x^1_1 
\cr
x^4_2 & x^3_2 & x^2_2 & x^1_2 \cr
x^4_3 & x^3_3 & x^2_3 & x^1_3\end{matrix} \right\| $. \normalsize 
\bigskip\\
\\
{\bf Example 4.19.} Let $r{=}3$. Then \bigskip \footnotesize

 $\| h_{{\leq}r} (x_i)|^{i{=}1,r} \|  =  \| h_{{\leq}3} (x_i)|^{i{=}1,r} 
\|  =  \left\| \begin{matrix}h_{{\leq}3} (x_1) \cr
h_{{\leq}3} (x_2) \cr
h_{{\leq}3} (x_3)\end{matrix} \right\|  =  \left\| \begin{matrix}h_1(x_1) 
& h_2(x_1) & h_3(x_1) \cr
h_1(x_2) & h_2(x_2) & h_3(x_2) \cr
h_1(x_3) & h_2(x_3) & h_3(x_3)\end{matrix} \right\| $. \bigskip

 $\bigwedge (h_{{\leq}r} (x)) =  \bigwedge (h_{{\leq}3} (x)) 
=  \bigwedge (h_1(x),h_2(x),h_3(x)) =  h_1(x){\wedge}h_2(x){\wedge}h_3(x)$. 
\normalsize \bigskip\\
\hfil\vfill\eject
 
 \noindent {\bf Example 4.20.} (to lemma 2.2) {\it Let ${\Delta} 
{=}2$, $r{=}3$, $d{=}r{-}1$. \medskip\\
Let $h_i(x) \in  {\bf R}[x]^{{\leq}d{+}{\Delta} } $ for 
$i{=}1,r$, i. e. \bigskip

 $h_i(x) =  \sum\limits^{d{+}{\Delta} } _{{\delta} {=}0} 
h_{i,{\delta} } {\cdot}x^{\delta} $ for $i{=}1,r$. \bigskip\\
Let \bigskip

 $S(x_1,x_2,x_3) =  {\det}\left\| \begin{matrix}\overline {\sigma} 
_0 &   & h_{1,4}  & h_{2,4}  & h_{3,4}  \cr
\overline {\sigma} _1 & \overline {\sigma} _0 & h_{1,3}  & 
h_{2,3}  & h_{3,3}  \cr
\overline {\sigma} _2 & \overline {\sigma} _1 & h_{1,2}  & 
h_{2,2}  & h_{3,2}  \cr
\overline {\sigma} _3 & \overline {\sigma} _2 & h_{1,1}  & 
h_{2,1}  & h_{3,1}  \cr
  & \overline {\sigma} _3 & h_{1,0}  & h_{2,0}  & h_{3,0} \end{matrix} 
\right\| $, \bigskip\\
Then \bigskip
\\
1) $S(x_1,x_2,x_3)$ is polynomially expressed via  $(\widetilde 
{{\sigma} }_1,\widetilde {{\sigma} }_2,\widetilde {{\sigma} }_3)$ 
by polynomial  of the degree  ${\leq}{\Delta} $; \medskip
\\
2) $S(x_1,x_2,x_3) \in  (\mathop{\sf ts}^r({\bf R}[x_1,x_2,x_3])^{{\leq}1} 
)^{\Delta} $; \medskip
\\
3) $S(x_1,x_2,x_3) \in  \mathop{\sf ts}^r({\bf R}[x_1,x_2,x_3])^{{\leq}{\Delta} 
} $. } \medskip\\
{\bf Proof.}  Determinant  polylineary    
depends from  ${\Delta} $  first   columns, elements  which  are 
 $(\widetilde {{\sigma} }_0,\widetilde {{\sigma} }_1,\widetilde 
{{\sigma} }_2,\widetilde {{\sigma} }_3)$,   hence,   
polynomially is expressed via  $(\widetilde {{\sigma} }_0,\widetilde 
{{\sigma} }_1,\widetilde {{\sigma} }_2,\widetilde {{\sigma} }_3)$ 
by homogeneouse  polynomial of the degree   ${\Delta} $.   
Then $S(x_1,x_2,x_3)$ polynomially is expressed via 
$(\widetilde {{\sigma} }_0,\widetilde {{\sigma} }_1,\widetilde 
{{\sigma} }_2,\widetilde {{\sigma} }_3)$ 
by homogeneouse  polynomial of the degree ${\Delta} $. 
Since $\widetilde {{\sigma} }_0{=}1(x_1,x_2,x_3)$,  
then  $S(x_1,x_2,x_3)$ is polynomially 
 expressed via $(\widetilde {{\sigma} }_1,\widetilde {{\sigma} 
}_2,\widetilde {{\sigma} }_3)$ 
by polynomial of the degree ${\leq}{\Delta} $. \medskip
\\
Since  $\widetilde {{\sigma} }_0,\widetilde {{\sigma} }_1,\widetilde 
{{\sigma} }_2,\widetilde {{\sigma} }_3  \in 
\mathop{\sf ts}^r({\bf R}[x_1,x_2,x_3])^{{\leq}1} $,  
and $S(x_1,x_2,x_3)$ is polynomially 
expressed via $(\widetilde {{\sigma} }_0,\widetilde 
{{\sigma} }_1,\widetilde {{\sigma} }_2,\widetilde {{\sigma} }_3)$ 
by homogeneouse polynomial of the degree 
${\leq}{\Delta} $, then $S(x_1,x_2,x_3)$ $\in 
 (\mathop{\sf ts}^r({\bf R}[x_1,x_2,x_3])^{{\leq}1} )^{\Delta} $. 
\medskip\\
Since by of proposition 1.25 \medskip

 $(\mathop{\sf ts}^r({\bf R}[x_1,x_2,x_3])^{{\leq}1} )^{\Delta} 
 \subseteq  \mathop{\sf ts}^r({\bf R}[x_1,x_2,x_3]^{{\leq}{\Delta} 
} $, \medskip
\\
then $S(x_1,x_2,x_3) \in  
\mathop{\sf ts}^r({\bf R}[x_1,x_2,x_3]^{{\leq}{\Delta} } $. \bigskip  
\\ 
 \noindent {\bf Example 4.21.} (to theorem 2.1) {\it Let ${\Delta} 
{=}3$, $r{=}3$, $d{=}r{-}1$. \medskip\\
1) Let $h_i(x) \in  {\bf R}[x]^{{\leq}d{+}{\Delta} } $ 
for $i{=}1,r$, i. e. \bigskip

 $h_i(x) =  \sum\limits^{d{+}{\Delta} } _{{\delta} {=}0} 
h_{i,{\delta} } {\cdot}x^{\delta} $ for $i{=}1,r$. \bigskip\\
Let \bigskip

 $S(x_1,x_2,x_3) =  {\det}\left\| \begin{matrix}\overline {\sigma} 
_0 &   &   & h_{1,5}  & h_{2,5}  & h_{3,5}  \cr
\overline {\sigma} _1 & \overline {\sigma} _0 &   & h_{1,4} 
 & h_{2,4}  & h_{3,4}  \cr
\overline {\sigma} _2 & \overline {\sigma} _1 & \overline {\sigma} 
_0 & h_{1,3}  & h_{2,3}  & h_{3,3}  \cr
\overline {\sigma} _3 & \overline {\sigma} _2 & \overline {\sigma} 
_1 & h_{1,2}  & h_{2,2}  & h_{3,2}  \cr
  & \overline {\sigma} _3 & \overline {\sigma} _2 & h_{1,1} 
 & h_{2,1}  & h_{3,1}  \cr
  &   & \overline {\sigma} _3 & h_{1,0}  & h_{2,0}  & h_{3,0} 
\end{matrix} \right\| $, \bigskip\\
then \bigskip

 $S(x_1,x_2,x_3) \in  \mathop{\sf ts}^r({\bf R}[x_1,x_2,x_3])^{{\leq}{\Delta} 
} $, \bigskip

 ${\det}\left\| \begin{matrix}h_1(x_1) & h_2(x_1) & h_3(x_1) 
\cr
h_1(x_2) & h_2(x_2) & h_3(x_2) \cr
h_1(x_3) & h_2(x_3) & h_3(x_3)\end{matrix} \right\|  =  S(x_1,x_2,x_3){\cdot}{\det}\left\| 
\begin{matrix}x^2_1 & x^1_1 & x^0_1 \cr
x^2_2 & x^1_2 & x^0_2 \cr
x^2_3 & x^1_3 & x^0_3\end{matrix} \right\| $, \bigskip\\
 that is equivalent to \bigskip

 $S(x{\otimes}1{\otimes}1,1{\otimes}x{\otimes}1,1{\otimes}1{\otimes}x) 
\in  \mathop{\sf TS}^r({\bf R}[x])^{{\leq}{\Delta} } $, \bigskip

 $h_1(x){\wedge}h_2(x){\wedge}h_3(x) =  S(x{\otimes}1{\otimes}1,1{\otimes}x{\otimes}1,1{\otimes}1{\otimes}x){\cdot}(x^2{\wedge}x^1{\wedge}x^0)$. 
}
\hfil\vfill\eject
 
 \noindent {\bf Proof.} Let ${\Delta} '{=}1$, 
${\Delta} ''{=}2$, then ${\Delta} {=}{\Delta} '{+}{\Delta} 
''$. There holds\\

 $\overline {{\sigma} }^{{\Delta} ''} _0{\cdot}\overline 
{{\sigma} }^{{\Delta} '} _3{\cdot}S(x_1,x_2,x_3){\cdot}{\det}\left\| 
\begin{matrix}x^2_1 & x^1_1 & x^0_1 \cr
x^2_2 & x^1_2 & x^0_2 \cr
x^2_3 & x^1_3 & x^0_3\end{matrix} \right\| $\\
\\

 \qquad $=  \overline {{\sigma} }^{{\Delta} ''} _0{\cdot}\overline 
{{\sigma} }^{{\Delta} '} _3{\cdot}{\det}\left\| \begin{matrix}\overline 
{\sigma} _0 &   &   & h_{1,5}  & h_{2,5}  & h_{3,5}  \cr
\overline {\sigma} _1 & \overline {\sigma} _0 &   & h_{1,4} 
 & h_{2,4}  & h_{3,4}  \cr
\overline {\sigma} _2 & \overline {\sigma} _1 & \overline {\sigma} 
_0 & h_{1,3}  & h_{2,3}  & h_{3,3}  \cr
\overline {\sigma} _3 & \overline {\sigma} _2 & \overline {\sigma} 
_1 & h_{1,2}  & h_{2,2}  & h_{3,2}  \cr
  & \overline {\sigma} _3 & \overline {\sigma} _2 & h_{1,1} 
 & h_{2,1}  & h_{3,1}  \cr
  &   & \overline {\sigma} _3 & h_{1,0}  & h_{2,0}  & h_{3,0} 
\end{matrix} \right\| {\cdot}{\det}\left\| \begin{matrix}x^2_1 
& x^1_1 & x^0_1 \cr
x^2_2 & x^1_2 & x^0_2 \cr
x^2_3 & x^1_3 & x^0_3\end{matrix} \right\| $\\
\\

 \qquad $\buildrel{ 1}\over  = ({-}1)^{r{\cdot}{\Delta} 
'} {\cdot}{\det}\left\| \begin{matrix}\overline {\sigma} _0 & 
  &   & h_{1,5}  & h_{2,5}  & h_{3,5}  \cr
\overline {\sigma} _1 & \overline {\sigma} _0 &   & h_{1,4} 
 & h_{2,4}  & h_{3,4}  \cr
\overline {\sigma} _2 & \overline {\sigma} _1 & \overline {\sigma} 
_0 & h_{1,3}  & h_{2,3}  & h_{3,3}  \cr
\overline {\sigma} _3 & \overline {\sigma} _2 & \overline {\sigma} 
_1 & h_{1,2}  & h_{2,2}  & h_{3,2}  \cr
  & \overline {\sigma} _3 & \overline {\sigma} _2 & h_{1,1} 
 & h_{2,1}  & h_{3,1}  \cr
  &   & \overline {\sigma} _3 & h_{1,0}  & h_{2,0}  & h_{3,0} 
\end{matrix} \right\| {\cdot}{\det}\left\| \begin{matrix}x^3_1 
& x^2_1 & x^1_1 \cr
x^3_2 & x^2_2 & x^1_2 \cr
x^3_3 & x^2_3 & x^1_3\end{matrix} \right\| $\\
\\

 \qquad $=  ({-}1)^{r{\cdot}{\Delta} '} {\cdot}{\det}\left\| 
\begin{matrix}\overline {\sigma} _0 &   &   & h_{1,5}  & h_{2,5} 
 & h_{3,5}  &   &   &   \cr
\overline {\sigma} _1 & \overline {\sigma} _0 &   & h_{1,4} 
 & h_{2,4}  & h_{3,4}  &   &   &   \cr
\overline {\sigma} _2 & \overline {\sigma} _1 & \overline {\sigma} 
_0 & h_{1,3}  & h_{2,3}  & h_{3,3}  & 1 &   &   \cr
\overline {\sigma} _3 & \overline {\sigma} _2 & \overline {\sigma} 
_1 & h_{1,2}  & h_{2,2}  & h_{3,2}  &   & 1 &   \cr
  & \overline {\sigma} _3 & \overline {\sigma} _2 & h_{1,1} 
 & h_{2,1}  & h_{3,1}  &   &   & 1 \cr
  &   & \overline {\sigma} _3 & h_{1,0}  & h_{2,0}  & h_{3,0} 
 &   &   &   \cr
     &      &      &        &        &        & x^3_1 & x^2_1 
& x^1_1 \cr
     &      &      &        &        &        & x^3_2 & x^2_2 
& x^1_2 \cr
     &      &      &        &        &        & x^3_3 & x^2_3 
& x^1_3\end{matrix} \right\| $\\
\\

 \qquad $\buildrel{ 2}\over  = ({-}1)^{r{\cdot}{\Delta} 
'} {\cdot}{\det}\left\| \begin{matrix}\overline {\sigma} _0 & 
  &   &  h_{1,5}  &  h_{2,5}  &  h_{3,5}  &   &   &   \cr
\overline {\sigma} _1 & \overline {\sigma} _0 &   &  h_{1,4} 
 &  h_{2,4}  &  h_{3,4}  &   &   &   \cr
\overline {\sigma} _2 & \overline {\sigma} _1 & \overline {\sigma} 
_0 &  h_{1,3}  &  h_{2,3}  &  h_{3,3}  & 1 &   &   \cr
\overline {\sigma} _3 & \overline {\sigma} _2 & \overline {\sigma} 
_1 &  h_{1,2}  &  h_{2,2}  &  h_{3,2}  &   & 1 &   \cr
  & \overline {\sigma} _3 & \overline {\sigma} _2 &  h_{1,1} 
 &  h_{2,1}  &  h_{3,1}  &   &   & 1 \cr
  &   & \overline {\sigma} _3 &  h_{1,0}  &  h_{2,0}  &  h_{3,0} 
 &   &   &   \cr
     &      &      & {-}h_1(x_1) & {-}h_2(x_1) & {-}h_3(x_1) 
&   &   &   \cr
     &      &      & {-}h_1(x_2) & {-}h_2(x_2) & {-}h_3(x_2) 
&   &   &   \cr
     &      &      & {-}h_1(x_3) & {-}h_2(x_3) & {-}h_3(x_3) 
&   &   &  \end{matrix} \right\| $\\
\hfil\vfill\eject

 \qquad $=  ({-}1)^{r{\cdot}{\Delta} '} {\cdot}{\det}\left\| 
\begin{matrix}\overline {\sigma} _0 &   &   &  h_{1,5}  &  h_{2,5} 
 &  h_{3,5}  &   &   &   \cr
\overline {\sigma} _1 & \overline {\sigma} _0 &   &  h_{1,4} 
 &  h_{2,4}  &  h_{3,4}  &   &   &   \cr
\overline {\sigma} _2 & \overline {\sigma} _1 & \overline {\sigma} 
_0 &  h_{1,3}  &  h_{2,3}  &  h_{3,3}  & 1 &   &   \cr
\overline {\sigma} _3 & \overline {\sigma} _2 & \overline {\sigma} 
_1 &  h_{1,2}  &  h_{2,2}  &  h_{3,2}  &   & 1 &   \cr
  & \overline {\sigma} _3 & \overline {\sigma} _2 &  h_{1,1} 
 &  h_{2,1}  &  h_{3,1}  &   &   & 1 \cr
  &   & \overline {\sigma} _3 &  h_{1,0}  &  h_{2,0}  &  h_{3,0} 
 &   &   &   \cr
     &      &      & {-}h_1(x_1) & {-}h_2(x_1) & {-}h_3(x_1) 
&   &   &   \cr
     &      &      & {-}h_1(x_2) & {-}h_2(x_2) & {-}h_3(x_2) 
&   &   &   \cr
     &      &      & {-}h_1(x_3) & {-}h_2(x_3) & {-}h_3(x_3) 
&   &   &  \end{matrix} \right\| $\\
\\

 \qquad $=  ({-}1)^r{\cdot}{\det}\left\| \begin{matrix}\overline 
{\sigma} _0 &   &   &   &   &   &  h_{1,5}  &  h_{2,5}  &  h_{3,5} 
 \cr
\overline {\sigma} _1 & \overline {\sigma} _0 &   &   &   & 
  &  h_{1,4}  &  h_{2,4}  &  h_{3,4}  \cr
\overline {\sigma} _2 & \overline {\sigma} _1 & 1 &   &   & 
\overline {\sigma} _0 &  h_{1,3}  &  h_{2,3}  &  h_{3,3}  \cr
\overline {\sigma} _3 & \overline {\sigma} _2 &   & 1 &   & 
\overline {\sigma} _1 &  h_{1,2}  &  h_{2,2}  &  h_{3,2}  \cr
  & \overline {\sigma} _3 &   &   & 1 & \overline {\sigma} 
_2 &  h_{1,1}  &  h_{2,1}  &  h_{3,1}  \cr
  &   &   &   &   & \overline {\sigma} _3 &  h_{1,0}  &  h_{2,0} 
 &  h_{3,0}  \cr
     &      &   &   &   &      & {-}h_1(x_1) & {-}h_2(x_1) 
& {-}h_3(x_1) \cr
     &      &   &   &   &      & {-}h_1(x_2) & {-}h_2(x_2) 
& {-}h_3(x_2) \cr
     &      &   &   &   &      & {-}h_1(x_3) & {-}h_2(x_3) 
& {-}h_3(x_3)\end{matrix} \right\| $\\
\\

 \qquad $=  ({-}1)^r{\cdot}{\det}\left\| \begin{matrix}\overline 
{\sigma} _0 &      &     &     &     &   \cr
\overline {\sigma} _1 & \overline {\sigma} _0 &   &   &   & 
  \cr
\overline {\sigma} _2 & \overline {\sigma} _1 & 1 &   &   & 
\overline {\sigma} _0 \cr
\overline {\sigma} _3 & \overline {\sigma} _2 &   & 1 &   & 
\overline {\sigma} _1 \cr
  & \overline {\sigma} _3 &   &   & 1 & \overline {\sigma} 
_2 \cr
  &   &     &     &     & \overline {\sigma} _3\end{matrix} 
\right\| {\cdot}{\det}\left\| \begin{matrix}{-}h_1(x_1) & {-}h_2(x_1) 
& {-}h_3(x_1) \cr
{-}h_1(x_2) & {-}h_2(x_2) & {-}h_3(x_2) \cr
{-}h_1(x_3) & {-}h_2(x_3) & {-}h_3(x_3)\end{matrix} \right\| 
$\\
\\

 \qquad $\buildrel{ 3}\over  = ({-}1)^r{\cdot}\overline {{\sigma} 
}^{{\Delta} ''} _0{\cdot}\overline {{\sigma} }^{{\Delta} 
'} _3{\cdot}{\det}\left\| \begin{matrix}{-}h_1(x_1) & {-}h_2(x_1) 
& {-}h_3(x_1) \cr
{-}h_1(x_2) & {-}h_2(x_2) & {-}h_3(x_2) \cr
{-}h_1(x_3) & {-}h_2(x_3) & {-}h_3(x_3)\end{matrix} \right\| 
$\\
\\

 \qquad $=  \overline {{\sigma} }^{{\Delta} ''} _0{\cdot}\overline 
{{\sigma} }^{{\Delta} '} _3{\cdot}{\det}\left\| \begin{matrix}h_1(x_1) 
& h_2(x_1) & h_3(x_1) \cr
h_1(x_2) & h_2(x_2) & h_3(x_2) \cr
h_1(x_3) & h_2(x_3) & h_3(x_3)\end{matrix} \right\| $.\\
\\
\\
Dividing both parts of the obtained equality by the common multiple 
we obtain the equality\\

 $S(x_1,x_2,x_3){\cdot}{\det}\left\| \begin{matrix}x^2_1 & 
x^1_1 & x^0_1 \cr
x^2_2 & x^1_2 & x^0_2 \cr
x^2_3 & x^1_3 & x^0_3\end{matrix} \right\|  =  {\det}\left\| 
\begin{matrix}h_1(x_1) & h_2(x_1) & h_3(x_1) \cr
h_1(x_2) & h_2(x_2) & h_3(x_2) \cr
h_1(x_3) & h_2(x_3) & h_3(x_3)\end{matrix} \right\| $.\\
\hfil\vfill\eject
 
 \noindent {\bf Example 4.22.} (to theorem 2.3) {\it Let ${\Delta} 
{=}3$, $r{=}3$, $d{=}3$. \medskip\\
1. Let $h_i(x) \in  {\bf R}[x]^{{\leq}d{+}{\Delta} } $ 
for $i{=}1,r$, i. e. \bigskip

 $h_i(x) =  \sum\limits^{d{+}{\Delta} } _{{\delta} {=}0} 
h_{i,{\delta} } {\cdot}x^{\delta} $ for $i{=}1,r$, \medskip\\
Then \bigskip

 $h_1(x){\wedge}h_2(x){\wedge}h_3(x)$  \bigskip

 \qquad $=  ({-}1)^r{\cdot}({-}1)^{r{\cdot}{\Delta} '} {\cdot}{\det}\left\| 
\begin{matrix}  & \tilde {\sigma} _0 &      &   &   &   &   & 
     & h_{1,6}  & h_{2,6}  & h_{3,6}  \cr
  & \tilde {\sigma} _1 & \tilde {\sigma} _0 &   &   &   &  
 &      & h_{1,5}  & h_{2,5}  & h_{3,5}  \cr
  & \tilde {\sigma} _2 & \tilde {\sigma} _1 & 1 &   &   &  
 &      & h_{1,4}  & h_{2,4}  & h_{3,4}  \cr
  & \tilde {\sigma} _3 & \tilde {\sigma} _2 &   & 1 &   &  
 & \tilde {\sigma} _0 & h_{1,3}  & h_{2,3}  & h_{3,3}  \cr
  &      & \tilde {\sigma} _3 &   &   & 1 &   & \tilde {\sigma} 
_1 & h_{1,2}  & h_{2,2}  & h_{3,2}  \cr
  &      &      &   &   &   & 1 & \tilde {\sigma} _2 & h_{1,1} 
 & h_{2,1}  & h_{3,1}  \cr
    &      &      &   &   &   &   & \tilde {\sigma} _3 & h_{1,0} 
 & h_{2,0}  & h_{3,0}  \cr
\wedge  &      &      & x^3 & x^2 & x^1 & x^0 &      &     
   &        &  \end{matrix} \right\| $.\\
\\
\\
This equality is equivalent to \\

 ${\det}\left\| \begin{matrix}h_1(x_1) & h_2(x_1) & h_3(x_1) 
\cr
h_1(x_2) & h_2(x_2) & h_3(x_2) \cr
h_1(x_3) & h_2(x_3) & h_3(x_3)\end{matrix} \right\| $\\
\\

 \qquad $=  ({-}1)^r{\cdot}({-}1)^{r{\cdot}{\Delta} '} {\cdot}{\det}\left\| 
\begin{matrix}\overline {\sigma} _0 &      &   &   &   &   & 
     & h_{1,6}  & h_{2,6}  & h_{3,6}  \cr
\overline {\sigma} _1 & \overline {\sigma} _0 &   &   &   & 
  &      & h_{1,5}  & h_{2,5}  & h_{3,5}  \cr
\overline {\sigma} _2 & \overline {\sigma} _1 & 1 &   &   & 
  &      & h_{1,4}  & h_{2,4}  & h_{3,4}  \cr
\overline {\sigma} _3 & \overline {\sigma} _2 &   & 1 &   & 
  & \overline {\sigma} _0 & h_{1,3}  & h_{2,3}  & h_{3,3}  \cr
     & \overline {\sigma} _3 &   &   & 1 &   & \overline {\sigma} 
_1 & h_{1,2}  & h_{2,2}  & h_{3,2}  \cr
     &      &   &   &   & 1 & \overline {\sigma} _2 & h_{1,1} 
 & h_{2,1}  & h_{3,1}  \cr
     &      &   &   &   &   & \overline {\sigma} _3 & h_{1,0} 
 & h_{2,0}  & h_{3,0}  \cr
     &      & x^3_1 & x^2_1 & x^1_1 & x^0_1 &      &       
 &        &   \cr
     &      & x^3_2 & x^2_2 & x^1_2 & x^0_2 &      &       
 &        &   \cr
     &      & x^3_3 & x^2_3 & x^1_3 & x^0_3 &      &       
 &        &  \end{matrix} \right\| $.\\
\\
2. There holds \bigskip

 $\bigwedge  ^r({\bf R}[x])^{{\leq}d{+}{\Delta} }  \subseteq 
 \mathop{\sf TS}({\bf R}[x])^{{\leq}{\Delta} } {\cdot}\bigwedge 
 ^r({\bf R}[x])^{{\leq}d} $. }\\
\hfil\vfill\eject
 
 \noindent {\bf Proof 1.} Let ${\Delta} '{=}1$, 
${\Delta} ''{=}2$, then ${\Delta} {=}{\Delta} '{+}{\Delta} 
''$. There holds\\

 $({-}1)^{{-}r{\cdot}{\Delta} '} {\cdot}\overline {{\sigma} 
}^{{\Delta} ''} _0{\cdot}\overline {{\sigma} }^{{\Delta} 
'} _3{\cdot}{\det}\left\| \begin{matrix}\overline {\sigma} _0 
&      &   &   &   &   &      & h_{1,6}  & h_{2,6}  & h_{3,6} 
 \cr
\overline {\sigma} _1 & \overline {\sigma} _0 &   &   &   & 
  &      & h_{1,5}  & h_{2,5}  & h_{3,5}  \cr
\overline {\sigma} _2 & \overline {\sigma} _1 & 1 &   &   & 
  &      & h_{1,4}  & h_{2,4}  & h_{3,4}  \cr
\overline {\sigma} _3 & \overline {\sigma} _2 &   & 1 &   & 
  & \overline {\sigma} _0 & h_{1,3}  & h_{2,3}  & h_{3,3}  \cr
     & \overline {\sigma} _3 &   &   & 1 &   & \overline {\sigma} 
_1 & h_{1,2}  & h_{2,2}  & h_{3,2}  \cr
     &      &   &   &   & 1 & \overline {\sigma} _2 & h_{1,1} 
 & h_{2,1}  & h_{3,1}  \cr
     &      &   &   &   &   & \overline {\sigma} _3 & h_{1,0} 
 & h_{2,0}  & h_{3,0}  \cr
     &      & x^3_1 & x^2_1 & x^1_1 & x^0_1 &      &       
 &        &   \cr
     &      & x^3_2 & x^2_2 & x^1_2 & x^0_2 &      &       
 &        &   \cr
     &      & x^3_3 & x^2_3 & x^1_3 & x^0_3 &      &       
 &        &  \end{matrix} \right\| $\\
\\

 \qquad $\buildrel{ 1}\over  = {\det}\left\| \begin{matrix}\overline 
{\sigma} _0 &      &   &   &   &   &      & h_{1,6}  & h_{2,6} 
 & h_{3,6}  \cr
\overline {\sigma} _1 & \overline {\sigma} _0 &   &   &   & 
  &      & h_{1,5}  & h_{2,5}  & h_{3,5}  \cr
\overline {\sigma} _2 & \overline {\sigma} _1 & 1 &   &   & 
  &      & h_{1,4}  & h_{2,4}  & h_{3,4}  \cr
\overline {\sigma} _3 & \overline {\sigma} _2 &   & 1 &   & 
  & \overline {\sigma} _0 & h_{1,3}  & h_{2,3}  & h_{3,3}  \cr
     & \overline {\sigma} _3 &   &   & 1 &   & \overline {\sigma} 
_1 & h_{1,2}  & h_{2,2}  & h_{3,2}  \cr
     &      &   &   &   & 1 & \overline {\sigma} _2 & h_{1,1} 
 & h_{2,1}  & h_{3,1}  \cr
     &      &   &   &   &   & \overline {\sigma} _3 & h_{1,0} 
 & h_{2,0}  & h_{3,0}  \cr
     &      & x^4_1 & x^3_1 & x^2_1 & x^1_1 &      &       
 &        &   \cr
     &      & x^4_2 & x^3_2 & x^2_2 & x^1_2 &      &       
 &        &   \cr
     &      & x^4_3 & x^3_3 & x^2_3 & x^1_3 &      &       
 &        &  \end{matrix} \right\| $\\
\\

 \qquad $\buildrel{ 2}\over  = {\det}\left\| \begin{matrix}\overline 
{\sigma} _0 &      &   &   &   &   &      &  h_{1,6}  &  h_{2,6} 
 &  h_{3,6}  \cr
\overline {\sigma} _1 & \overline {\sigma} _0 &   &   &   & 
  &      &  h_{1,5}  &  h_{2,5}  &  h_{3,5}  \cr
\overline {\sigma} _2 & \overline {\sigma} _1 & 1 &   &   & 
  &      &  h_{1,4}  &  h_{2,4}  &  h_{3,4}  \cr
\overline {\sigma} _3 & \overline {\sigma} _2 &   & 1 &   & 
  & \overline {\sigma} _0 &  h_{1,3}  &  h_{2,3}  &  h_{3,3} 
 \cr
     & \overline {\sigma} _3 &   &   & 1 &   & \overline {\sigma} 
_1 &  h_{1,2}  &  h_{2,2}  &  h_{3,2}  \cr
     &      &   &   &   & 1 & \overline {\sigma} _2 &  h_{1,1} 
 &  h_{2,1}  &  h_{3,1}  \cr
     &      &   &   &   &   & \overline {\sigma} _3 &  h_{1,0} 
 &  h_{2,0}  &  h_{3,0}  \cr
     &      &   &   &   &   &      & {-}h_1(x_1) & {-}h_2(x_1) 
& {-}h_3(x_1) \cr
     &      &   &   &   &   &      & {-}h_1(x_2) & {-}h_2(x_2) 
& {-}h_3(x_2) \cr
     &      &   &   &   &   &      & {-}h_1(x_3) & {-}h_2(x_3) 
& {-}h_3(x_3)\end{matrix} \right\| $\\
\\

 \qquad $=  {\det}\left\| \begin{matrix}\overline {\sigma} 
_0 &      &   &   &   &   &   \cr
\overline {\sigma} _1 & \overline {\sigma} _0 &   &   &   & 
  &   \cr
\overline {\sigma} _2 & \overline {\sigma} _1 & 1 &   &   & 
  &   \cr
\overline {\sigma} _3 & \overline {\sigma} _2 &   & 1 &   & 
  & \overline {\sigma} _0 \cr
     & \overline {\sigma} _3 &   &   & 1 &   & \overline {\sigma} 
_1 \cr
     &      &   &   &   & 1 & \overline {\sigma} _2 \cr
     &      &   &   &   &   & \overline {\sigma} _3\end{matrix} 
\right\| {\cdot}{\det}\left\| \begin{matrix}{-}h_1(x_1) & {-}h_2(x_1) 
& {-}h_3(x_1) \cr
{-}h_1(x_2) & {-}h_2(x_2) & {-}h_3(x_2) \cr
{-}h_1(x_3) & {-}h_2(x_3) & {-}h_3(x_3)\end{matrix} \right\| 
$\\
\hfil\vfill\eject

 \qquad $=  {\det}\left\| \begin{matrix}\overline {\sigma} 
_0 &      &   &   &   &   &   \cr
\overline {\sigma} _1 & \overline {\sigma} _0 &   &   &   & 
  &   \cr
\overline {\sigma} _2 & \overline {\sigma} _1 & 1 &   &   & 
  &   \cr
\overline {\sigma} _3 & \overline {\sigma} _2 &   & 1 &   & 
  & \overline {\sigma} _0 \cr
     & \overline {\sigma} _3 &   &   & 1 &   & \overline {\sigma} 
_1 \cr
     &      &   &   &   & 1 & \overline {\sigma} _2 \cr
     &      &   &   &   &   & \overline {\sigma} _3\end{matrix} 
\right\| {\cdot}{\det}\left\| \begin{matrix}{-}h_1(x_1) & {-}h_2(x_1) 
& {-}h_3(x_1) \cr
{-}h_1(x_2) & {-}h_2(x_2) & {-}h_3(x_2) \cr
{-}h_1(x_3) & {-}h_2(x_3) & {-}h_3(x_3)\end{matrix} \right\| 
$\\
\\

 \qquad $\buildrel{ 3}\over  = \overline {{\sigma} }^{{\Delta} 
''} _0{\cdot}\overline {{\sigma} }^{{\Delta} '} _3{\cdot}{\det}\left\| 
\begin{matrix}{-}h_1(x_1) & {-}h_2(x_1) & {-}h_3(x_1) \cr
{-}h_1(x_2) & {-}h_2(x_2) & {-}h_3(x_2) \cr
{-}h_1(x_3) & {-}h_2(x_3) & {-}h_3(x_3)\end{matrix} \right\| 
$\\
\\

 \qquad $=  \overline {{\sigma} }^{{\Delta} ''} _0{\cdot}\overline 
{{\sigma} }^{{\Delta} '} _3{\cdot}({-}1)^r{\cdot}{\det}\left\| 
\begin{matrix}h_1(x_1) & h_2(x_1) & h_3(x_1) \cr
h_1(x_2) & h_2(x_2) & h_3(x_2) \cr
h_1(x_3) & h_2(x_3) & h_3(x_3)\end{matrix} \right\| $.\\
\\
Dividing both parts of the obtained equality by the common multiple 
we obtain the equality\\
\\

 ${\det}\left\| \begin{matrix}\overline {\sigma} _0 &      
&   &   &   &   &      & h_{1,6}  & h_{2,6}  & h_{3,6}  \cr
\overline {\sigma} _1 & \overline {\sigma} _0 &   &   &   & 
  &      & h_{1,5}  & h_{2,5}  & h_{3,5}  \cr
\overline {\sigma} _2 & \overline {\sigma} _1 & 1 &   &   & 
  &      & h_{1,4}  & h_{2,4}  & h_{3,4}  \cr
\overline {\sigma} _3 & \overline {\sigma} _2 &   & 1 &   & 
  & \overline {\sigma} _0 & h_{1,3}  & h_{2,3}  & h_{3,3}  \cr
     & \overline {\sigma} _3 &   &   & 1 &   & \overline {\sigma} 
_1 & h_{1,2}  & h_{2,2}  & h_{3,2}  \cr
     &      &   &   &   & 1 & \overline {\sigma} _2 & h_{1,1} 
 & h_{2,1}  & h_{3,1}  \cr
     &      &   &   &   &   & \overline {\sigma} _3 & h_{1,0} 
 & h_{2,0}  & h_{3,0}  \cr
     &      & x^3_1 & x^2_1 & x^1_1 & x^0_1 &      &       
 &        &   \cr
     &      & x^3_2 & x^2_2 & x^1_2 & x^0_2 &      &       
 &        &   \cr
     &      & x^3_3 & x^2_3 & x^1_3 & x^0_3 &      &       
 &        &  \end{matrix} \right\| $\\
\\

 \qquad $=  ({-}1)^{r{\cdot}{\Delta} '} {\cdot}({-}1)^r{\cdot}{\det}\left\| 
\begin{matrix}h_1(x_1) & h_2(x_1) & h_3(x_1) \cr
h_1(x_2) & h_2(x_2) & h_3(x_2) \cr
h_1(x_3) & h_2(x_3) & h_3(x_3)\end{matrix} \right\| $.\\
\hfil\vfill\eject
 
 \noindent {\bf Proof 2.} There holds\\

 $h_1(x){\wedge}h_2(x){\wedge}h_3(x)$\\

 \qquad $=  ({-}1)^r{\cdot}({-}1)^{r{\cdot}{\Delta} '} {\cdot}{\det}\left\| 
\begin{matrix}  & \tilde {\sigma} _0 &      &   &   &   &   & 
     & h_{1,6}  & h_{2,6}  & h_{3,6}  \cr
  & \tilde {\sigma} _1 & \tilde {\sigma} _0 &   &   &   &  
 &      & h_{1,5}  & h_{2,5}  & h_{3,5}  \cr
  & \tilde {\sigma} _2 & \tilde {\sigma} _1 & 1 &   &   &  
 &      & h_{1,4}  & h_{2,4}  & h_{3,4}  \cr
  & \tilde {\sigma} _3 & \tilde {\sigma} _2 &   & 1 &   &  
 & \tilde {\sigma} _0 & h_{1,3}  & h_{2,3}  & h_{3,3}  \cr
  &      & \tilde {\sigma} _3 &   &   & 1 &   & \tilde {\sigma} 
_1 & h_{1,2}  & h_{2,2}  & h_{3,2}  \cr
  &      &      &   &   &   & 1 & \tilde {\sigma} _2 & h_{1,1} 
 & h_{2,1}  & h_{3,1}  \cr
    &      &      &   &   &   &   & \tilde {\sigma} _3 & h_{1,0} 
 & h_{2,0}  & h_{3,0}  \cr
\wedge  &      &      & x^3 & x^2 & x^1 & x^0 &      &     
   &        &  \end{matrix} \right\| $\\
\\
\\

 \qquad $=  {+}{\det}\left\| \begin{matrix}\tilde {\sigma} 
_0 &      &   &      & h_{1,6}  & h_{2,6}  & h_{3,6}  \cr
\tilde {\sigma} _1 & \tilde {\sigma} _0 &   &      & h_{1,5} 
 & h_{2,5}  & h_{3,5}  \cr
\tilde {\sigma} _2 & \tilde {\sigma} _1 & 1 &      & h_{1,4} 
 & h_{2,4}  & h_{3,4}  \cr
\tilde {\sigma} _3 & \tilde {\sigma} _2 &   & \tilde {\sigma} 
_0 & h_{1,3}  & h_{2,3}  & h_{3,3}  \cr
     & \tilde {\sigma} _3 &   & \tilde {\sigma} _1 & h_{1,2} 
 & h_{2,2}  & h_{3,2}  \cr
     &      &   & \tilde {\sigma} _2 & h_{1,1}  & h_{2,1}  
& h_{3,1}  \cr
     &      &   & \tilde {\sigma} _3 & h_{1,0}  & h_{2,0}  
& h_{3,0} \end{matrix} \right\| {\cdot}(x^2{\wedge}x^1{\wedge}x^0)$
\\
\\

 \qquad \qquad ${-}{\det}\left\| \begin{matrix}\tilde {\sigma} 
_0 &      &   &      & h_{1,6}  & h_{2,6}  & h_{3,6}  \cr
\tilde {\sigma} _1 & \tilde {\sigma} _0 &   &      & h_{1,5} 
 & h_{2,5}  & h_{3,5}  \cr
\tilde {\sigma} _2 & \tilde {\sigma} _1 &   &      & h_{1,4} 
 & h_{2,4}  & h_{3,4}  \cr
\tilde {\sigma} _3 & \tilde {\sigma} _2 & 1 & \tilde {\sigma} 
_0 & h_{1,3}  & h_{2,3}  & h_{3,3}  \cr
     & \tilde {\sigma} _3 &   & \tilde {\sigma} _1 & h_{1,2} 
 & h_{2,2}  & h_{3,2}  \cr
     &      &   & \tilde {\sigma} _2 & h_{1,1}  & h_{2,1}  
& h_{3,1}  \cr
     &      &   & \tilde {\sigma} _3 & h_{1,0}  & h_{2,0}  
& h_{3,0} \end{matrix} \right\| {\cdot}(x^3{\wedge}x^1{\wedge}x^0)$
\\
\\

 \qquad \qquad ${+}{\det}\left\| \begin{matrix}\tilde {\sigma} 
_0 &      &   &      & h_{1,6}  & h_{2,6}  & h_{3,6}  \cr
\tilde {\sigma} _1 & \tilde {\sigma} _0 &   &      & h_{1,5} 
 & h_{2,5}  & h_{3,5}  \cr
\tilde {\sigma} _2 & \tilde {\sigma} _1 &   &      & h_{1,4} 
 & h_{2,4}  & h_{3,4}  \cr
\tilde {\sigma} _3 & \tilde {\sigma} _2 &   & \tilde {\sigma} 
_0 & h_{1,3}  & h_{2,3}  & h_{3,3}  \cr
     & \tilde {\sigma} _3 & 1 & \tilde {\sigma} _1 & h_{1,2} 
 & h_{2,2}  & h_{3,2}  \cr
     &      &   & \tilde {\sigma} _2 & h_{1,1}  & h_{2,1}  
& h_{3,1}  \cr
     &      &   & \tilde {\sigma} _3 & h_{1,0}  & h_{2,0}  
& h_{3,0} \end{matrix} \right\| {\cdot}(x^3{\wedge}x^2{\wedge}x^0)$
\\
\\

 \qquad \qquad ${-}{\det}\left\| \begin{matrix}\tilde {\sigma} 
_0 &      &   &      & h_{1,6}  & h_{2,6}  & h_{3,6}  \cr
\tilde {\sigma} _1 & \tilde {\sigma} _0 &   &      & h_{1,5} 
 & h_{2,5}  & h_{3,5}  \cr
\tilde {\sigma} _2 & \tilde {\sigma} _1 &   &      & h_{1,4} 
 & h_{2,4}  & h_{3,4}  \cr
\tilde {\sigma} _3 & \tilde {\sigma} _2 &   & \tilde {\sigma} 
_0 & h_{1,3}  & h_{2,3}  & h_{3,3}  \cr
     & \tilde {\sigma} _3 &   & \tilde {\sigma} _1 & h_{1,2} 
 & h_{2,2}  & h_{3,2}  \cr
     &      & 1 & \tilde {\sigma} _2 & h_{1,1}  & h_{2,1}  
& h_{3,1}  \cr
     &      &   & \tilde {\sigma} _3 & h_{1,0}  & h_{2,0}  
& h_{3,0} \end{matrix} \right\| {\cdot}(x^3{\wedge}x^2{\wedge}x^1)$.\bigskip 
\\
Here $(x^2{\wedge}x^1{\wedge}x^0)$, $(x^3{\wedge}x^1{\wedge}x^0)$, 
$(x^3{\wedge}x^2{\wedge}x^0)$, $(x^3{\wedge}x^1{\wedge}x^1)  
\in  \bigwedge  ^3({\bf R}[x])^{{\leq}3}   =  \bigwedge  ^r({\bf 
R}[x])^{{\leq}d} $. Coefficients of their is determinants, which $\in   (\mathop{\sf 
TS}({\bf R}[x])^{{\leq}1} )^{\Delta}   \subseteq   \mathop{\sf 
TS}({\bf R}[x])^{{\leq}{\Delta} } $,  since  $\widetilde 
{{\sigma} }_0,\widetilde {{\sigma} }_1,\widetilde {{\sigma} }_2,\widetilde 
{{\sigma} }_3  \in $ $\mathop{\sf TS}({\bf R}[x])^{{\leq}1} $, 
in all determinants the number of columns  with $\widetilde {{\sigma} }_0,\widetilde {{\sigma} 
}_1,\widetilde {{\sigma} }_2,\widetilde {{\sigma} }_3$ is equal to 
${\Delta} ''{+}{\Delta} '{=}{\Delta} $, inclusion $\subseteq 
$ holds  by of proposition 1.24.  There holds  $h_1(x){\wedge}h_2(x){\wedge}h_3(x) 
 \in   \bigwedge  ^r({\bf R}[x])^{{\leq}d{+}{\Delta} } $, 
 since $h_1(x),h_2(x),h_3(x){\in}{\bf R}[x]^{{\leq}d{+}{\Delta} 
} $, $r{=}3$. Hence, \bigskip

 $\bigwedge  ^r({\bf R}[x])^{{\leq}d{+}{\Delta} }  \subseteq 
 \mathop{\sf TS}({\bf R}[x])^{{\leq}{\Delta} } {\cdot}\bigwedge 
 ^r({\bf R}[x])^{{\leq}d} $, \bigskip\\
since $\bigwedge  ^r({\bf R}[x])^{{\leq}d{+}{\Delta} } 
$ lineary generated by elements  of the form  $h_1(x){\wedge}h_2(x){\wedge}h_3(x)$, 
 where $h_1(x)$, $h_2(x)$, $h_3(x) \in  {\bf R}[x]^{{\leq}d{+}{\Delta} 
} $ for $r{=}3$.
\hfil\vfill\eject

 \noindent {\bf Example 4.23.} (to theorem 2.2)  Let ${\Delta} {=}3$, 
$r{=}3$, $d{=}r{-}1$,\bigskip 

 $F(x) =  F_2{\cdot}x^2{+}F_1{\cdot}x^1{+}F_0{\cdot}x^0$, \bigskip

 $S(x_1,x_2,x_3) =  F(x_1){\cdot}F(x_2){\cdot}F(x_3)$. \bigskip\\
Tnen \bigskip

 $S(x{\otimes}1{\otimes}1,1{\otimes}x{\otimes}1,1{\otimes}1{\otimes}x){\cdot}(x^2{\wedge}x^1{\wedge}x^0) 
=  (F(x){\otimes}F(x){\otimes}F(x)){\cdot}(x^2{\wedge}x^1{\wedge}x^0)$ \bigskip

 \qquad $=  (F(x){\cdot}x^2){\wedge}(F(x){\cdot}x^1{\wedge}(F(x){\cdot}x^0)$,
\\

 $S(x_1,x_2,x_3) =  {\det}\left\| \begin{matrix}{\sigma} _0(x_1,x_2,x_3) 
&   & F_2 &      &   \cr
{\sigma} _1(x_1,x_2,x_3) & {\sigma} _0(x_1,x_2,x_3) & F_1 & 
F_2 &   \cr
{\sigma} _2(x_1,x_2,x_3) & {\sigma} _1(x_1,x_2,x_3) & F_0 & 
F_1 & F_2 \cr
{\sigma} _3(x_1,x_2,x_3) & {\sigma} _2(x_1,x_2,x_3) &      
& F_0 & F_1 \cr
  & {\sigma} _3(x_1,x_2,x_3) &      &      & F_0\end{matrix} 
\right\|  =  \mathop{\rm Res}({\sigma} (x_1,x_2,x_3),f)$.\\
\\
\\
Let \bigskip

 $f(x) =  (x{-}{\lambda} _1){\cdot}(x{-}{\lambda} _2){\cdot}(x{-}{\lambda} _3) =  f_3{\cdot}x^3{+}f_2{\cdot}x^2{+}f_1{\cdot}x^1{+}f_0{\cdot}x^0$.  
\bigskip
\\
Then \bigskip

 $f(x) =  {\sigma} _0({\lambda} _1,{\lambda} _2,{\lambda} _3){\cdot}x^3{+}{\sigma} 
_1({\lambda} _1,{\lambda} _2,{\lambda} _3){\cdot}x^2{+}{\sigma} 
_2({\lambda} _1,{\lambda} _2,{\lambda} _3){\cdot}x^1{+}{\sigma} 
_3({\lambda} _1,{\lambda} _2,{\lambda} _3){\cdot}x^0$, \bigskip
\\
hence, \bigskip

 ${\sigma} _0({\lambda} _1,{\lambda} _2,{\lambda} _3) =  f_3$, 
${\sigma} _1({\lambda} _1,{\lambda} _2,{\lambda} _3) =  f_2$,  
${\sigma} _2({\lambda} _1,{\lambda} _2,{\lambda} _3) =  f_1$, 
${\sigma} _3({\lambda} _1,{\lambda} _2,{\lambda} _3) =  f_0$, 
\bigskip

 ${\sigma} ({\lambda} _1,{\lambda} _2,{\lambda} _3)(x) =  {\sigma} 
_0({\lambda} _1,{\lambda} _2,{\lambda} _3){\cdot}x^3{+}{\sigma} 
_0({\lambda} _1,{\lambda} _2,{\lambda} _3){\cdot}x^2{+}{\sigma} 
_0({\lambda} _1,{\lambda} _2,{\lambda} _3){\cdot}x^1{+}{\sigma} 
_0({\lambda} _1,{\lambda} _2,{\lambda} _3){\cdot}x^0$  \bigskip

 \qquad $=  f_3{\cdot}x^3{+}f_2{\cdot}x^2{+}f_1{\cdot}x^1{+}f_0{\cdot}x^0 =  f(x)$. \bigskip\\
We have \bigskip

 $S({\lambda} _1,{\lambda} _2,{\lambda} _3) =  F({\lambda} 
_1){\cdot}F({\lambda} _2){\cdot}F({\lambda} _3)$,\\
\\

 $S({\lambda} _1,{\lambda} _2,{\lambda} _3) =  {\det}\left\| 
\begin{matrix}{\sigma} _0({\lambda} _1,{\lambda} _2,{\lambda} 
_3) &   & F_2 &      &   \cr
{\sigma} _1({\lambda} _1,{\lambda} _2,{\lambda} _3) & {\sigma} 
_0({\lambda} _1,{\lambda} _2,{\lambda} _3) & F_1 & F_2 &   \cr
{\sigma} _2({\lambda} _1,{\lambda} _2,{\lambda} _3) & {\sigma} 
_1({\lambda} _1,{\lambda} _2,{\lambda} _3) & F_0 & F_1 & F_2 
\cr
{\sigma} _3({\lambda} _1,{\lambda} _2,{\lambda} _3) & {\sigma} 
_2({\lambda} _1,{\lambda} _2,{\lambda} _3) &      & F_0 & F_1 
\cr
  & {\sigma} _3({\lambda} _1,{\lambda} _2,{\lambda} _3) &  
    &      & F_0\end{matrix} \right\| $\\
\\

 \qquad $=  {\det}\left\| \begin{matrix}f_3 &   & F_2 &    
  &   \cr
f_2 & f_0 & F_1 & F_2 &   \cr
f_1 & f_1 & F_0 & F_1 & F_2 \cr
f_0 & f_2 &      & F_0 & F_1 \cr
  & f_3 &      &      & F_0\end{matrix} \right\|  =  \mathop{\rm 
Res}({\sigma} ({\lambda} _1,{\lambda} _2,{\lambda} _3),F) =  
\mathop{\rm Res}(f,F)$.

\hfil\vfill\eject
 \noindent {\bf Example 4.24.} (Induction step $E(4) \Rightarrow 
 E(5)$ in proof of 2 of lemma 3.1).\\

 ${\det}\left\| \begin{matrix}x^4_1 & x^3_1 & x^2_1 & x^1_1 
& x^0_1 \cr
x^4_2 & x^3_2 & x^2_2 & x^1_2 & x^0_2 \cr
x^4_3 & x^3_3 & x^2_3 & x^1_3 & x^0_3 \cr
x^4_4 & x^3_4 & x^2_4 & x^1_4 & x^0_4 \cr
x^4_5 & x^3_5 & x^2_5 & x^1_5 & x^0_5\end{matrix} \right\| 
 =  {\det}\left\| \begin{matrix}x^4_1 & x^3_1 & x^2_1 & x^1_1 
& x^0_1 \cr
x^4_2 & x^3_2 & x^2_2 & x^1_2 & x^0_2 \cr
x^4_3 & x^3_3 & x^2_3 & x^1_3 & x^0_3 \cr
x^4_4 & x^3_4 & x^2_4 & x^1_4 & x^0_4 \cr
x^4_5 & x^3_5 & x^2_5 & x^1_5 & x^0_5\end{matrix} \right\| 
{\cdot}\left\| \begin{matrix}1\h &       &       &       &   
\cr
{-}x_5 & 1\h &       &       &   \cr
      & {-}x_5 & 1\h &       &   \cr
      &       & {-}x_5 & 1\h &   \cr
      &       &       & {-}x_5 & 1\end{matrix} \right\| $\\
\\

 \qquad $=  {\det}\left\| \begin{matrix}x^4_1{-}x_5{\cdot}x^3_1 
& x^3_1{-}x_5{\cdot}x^2_1 & x^2_1{-}x_5{\cdot}x^1_1 & x^1_1{-}x_5{\cdot}x^0_1 
& x^0_1 \cr
x^4_2{-}x_5{\cdot}x^3_2 & x^3_2{-}x_5{\cdot}x^2_2 & x^2_2{-}x_5{\cdot}x^1_2 
& x^1_2{-}x_5{\cdot}x^0_2 & x^0_2 \cr
x^4_3{-}x_5{\cdot}x^3_3 & x^3_3{-}x_5{\cdot}x^2_3 & x^2_3{-}x_5{\cdot}x^1_3 
& x^1_3{-}x_5{\cdot}x^0_3 & x^0_3 \cr
x^4_4{-}x_5{\cdot}x^3_4 & x^3_4{-}x_5{\cdot}x^2_4 & x^2_4{-}x_5{\cdot}x^1_4 
& x^1_4{-}x_5{\cdot}x^0_4 & x^0_4 \cr
x^4_5{-}x_5{\cdot}x^3_5 & x^3_5{-}x_5{\cdot}x^2_5 & x^2_5{-}x_5{\cdot}x^1_5 
& x^1_5{-}x_5{\cdot}x^0_5 & x^0_5\end{matrix} \right\| $\\
\\

 \qquad $=  {\det}\left\| \begin{matrix}(x_1{-}x_5){\cdot}x^3_1 
& (x_1{-}x_5){\cdot}x^2_1 & (x_1{-}x_5){\cdot}x^1_1 & (x_1{-}x_5){\cdot}x^0_1 
& x^0_1 \cr
(x_2{-}x_5){\cdot}x^3_2 & (x_2{-}x_5){\cdot}x^2_2 & (x_2{-}x_5){\cdot}x^1_2 
& (x_2{-}x_5){\cdot}x^0_2 & x^0_2 \cr
(x_3{-}x_5){\cdot}x^3_3 & (x_3{-}x_5){\cdot}x^2_3 & (x_3{-}x_5){\cdot}x^1_3 
& (x_3{-}x_5){\cdot}x^0_3 & x^0_3 \cr
(x_4{-}x_5){\cdot}x^3_4 & (x_4{-}x_5){\cdot}x^2_4 & (x_4{-}x_5){\cdot}x^1_4 
& (x_4{-}x_5){\cdot}x^0_4 & x^0_4 \cr
(x_5{-}x_5){\cdot}x^3_5 & (x_5{-}x_5){\cdot}x^2_5 & (x_5{-}x_5){\cdot}x^1_5 
& (x_5{-}x_5){\cdot}x^0_5 & x^0_5\end{matrix} \right\| $\\
\\

 \qquad $=  {\det}\left\| \begin{matrix}(x_1{-}x_5){\cdot}x^3_1 
& (x_1{-}x_5){\cdot}x^2_1 & (x_1{-}x_5){\cdot}x^1_1 & (x_1{-}x_5){\cdot}x^0_1 
& 1 \cr
(x_2{-}x_5){\cdot}x^3_2 & (x_2{-}x_5){\cdot}x^2_2 & (x_2{-}x_5){\cdot}x^1_2 
& (x_2{-}x_5){\cdot}x^0_2 & 1 \cr
(x_3{-}x_5){\cdot}x^3_3 & (x_3{-}x_5){\cdot}x^2_3 & (x_3{-}x_5){\cdot}x^1_3 
& (x_3{-}x_5){\cdot}x^0_3 & 1 \cr
(x_4{-}x_5){\cdot}x^3_4 & (x_4{-}x_5){\cdot}x^2_4 & (x_4{-}x_5){\cdot}x^1_4 
& (x_4{-}x_5){\cdot}x^0_4 & 1 \cr
             &              &              &              & 
1\end{matrix} \right\| $\\
\\

 \qquad $=  {\det}\left\| \begin{matrix}(x_1{-}x_5){\cdot}x^3_1 
& (x_1{-}x_5){\cdot}x^2_1 & (x_1{-}x_5){\cdot}x^1_1 & (x_1{-}x_5){\cdot}x^0_1 
\cr
(x_2{-}x_5){\cdot}x^3_2 & (x_2{-}x_5){\cdot}x^2_2 & (x_2{-}x_5){\cdot}x^1_2 
& (x_2{-}x_5){\cdot}x^0_2 \cr
(x_3{-}x_5){\cdot}x^3_3 & (x_3{-}x_5){\cdot}x^2_3 & (x_3{-}x_5){\cdot}x^1_3 
& (x_3{-}x_5){\cdot}x^0_3 \cr
(x_4{-}x_5){\cdot}x^3_4 & (x_4{-}x_5){\cdot}x^2_4 & (x_4{-}x_5){\cdot}x^1_4 
& (x_4{-}x_5){\cdot}x^0_4\end{matrix} \right\| $\\
\\

 \qquad $=  (x_1{-}x_5){\cdot}(x_2{-}x_5){\cdot}(x_3{-}x_5){\cdot}(x_4{-}x_5){\cdot}{\det}\left\| 
\begin{matrix}x^3_1 & x^2_1 & x^1_1 & x^0_1 \cr
x^3_2 & x^2_2 & x^1_2 & x^0_2 \cr
x^3_3 & x^2_3 & x^1_3 & x^0_3 \cr
x^3_4 & x^2_4 & x^1_4 & x^0_4\end{matrix} \right\| $\\
\\

 \qquad $=  ((x_1{-}x_5){\cdot}(x_2{-}x_5){\cdot}(x_3{-}x_5){\cdot}(x_4{-}x_5))$ 
 \smallskip

 \qquad \qquad ${\cdot}((x_1{-}x_4){\cdot}(x_2{-}x_4){\cdot}(x_3{-}x_4)){\cdot}((x_1{-}x_3){\cdot}(x_2{-}x_3)){\cdot}(x_1{-}x_2)$.
\hfil\vfill\eject
 
 \noindent {\bf Example 4.25.} (Induction step $E(3) \Rightarrow 
 E(4)$ in proof of 1 of lemma 3.1).\\

 ${\det}\left\| \begin{matrix}F(x_1) & x^3_1 & x^2_1 & x^1_1 
& x^0_1 \cr
F(x_2) & x^3_2 & x^2_2 & x^1_2 & x^0_2 \cr
F(x_3) & x^3_3 & x^2_3 & x^1_3 & x^0_3 \cr
F(x_4) & x^3_4 & x^2_4 & x^1_4 & x^0_4 \cr
F(x_5) & x^3_5 & x^2_5 & x^1_5 & x^0_5\end{matrix} \right\| 
 =  {\det}\left\| \begin{matrix}F(x_1) & x^3_1 & x^2_1 & x^1_1 
& x^0_1 \cr
F(x_2) & x^3_2 & x^2_2 & x^1_2 & x^0_2 \cr
F(x_3) & x^3_3 & x^2_3 & x^1_3 & x^0_3 \cr
F(x_4) & x^3_4 & x^2_4 & x^1_4 & x^0_4 \cr
F(x_5) & x^3_5 & x^2_5 & x^1_5 & x^0_5\end{matrix} \right\| 
{\cdot}\left\| \begin{matrix} 1\h &        &        &        
&   \cr
         &  1\h &        &        &   \cr
         & {-}x_5 &  1\h &        &   \cr
         &        & {-}x_5 &  1\h &   \cr
{-}F(x_5) &        &        & {-}x_5 & 1\end{matrix} \right\| 
$\\
\\

 \qquad $=  {\det}\left\| \begin{matrix}F(x_1){-}F(x_5){\cdot}x^0_1 
& x^3_1{-}x_5{\cdot}x^2_1 & x^2_1{-}x_5{\cdot}x^1_1 & x^1_1{-}x_5{\cdot}x^0_1 
& x^0_1 \cr
F(x_2){-}F(x_5){\cdot}x^0_2 & x^3_2{-}x_5{\cdot}x^2_2 & x^2_2{-}x_5{\cdot}x^1_2 
& x^1_2{-}x_5{\cdot}x^0_2 & x^0_2 \cr
F(x_3){-}F(x_5){\cdot}x^0_3 & x^3_3{-}x_5{\cdot}x^2_3 & x^2_3{-}x_5{\cdot}x^1_3 
& x^1_3{-}x_5{\cdot}x^0_3 & x^0_3 \cr
F(x_4){-}F(x_5){\cdot}x^0_4 & x^3_4{-}x_5{\cdot}x^2_4 & x^2_4{-}x_5{\cdot}x^1_4 
& x^1_4{-}x_5{\cdot}x^0_4 & x^0_4 \cr
F(x_5){-}F(x_5){\cdot}x^0_5 & x^3_5{-}x_5{\cdot}x^2_5 & x^2_5{-}x_5{\cdot}x^1_5 
& x^1_5{-}x_5{\cdot}x^0_5 & x^0_5\end{matrix} \right\| $\\
\\

 \qquad $=  {\det}\left\| \begin{matrix}F(x_1){-}F(x_5) & (x_1{-}x_5){\cdot}x^2_1 
& (x_1{-}x_5){\cdot}x^1_1 & (x_1{-}x_5){\cdot}x^0_1 & 1 \cr
F(x_2){-}F(x_5) & (x_2{-}x_5){\cdot}x^2_2 & (x_2{-}x_5){\cdot}x^1_2 
& (x_2{-}x_5){\cdot}x^0_2 & 1 \cr
F(x_3){-}F(x_5) & (x_3{-}x_5){\cdot}x^2_3 & (x_3{-}x_5){\cdot}x^1_3 
& (x_3{-}x_5){\cdot}x^0_3 & 1 \cr
F(x_4){-}F(x_5) & (x_4{-}x_5){\cdot}x^2_4 & (x_4{-}x_5){\cdot}x^1_4 
& (x_4{-}x_5){\cdot}x^0_4 & 1 \cr
              &              &              &              
& 1\end{matrix} \right\| $\\
\\

 \qquad $=  {\det}\left\| \begin{matrix}(x_1{-}x_5){\cdot}(\nabla 
(x_1,x_5;x)_{*}\, F(x)) & (x_1{-}x_5){\cdot}x^2_1 & (x_1{-}x_5){\cdot}x^1_1 
& (x_1{-}x_5){\cdot}x^0_1 & 1 \cr
(x_2{-}x_5){\cdot}(\nabla (x_2,x_5;x)_{*}\, F(x)) & (x_2{-}x_5){\cdot}x^2_2 
& (x_2{-}x_5){\cdot}x^1_2 & (x_2{-}x_5){\cdot}x^0_2 & 1 \cr
(x_3{-}x_5){\cdot}(\nabla (x_3,x_5;x)_{*}\, F(x)) & (x_3{-}x_5){\cdot}x^2_3 
& (x_3{-}x_5){\cdot}x^1_3 & (x_3{-}x_5){\cdot}x^0_3 & 1 \cr
(x_4{-}x_5){\cdot}(\nabla (x_4,x_5;x)_{*}\, F(x)) & (x_4{-}x_5){\cdot}x^2_4 
& (x_4{-}x_5){\cdot}x^1_4 & (x_4{-}x_5){\cdot}x^0_4 & 1 \cr
                               &              &            
  &              & 1\end{matrix} \right\| $\\
\\

 \qquad $=  {\det}\left\| \begin{matrix}(x_1{-}x_5){\cdot}(\nabla 
(x_1,x_5;x)_{*}\, F(x)) & (x_1{-}x_5){\cdot}x^2_1 & (x_1{-}x_5){\cdot}x^1_1 
& (x_1{-}x_5){\cdot}x^0_1 \cr
(x_2{-}x_5){\cdot}(\nabla (x_2,x_5;x)_{*}\, F(x)) & (x_2{-}x_5){\cdot}x^2_2 
& (x_2{-}x_5){\cdot}x^1_2 & (x_2{-}x_5){\cdot}x^0_2 \cr
(x_3{-}x_5){\cdot}(\nabla (x_3,x_5;x)_{*}\, F(x)) & (x_3{-}x_5){\cdot}x^2_3 
& (x_3{-}x_5){\cdot}x^1_3 & (x_3{-}x_5){\cdot}x^0_3 \cr
(x_4{-}x_5){\cdot}(\nabla (x_4,x_5;x)_{*}\, F(x)) & (x_4{-}x_5){\cdot}x^2_4 
& (x_4{-}x_5){\cdot}x^1_4 & (x_4{-}x_5){\cdot}x^0_4\end{matrix} 
\right\| $\\
\\

 \qquad $=  (x_1{-}x_5){\cdot}(x_2{-}x_5){\cdot}(x_3{-}x_5){\cdot}(x_4{-}x_5){\cdot}{\det}\left\| 
\begin{matrix}\nabla (x_1,x_5;x)_{*}\, F(x) & x^2_1 & x^1_1 & 
x^0_1 \cr
\nabla (x_2,x_5;x)_{*}\, F(x) & x^2_2 & x^1_2 & x^0_2 \cr
\nabla (x_3,x_5;x)_{*}\, F(x) & x^2_3 & x^1_3 & x^0_3 \cr
\nabla (x_4,x_5;x)_{*}\, F(x) & x^2_4 & x^1_4 & x^0_4\end{matrix} 
\right\| $\\
\\

 \qquad $=  (x_1{-}x_5){\cdot}(x_2{-}x_5){\cdot}(x_3{-}x_5){\cdot}(x_4{-}x_5)$
\\

 \qquad \qquad \qquad ${\cdot}{\det}\left\| \begin{matrix}x^3_1 
& x^2_1 & x^1_1 & x^0_1 \cr
x^3_2 & x^2_2 & x^1_2 & x^0_2 \cr
x^3_3 & x^2_3 & x^1_3 & x^0_3 \cr
x^3_4 & x^2_4 & x^1_4 & x^0_4\end{matrix} \right\| {\cdot}(\nabla 
(x_1,x_2,x_3,x_4;x)_{*}\, \nabla (x,x_5;x)_{*}\, F(x))$\\
\\

 \qquad $=  {\det}\left\| \begin{matrix}x^4_1 & x^3_1 & x^2_1 
& x^1_1 & x^0_1 \cr
x^4_2 & x^3_2 & x^2_2 & x^1_2 & x^0_2 \cr
x^4_3 & x^3_3 & x^2_3 & x^1_3 & x^0_3 \cr
x^4_4 & x^3_4 & x^2_4 & x^1_4 & x^0_4 \cr
x^4_5 & x^3_5 & x^2_5 & x^1_5 & x^0_5\end{matrix} \right\| 
{\cdot}(\nabla (x_1,x_2,x_3,x_4,x_5;x)_{*}\, F(x))$.\\
\hfil\vfill\eject
 
 \noindent {\bf Example 4.26.} (to corollary 3.1) Let \medskip

 $F(x) =  \sum\limits^{d{+}{\Delta} } _{{\delta} {=}0} F_{\delta} 
{\cdot}x^{\delta} $.  \medskip\\
Let $d{=}3$, $r{=}4$, ${\Delta} {=}r{-}1$. By of 2 of lemma 
3.2\\

 $(\nabla (x_1,x_2,x_3,x_4;x)_{*}\, F(x)){\cdot}{\det}\left\| 
\begin{matrix}x^3_1 & x^2_1 & x^1_1 & x^0_1 \cr
x^3_2 & x^2_2 & x^1_2 & x^0_2 \cr
x^3_3 & x^2_3 & x^1_3 & x^0_3 \cr
x^3_4 & x^2_4 & x^1_4 & x^0_4\end{matrix} \right\|  =  {\det}\left\| 
\begin{matrix}F(x_1) & x^3_1 & x^2_1 & x^1_1 & x^0_1 \cr
F(x_2) & x^3_2 & x^2_2 & x^1_2 & x^0_2 \cr
F(x_3) & x^3_3 & x^2_3 & x^1_3 & x^0_3 \cr
F(x_4) & x^3_4 & x^2_4 & x^1_4 & x^0_4 \cr
F(x_5) & x^3_5 & x^2_5 & x^1_5 & x^0_5\end{matrix} \right\| 
$.\\
\\
 that is equivalent to \bigskip

 $(\nabla ({\bf x}_1,{\bf x}_2,{\bf x}_3,{\bf x}_4;x)_{*}\, 
F(x)){\cdot}(x^3{\wedge}x^2{\wedge}x^1{\wedge}x^0) =  F(x){\wedge}x^2{\wedge}x^1{\wedge}x^0$, 
\bigskip\\
and by of 1 of theorem 2.2 \bigskip

 $\nabla (x_1,x_2,x_3,x_4;x)_{*}\, F(x)$  \bigskip

 \qquad $=  {\det}\left\| \begin{matrix}\overline {\sigma} 
_0 &   &   & F_6 &   &   &   \cr
\overline {\sigma} _1 & \overline {\sigma} _0 &   & F_5 &  
 &   &   \cr
\overline {\sigma} _2 & \overline {\sigma} _1 & \overline {\sigma} 
_0 & F_4 &   &   &   \cr
\overline {\sigma} _3 & \overline {\sigma} _2 & \overline {\sigma} 
_1 & F_3 &   &   &   \cr
\overline {\sigma} _4 & \overline {\sigma} _3 & \overline {\sigma} 
_2 & F_2 & 1 &   &   \cr
  & \overline {\sigma} _4 & \overline {\sigma} _3 & F_1 &  
 & 1 &   \cr
  &   & {\sigma} _4 & F_0 &   &   & 1\end{matrix} \right\| 
 =  {\det}\left\| \begin{matrix}\overline {\sigma} _0 &   &  
 & F_6 \cr
\overline {\sigma} _1 & \overline {\sigma} _0 &   & F_5 \cr
\overline {\sigma} _2 & \overline {\sigma} _1 & \overline {\sigma} 
_0 & F_4 \cr
\overline {\sigma} _3 & \overline {\sigma} _2 & \overline {\sigma} 
_1 & F_3\end{matrix} \right\| $.\\
\\
 
 \noindent Let $d{=}4$, $r{=}3$, then\\

 $\nabla (x_1,x_2,x_3;x)_{*}\, F(x) =  \frac{  {\det}\left\| 
\begin{matrix}F(x_1) & x^1_1 & x^0_1 \cr
F(x_2) & x^1_2 & x^0_2 \cr
F(x_3) & x^1_3 & x^0_3\end{matrix} \right\| } {  {\det}\left\| 
\begin{matrix}x^2_1 & x^1_1 & x^0_1 \cr
x^2_2 & x^1_2 & x^0_2 \cr
x^2_3 & x^1_3 & x^0_3\end{matrix} \right\| } $\\
\\

 \qquad $=  {\det}\left\| \begin{matrix}{\sigma} _0 &   & F_4 
&   &   \cr
{\sigma} _1 & {\sigma} _0 & F_3 &   &   \cr
{\sigma} _2 & {\sigma} _1 & F_2 &   &   \cr
{\sigma} _3 & {\sigma} _2 & F_1 & 1 &   \cr
  & {\sigma} _3 & F_0 &   & 1\end{matrix} \right\|  =  {\det}\left\| 
\begin{matrix}{\sigma} _0 &   & F_4 \cr
{\sigma} _1 & {\sigma} _0 & F_3 \cr
{\sigma} _2 & {\sigma} _1 & F_2\end{matrix} \right\| $\\
\\
\\
Let $d{=}2$, $r{=}3$, then\\
\\

 $\nabla (x_1,x_2,x_3;x)_{*}\, F(x) =  {\det}\left\| \begin{matrix}F_2 
&   &   \cr
F_1 & 1 &   \cr
F_0 &   & 1\end{matrix} \right\|  =  {\det}\left\| F_2\right\| 
$.\\
\\
\\
Let $d{=}3$, $r{=}3$, then\\
\\

 $\nabla (x_1,x_2,x_3;x)_{*}\, F(x) =  {\det}\left\| \begin{matrix}{\sigma} 
_0 & F_3 &   &   \cr
{\sigma} _1 & F_2 &   &   \cr
{\sigma} _2 & F_1 & 1 &   \cr
{\sigma} _3 & F_0 &   & 1\end{matrix} \right\|  =  {\det}\left\| 
\begin{matrix}{\sigma} _0 & F_3 \cr
{\sigma} _1 & F_2\end{matrix} \right\| $.
\hfil\vfill\eject
 
 \noindent Let $d{=}3$, $r{=}4$, then\\
\\

 $\nabla (x_1,x_2,x_3,x_4;x)_{*}\, F(x) =  {\det}\left\| \begin{matrix}F_3 
&   &   &   \cr
F_2 & 1 &   &   \cr
F_1 &   & 1 &   \cr
F_0 &   &   & 1\end{matrix} \right\|  =  {\det}\left\| F_3\right\| 
$.\\
\\
\\
Let $d{=}4$, $r{=}4$, then\\
\\

 $\nabla (x_1,x_2,x_3,x_4;x)_{*}\, F(x) =  {\det}\left\| \begin{matrix}{\sigma} 
_0 & F_4 &   &   &   \cr
{\sigma} _1 & F_3 &   &   &   \cr
{\sigma} _2 & F_2 & 1 &   &   \cr
{\sigma} _3 & F_1 &   & 1 &   \cr
{\sigma} _4 & F_0 &   &   & 1\end{matrix} \right\|  =  {\det}\left\| 
\begin{matrix}{\sigma} _0 & F_4 \cr
{\sigma} _1 & F_3\end{matrix} \right\| $.\\
\\
\\
Let $d{=}5$, $r{=}4$, then\\
\\

 $\nabla (x_1,x_2,x_3,x_4;x)_{*}\, F(x) =  {\det}\left\| \begin{matrix}{\sigma} 
_0 &   & F_5 &   &   &   \cr
{\sigma} _1 & {\sigma} _0 & F_4 &   &   &   \cr
{\sigma} _2 & {\sigma} _1 & F_3 &   &   &   \cr
{\sigma} _3 & {\sigma} _2 & F_2 & 1 &   &   \cr
{\sigma} _4 & {\sigma} _3 & F_1 &   & 1 &   \cr
  & {\sigma} _4 & F_0 &   &   & 1\end{matrix} \right\|  =  
{\det}\left\| \begin{matrix}{\sigma} _0 &   & F_5 \cr
{\sigma} _1 & {\sigma} _0 & F_4 \cr
{\sigma} _2 & {\sigma} _1 & F_3\end{matrix} \right\| $.\\
\\
\\
Let $d{=}6$, $r{=}4$, then\\
\\

 $\nabla (x_1,x_2,x_3,x_4;x)_{*}\, F(x) =  {\det}\left\| \begin{matrix}{\sigma} 
_0 &   &   & F_6 &   &   &   \cr
{\sigma} _1 & {\sigma} _0 &   & F_5 &   &   &   \cr
{\sigma} _2 & {\sigma} _1 & {\sigma} _0 & F_4 &   &   &   \cr
{\sigma} _3 & {\sigma} _2 & {\sigma} _1 & F_3 &   &   &   \cr
{\sigma} _4 & {\sigma} _3 & {\sigma} _2 & F_2 & 1 &   &   \cr
  & {\sigma} _4 & {\sigma} _3 & F_1 &   & 1 &   \cr
  &   & {\sigma} _4 & F_0 &   &   & 1\end{matrix} \right\| 
 =  {\det}\left\| \begin{matrix}{\sigma} _0 &   &   & F_6 \cr
{\sigma} _1 & {\sigma} _0 &   & F_5 \cr
{\sigma} _2 & {\sigma} _1 & {\sigma} _0 & F_4 \cr
{\sigma} _3 & {\sigma} _2 & {\sigma} _1 & F_3\end{matrix} \right\| $.
\hfil\vfill\eject
 \noindent Let $d{=}1$, $r{=}2$, then\\
\\

 $\nabla (x_1,x_2;x)_{*}\, F(x) =  {\det}\left\| \begin{matrix}F_1 
&   \cr
F_0 & 1\end{matrix} \right\|  =  {\det}\left\| F_1\right\| 
$.\\
\\
\\
Let $d{=}2$, $r{=}2$, then\\
\\

 $\nabla (x_1,x_2;x)_{*}\, F(x) =  {\det}\left\| \begin{matrix}{\sigma} 
_0 & F_2 &   \cr
{\sigma} _1 & F_1 &   \cr
{\sigma} _2 & F_0 & 1\end{matrix} \right\|  =  {\det}\left\| 
\begin{matrix}{\sigma} _0 & F_2 \cr
{\sigma} _1 & F_1\end{matrix} \right\| $.\\
\\
\\
Let $d{=}3$, $r{=}2$, then\\
\\

 $\nabla (x_1,x_2;x)_{*}\, F(x) =  {\det}\left\| \begin{matrix}{\sigma} 
_0 &   & F_3 &   \cr
{\sigma} _1 & {\sigma} _0 & F_2 &   \cr
{\sigma} _2 & {\sigma} _1 & F_1 &   \cr
  & {\sigma} _2 & F_0 & 1\end{matrix} \right\|  =  {\det}\left\| 
\begin{matrix}{\sigma} _0 &   & F_3 \cr
{\sigma} _1 & {\sigma} _0 & F_2 \cr
{\sigma} _2 & {\sigma} _1 & F_1\end{matrix} \right\| $.\\
\\
\\
Let $d{=}4$, $r{=}2$, then\\
\\

 $\nabla (x_1,x_2;x)_{*}\, F(x) =  {\det}\left\| \begin{matrix}{\sigma} 
_0 &   &   & F_4 &   \cr
{\sigma} _1 & {\sigma} _0 &   & F_3 &   \cr
{\sigma} _2 & {\sigma} _1 & {\sigma} _0 & F_2 &   \cr
  & {\sigma} _2 & {\sigma} _1 & F_1 &   \cr
  &   & {\sigma} _2 & F_0 & 1\end{matrix} \right\|  =  {\det}\left\| 
\begin{matrix}{\sigma} _0 &   &   & F_4 \cr
{\sigma} _1 & {\sigma} _0 &   & F_3 \cr
{\sigma} _2 & {\sigma} _1 & {\sigma} _0 & F_2 \cr
  & {\sigma} _2 & {\sigma} _1 & F_1\end{matrix} \right\| $.
\\
\\
\hfil\vfill\eject
 
 \noindent Let $d{=}0$, $r{=}1$, then\\

 $\nabla (x_1;x)_{*}\, F(x) =  {\det}\left\| F_0\right\|  = 
 {\det}\left\| F_0\right\|  =  F(x_1)$.\\
\\
\\
Let $d{=}1$, $r{=}1$, then\\
\\

 $\nabla (x_1;x)_{*}\, F(x) =  {\det}\left\| \begin{matrix}{\sigma} 
_0 & F_1 \cr
{\sigma} _1 & F_0\end{matrix} \right\|  =  {\det}\left\| \begin{matrix}{\sigma} 
_0 & F_1 \cr
{\sigma} _1 & F_0\end{matrix} \right\|  =  {\det}\left\| \begin{matrix}1\h 
& F_1 \cr
{-}x_1 & F_0\end{matrix} \right\|  =  F(x_1)$.\\
\\
\\
Let $d{=}2$, $r{=}1$, then\\
\\

 $\nabla (x_1;x)_{*}\, F(x) =  {\det}\left\| \begin{matrix}{\sigma} 
_0 &   & F_2 \cr
{\sigma} _1 & {\sigma} _0 & F_1 \cr
  & {\sigma} _1 & F_0\end{matrix} \right\|  =  {\det}\left\| 
\begin{matrix}{\sigma} _0 &   & F_2 \cr
{\sigma} _1 & {\sigma} _0 & F_1 \cr
  & {\sigma} _1 & F_0\end{matrix} \right\|  =  {\det}\left\| 
\begin{matrix}1\h &   & F_2 \cr
{-}x_1 & 1\h & F_1 \cr
  & {-}x_1 & F_0\end{matrix} \right\|  =  F(x_1)$.\\
\\
\\
Let $d{=}3$, $r{=}1$, then\\

 $\nabla (x_1;x)_{*}\, F(x) =  {\det}\left\| \begin{matrix}{\sigma} 
_0 &   &   & F_3 \cr
{\sigma} _1 & {\sigma} _0 &   & F_2 \cr
  & {\sigma} _1 & {\sigma} _0 & F_1 \cr
  &   & {\sigma} _1 & F_0\end{matrix} \right\|  =  {\det}\left\| 
\begin{matrix}{\sigma} _0 &   &   & F_3 \cr
{\sigma} _1 & {\sigma} _0 &   & F_2 \cr
  & {\sigma} _1 & {\sigma} _0 & F_1 \cr
  &   & {\sigma} _1 & F_0\end{matrix} \right\| $\\
\\

 \qquad \qquad \qquad \qquad \qquad  $=  {\det}\left\| \begin{matrix}1\h 
&   &   & F_3 \cr
{-}x_1 & 1\h &   & F_2 \cr
  & {-}x_1 & 1\h & F_1 \cr
  &   & {-}x_1 & F_0\end{matrix} \right\|  =  F(x_1)$.\\
\\
\\
Let $d{=}4$, $r{=}1$, then\\
\\

 $\nabla (x_1;x)_{*}\, F(x) =  {\det}\left\| \begin{matrix}{\sigma} 
_0 &   &   &   & F_4 \cr
{\sigma} _1 & {\sigma} _0 &   &   & F_3 \cr
  & {\sigma} _1 & {\sigma} _0 &   & F_2 \cr
  &   & {\sigma} _1 & {\sigma} _0 & F_1 \cr
  &   &   & {\sigma} _1 & F_0\end{matrix} \right\|  =  {\det}\left\| 
\begin{matrix}{\sigma} _0 &   &   &   & F_4 \cr
{\sigma} _1 & {\sigma} _0 &   &   & F_3 \cr
  & {\sigma} _1 & {\sigma} _0 &   & F_2 \cr
  &   & {\sigma} _1 & {\sigma} _0 & F_1 \cr
  &   &   & {\sigma} _1 & F_0\end{matrix} \right\| $\\
\\

 \qquad \qquad \qquad \qquad \qquad $=  {\det}\left\| \begin{matrix}1\h 
&   &   &   & F_4 \cr
{-}x_1 & 1\h &   &   & F_3 \cr
  & {-}x_1 & 1\h &   & F_2 \cr
  &   & {-}x_1 & 1\h & F_1 \cr
  &   &   & {-}x_1 & F_0\end{matrix} \right\|  =  F(x_1)$.
\\

\hfil\vfill\eject

 \noindent {\bf Example 4.27.} (to corollary 3.2) Let $d{=}6$, $r{=}4$, \medskip

 $F(x) =  \sum\limits^{ d} _{{\delta} {=}0} F_{\delta} {\cdot}x^{\delta} 
$, \medskip

 $f(x) =  (x{-}{\lambda} _1){\cdot}(x{-}{\lambda} _2){\cdot}(x{-}{\lambda} _3){\cdot}(x{-}{\lambda} _4)$   \medskip

 \qquad $=  f_4{\cdot}x^4{+}f_3{\cdot}x^3{+}f_2{\cdot}x^2{+}f_1{\cdot}x^1{+}f_0{\cdot}x^0$. 
\medskip
\\
Tnen\\
\\

 $\nabla ({\lambda} _1,{\lambda} _2,{\lambda} _3,{\lambda} 
_4;x)_{*}\, F(x) =  {\det}\left\| \begin{matrix}f_4 &   &   & 
F_6 &   &   &   \cr
f_3 & f_4 &   & F_5 &   &   &   \cr
f_2 & f_3 & f_4 & F_4 &   &   &   \cr
f_2 & f_2 & f_3 & F_3 &   &   &   \cr
f_0 & f_2 & f_2 & F_2 & 1 &   &   \cr
  & f_0 & f_2 & F_1 &   & 1 &   \cr
  &   & f_0 & F_0 &   &   & 1\end{matrix} \right\|  =  {\det}\left\| 
\begin{matrix}f_4 &   &   & F_6 \cr
f_3 & f_4 &   & F_5 \cr
f_2 & f_3 & f_4 & F_4 \cr
f_2 & f_2 & f_3 & F_3\end{matrix} \right\| $.\\
\\

 $\nabla ({\lambda} _1,{\lambda} _2,{\lambda} _3,{\lambda} 
_4;x)_{*} \equiv  {\det}\left\| \begin{matrix}f_4 &   &   & [x^6]_{*} 
&   &   &   \cr
f_3 & f_4 &   & [x^5]_{*} &   &   &   \cr
f_2 & f_3 & f_4 & [x^4]_{*} &   &   &   \cr
f_2 & f_2 & f_3 & [x^3]_{*} &   &   &   \cr
f_0 & f_2 & f_2 & [x^2]_{*} & 1 &   &   \cr
  & f_0 & f_2 & [x^1]_{*} &   & 1 &   \cr
  &   & f_0 & [x^0]_{*} &   &   & 1\end{matrix} \right\|  \equiv 
 {\det}\left\| \begin{matrix}f_4 &   &   & [x^6]_{*} \cr
f_3 & f_4 &   & [x^5]_{*} \cr
f_2 & f_3 & f_4 & [x^4]_{*} \cr
f_2 & f_2 & f_3 & [x^3]_{*}\end{matrix} \right\| $ in ${\bf 
R}[x]^{{\leq}d} $.\\
\\
\\


\footnotesize

\bibliographystyle{amsplain}

\end{document}